\documentclass{elsart}
\usepackage[english]{babel}
\usepackage{latexsym,amsfonts,amssymb,amsmath,longtable,mathrsfs}
\usepackage[unicode,final,hyperindex,hypertex]{hyperref}%pdf
\usepackage{bbm}
\usepackage{lineno}
\usepackage{cite}
\usepackage[varg]{txfonts}
\usepackage{stmaryrd}
\renewcommand{\le}{\leqslant}
\renewcommand{\ge}{\geqslant}
\renewcommand{\leq}{\leqslant}
\renewcommand{\geq}{\geqslant}
\renewcommand{\unlhd}{\trianglelefteqslant}

\newtheorem{Lemma}{{\bfseries Lemma}}[section]
\newtheorem{Cor}[Lemma]{{\bfseries Corollary}}
\newtheorem{Theo}[Lemma]{{\bfseries Theorem}}
\theoremstyle{definition}

\newtheorem{ex}[Lemma]{{\bfseries Example}}

\newcommand{\Sym}{{\mathrm{Sym}}}
\newcommand{\Alt}{{\mathrm{Alt}}}

\newcommand{\Aut}{{\mathrm{Aut}}}
\newcommand{\Out}{{\mathrm{Out}}}
\newcommand{\Inn}{{\mathrm{Inn}}}

\newcommand{\F}{\mathbbmss{F}}

\newcommand{\GL}{{\mathrm{GL}}}
\newcommand{\GU}{{\mathrm{GU}}}
\newcommand{\GO}{{\mathrm{GO}}}
\newcommand{\SO}{{\mathrm{SO}}}
\newcommand{\SL}{{\mathrm{SL}}}
\newcommand{\PSL}{{\mathrm{PSL}}}
\newcommand{\PSp}{{\mathrm{PSp}}}
\newcommand{\PGL}{{\mathrm{PGL}}}
\newcommand{\SU}{{\mathrm{SU}}}
\newcommand{\PSU}{{\mathrm{PSU}}}
\newcommand{\Sp}{{\mathrm{Sp}}}

\newcommand{\ov}{\overline}
\renewcommand{\P}{\mathrm{P}}

\newcommand{\Hall}{\mathrm{Hall}}
\newcommand{\splitext}{\,\colon\!}
\newcommand{\arbitraryext}{\,\ldotp}
\newcommand{\nonsplitext}{\!{}^{\text{\normalsize{\textperiodcentered}}}}
\newcommand{\Oo}{\mathop{\rm O}\nolimits}

\newcommand{\ti}{{\mbox{}_3}}
\newcommand{\tw}{{\mbox{}_2}}
\renewcommand{\emptyset}{\varnothing}

\begin{document}
%\linenumbers
%\sffamily
%\changenotsign
\begin{frontmatter}
\title{On the number of classes of conjugate Hall subgroups in finite simple groups\thanksref{general}}
\thanks[general]{The work is supported by  RFBR, projects 08-01-00322, 10-01-00391, and 10-01-90007,  ADTP ``Development
of the Scientific Potential of
Higher
School'' of the Russian Federal Agency for Education (Grant
2.1.1.419), and Federal Target Grant "Scientific and
educational personnel of innovation Russia" for 2009-2013 (government
contract No. 02.740.11.5191 and No. 14.740.11.0346).}
\author{D.O.Revin}
\ead{revin@math.nsc.ru}
\author{E.P.Vdovin\corauthref{Coraut}\thanksref{vdovin}}
\ead{vdovin@math.nsc.ru}
\thanks[vdovin]{The second author gratefully acknowledges the support from
Deligne 2004 Balzan prize in mathematics.}
\corauth[Coraut]{Corresponding author}
\address{Sobolev Institute of mathematics, Novosibirsk}
\begin{abstract}
In this paper we find the number of classes of conjugate $\pi$-Hall subgroups in all finite almost simple groups. We also complete the classification
of $\pi$-Hall subgroups in finite simple groups and correct some mistakes from our previous paper.
\end{abstract}

\end{frontmatter}

\section{Introduction}

Let~$\pi$ be a set of primes. We denote by~$\pi'$ the set of all primes not in
$\pi$, by~$\pi(n)$ the set of all prime divisors of a positive integer~$n$,
for a finite group~$G$ we denote~$\pi(|G|)$ by~$\pi(G)$. A positive integer
$n$ with~$\pi(n)\subseteq\pi$ is called a~$\pi$-num\-ber,  a group~$G$ with
$\pi(G)\subseteq \pi$ is called a~$\pi$-group. Given positive integer~$n$ denote by~$n_\pi$ the maximal divisor~$t$ of~$n$ with
$\pi(t)\subseteq\pi$. A subgroup~$H$ of~$G$ is called
a {\it~$\pi$-Hall subgroup}, if~$\pi(H)\subseteq\pi$ and~$\pi(|G:H|)\subseteq
\pi'$.  According to \cite{HallProperties} we  say that $G$ {\it satisfies
$E_\pi$} (or briefly $G\in E_\pi$), if $G$ possesses a $\pi$-Hall subgroup. If
$G\in E_\pi$ and every two $\pi$-Hall subgroups are conjugate, then we say
that  $G$ {\it satisfies $C_\pi$} ($G\in C_\pi$). If $G\in C_\pi$ and
each $\pi$-subgroup of $G$ is included in a  $\pi$-Hall
subgroup of $G$, then we  say that $G$ {\it satisfies $D_\pi$} ($G\in
D_\pi$). Thus $G\in D_\pi$ means that a complete analogue of the Sylow theorem for $\pi$-Hall subgroups of $G$ holds. A group satisfying
$E_\pi$ ($C_\pi$, $D_\pi$) is also called
an $E_\pi$-gro\-up
(respectively  $C_\pi$-group,  $D_\pi$-group).

In \cite[Theorem~7.7]{RevVdoContemp} the authors proved that a finite group $G$ satisfies $D_\pi$ if and only if each composition factor of
$G$ satisfies $D_\pi$. In the next series of papers \cite{OddDpi,DAN,Linuni,Two} for every finite simple group $S$ and for every set
$\pi$ of primes pure arithmetic necessary and sufficient condition for $S$ to satisfy  $D_\pi$ were found.

The authors intend to write a series of papers, where arithmetic criteria for $E_\pi$ and $C_\pi$ will be obtained. The present paper
is the first one in this series. Since, in contrast with $D_\pi$-groups, the class of $E_\pi$-group is not closed under extensions, while
the class of $C_\pi$-groups is not closed under normal subgroups, the general theory for $E_\pi$ and $C_\pi$
is more
complicated than the theory for $D_\pi$.  In particular, the answer to the question, whether given group $G$
satisfies $E_\pi$ or
$C_\pi$, cannot be obtained in terms of composition factors of $G$, i.e., this question cannot be reduced to  similar
questions about
simple groups.  We intend to reduce the question to similar questions for almost simple groups. Recall that a finite
group $G$ is called
{\em almost simple}, if its generalized Fitting subgroup $F^\ast(G)$ is a nonabelian simple group $S$, or, equivalently,
if $\Inn(S)\leq
G\leq \Aut(S)$ for a nonabelian finite simple group $S$. An important step in this direction was made by F.Gross
in~\cite{Gr0}.

In this paper we prove (by using the classification of finite simple groups) an important theorem on the number of
classes of conjugate
$\pi$-Hall subgroups in finite simple groups. We show that this number is bounded and is a $\pi$-number. This result (in
a more general form,
which will be used in future research) can be found in Theorem \ref{IndClassNumb} below. First we need to introduce some
notations.

We denote by $\mathrm{Hall}_\pi(G)$ the set of all $\pi$-Hall subgroups
of  $G$ (notice that this set can be
empty). Assume that $A\unlhd G$. Then define $\Hall_\pi^G(A)=\{H\cap A\mid
H\in\Hall_\pi(G)\}.$ By Lemma  \ref{simpleepi}(a) (see  below), $\Hall_\pi^G(A)\subseteq\mathrm{Hall}_\pi(A)$. The
elements
of $\Hall_\pi^G(A)$ are called  $G$-induced $\pi$-Hall subgroups of $A$. Clearly $A$ acts by conjugation on both
$\mathrm{Hall}_\pi(A)$ and $\Hall_\pi^G(A)$. Denote by  $k_\pi(A)$ and $k_\pi^G(A)$ the number of orbits under this
action, respectively.
Thus $k_\pi(A)$ is the number of classes of conjugate $\pi$-Hall subgroups of $A$,  $k_\pi^G(A)$ is the number of
classes of conjugate
$G$-induced $\pi$-Hall subgroups of $A$, and~${k_\pi^G(A)\leqslant k_\pi(A)}$.

\begin{Theo}\label{IndClassNumb}
Let $\pi$ be a set of primes, $G$ a finite almost simple group with nonabelian simple socle $S$. Then the following
statements hold{\em:}
\begin{itemize}
 \item[{\em (a)}] if $2\not\in\pi$, then $k_\pi^G(S)\in\{0,1\}${\em;}
 \item[{\em (b)}] if $3\not\in\pi$, then $k_\pi^G(S)\in\{0,1,2\}${\em;}
 \item[{\em (c)}] if $2,3\in\pi$, then $k_\pi^G(S)\in\{0,1,2,3,4,9\}$.
\end{itemize}
In particular, $k_\pi^G(S)$ is bounded and, if $G\in E_\pi$, then $k_\pi^G(S)$ is a $\pi$-number.
\end{Theo}

Setting $G=S$ we obtain that $k_\pi^G(S)=k_\pi(S)$, so the same statement on the number of classes of conjugate
$\pi$-Hall subgroups is
true for every simple group~$S$.

Theorem \ref{IndClassNumb}  generalizes a result by F.Gross \cite[Theorem~B]{GrossConjugacy}, which states that for
every set $\pi$ of odd
primes every finite simple $E_\pi$-group satisfies $C_\pi$ (equivalently $k_\pi(S)=1$). Since by Chunikhin's theorem the
class of
$C_\pi$-groups is closed under extension (see Lemma \ref{simpleepi}(f) below), it follows that $k_\pi(G)=1$ for every
finite group
$G$, see~\cite[Theorem~A]{GrossConjugacy}.

In contrast with the case $2\not\in\pi$, Example \ref{ExampleClassNumbIsNotPiNumb} shows that in case $2\in\pi$
Theorem~\ref{IndClassNumb}
cannot be generalized to arbitrary (nonsimple) groups.

\begin{ex}\label{ExampleClassNumbIsNotPiNumb}
{\em Assume that $X\in E_\pi$ is such that $k=k_\pi(X)>1$ (in particular, $\pi$ is not equal to the set of all primes).
Suppose $p\in\pi'$.
Denote a cyclic
subgroup of order $p$ of  $\Sym_p$ by $Y$.  Consider $G=X\wr Y$ and let $$M\simeq \underbrace{X\times\dots\times
X}_{\displaystyle p
\mbox{ \rm times}}$$ be the base of the wreath product. It is clear that $k_\pi(M)=k^p$. Since $M$ is a normal subgroup
of $G$ and
$|G:M|=p$ is a $\pi'$-number, we have ${\rm Hall}_\pi(G)={\rm Hall}_\pi(M)$. The subgroup $Y$ acts on the set of classes
of conjugate
$\pi$-Hall subgroups of $M$. Applying the Burnside formula to this action it is easy to show that
$$k_\pi(G)=\displaystyle\frac{k^p+(p-1)k}{p}.$$ Now assume that  $\pi=\{2,3\}$ and $X=\SL_3(2)$. Then
\cite[Theorem~1.2]{RevinSIBAM}
implies that $k_\pi(X)=2$. Since $p\in\pi'$ can be taken arbitrarily large and
$(2^p-2)/p+2$ tends to infinity as $p$ tends to infinity, we obtain that for a nonsimple group  $G$ the number
$k_\pi(G)$ is not bounded in
general. Moreover, if we take  $p=7$, then $k_\pi(G)=20$, whence it is possible that $k_\pi(G)$ is not a $\pi$-number. }
\end{ex}

Although Theorem \ref{IndClassNumb} is not true for arbitrary groups, it can be used in order to obtain important
results for finite
groups. As an example of this using we give a short solution to Problem  13.33 from ``Kourovka notebook'' \cite{Kour}.
Earlier this problem
was solved by the authors in \cite{rev2,RevVdoContemp} by using other arguments.

\begin{Cor}\label{DpiNorm}
Let $\pi$ be a set of primes,  $A$ a normal subgroup of a finite $D_\pi$-group $G$. Then~${A\in D_\pi}$.
\end{Cor}

Theorem \ref{IndClassNumb} follows from the classification of Hall subgroups in finite simple groups. We briefly outline
the history.

Hall subgroups in  finite groups  close to  simple
were investigated by many authors. In
\cite{HallProperties} P.Hall found solvable $\pi$-Hall subgroups in symmetric
groups. Together with the famous Odd Order Theorem (see \cite{FeitThompsonOddOrder})
this result implies the classification of Hall subgroups of odd order in
alternating groups. Later J.Thompson in \cite{ThomHallSymmetric} found
non-solvable $\pi$-Hall subgroups in symmetric groups. The problem of
classification of Hall subgroups of even order in alternating groups remained
open for quite a long period. This classification was completed in
\cite{RevVdoContemp}. The classification of Hall subgroups in the alternating and
symmetric groups is given in Lemma \ref{HallSymmetric} below.

The Hall subgroups of odd order in sporadic simple groups are classified by F.Gross
in \cite{GrossConjectureHall}. The classification of Hall
subgroups in the sporadic groups was completed by the first author in~\cite{rev2}.

The classification of $\pi$-Hall subgroups in groups of Lie type in characteristic $p$ with $p\in
\pi$ was obtained by F.Gross (in case $2\not\in\pi$, \cite{GrossConjectureHall} and \cite{GrossNotDivisibleBy3}) and by
the first author (in case
$2\in\pi$, \cite{RevinSIBAM}). The classification of $\pi$-Hall subgroups in groups of Lie type
with $2,p\not\in\pi$ was obtained by F.Gross in classical groups \cite{GrossOddOrder}, and by
the authors in exceptional groups \cite{vdorev}. The classification
of $\pi$-Hall subgroups with $2\in\pi$ and $3,p\not\in\pi$ was obtained by F.Gross
in linear and symplectic groups \cite{GrossNotDivisibleBy3} and by the authors in the remaining cases
\cite{RevVdoContemp}.

In \cite{RevVdoContemp} the authors announced the classification of $\pi$-Hall
subgroups in finite simple groups. Unfortunately, Lemma 3.14 in
\cite{RevVdoContemp} contains a wrong statement. Namely, it states in item (a)
that if $G\simeq \PSL_2(q)$ and $\Alt_4\simeq H$ is a $\{2,3\}$-Hall subgroup
of $G$, then $\PGL_2(q)$ does not possesses a $\{2,3\}$-Hall subgroup
$H_1$ such that $H_1\cap G=H$, but this is not true. Due to this gap we lost
several series of $\pi$-Hall subgroups in groups of Lie type with $2,3\in\pi$. In the
present paper we correct this mistake and complete the classification of
$\pi$-Hall subgroups in finite simple groups (see Lemmas \ref{SLU23dim2}, \ref{HallSubgroupsOfLinearAndUnitaryGroups},
\ref{HallSubgroupsOfSymplecticGroups}, \ref{HallSubgroupsOfOrthogonalGroupsOfEvenDimension},
\ref{HallG2}--\ref{Hall3D4}). We correct also
other known minor
mistakes in proofs and statements from \cite{RevVdoContemp} and make some proofs
more elementary. The description of $\pi$-Hall subgroups in the groups of Lie type over a field of characteristic $p$ in
case $2,3\in\pi$ and
$p\not\in\pi$ takes a significant part of the paper. We also prove that in this case the groups of Lie type do not
satisfy $D_\pi$. This fact
means that all results from \cite{RevVdoContemp} concerning  $D_\pi$ (in particular, the above-mentioned
\cite[Theorem~7.7]{RevVdoContemp}),
and the results from \cite{DAN,OddDpi,Linuni,Two} remain valid and their proofs need no corrections.

\section{Notation and preliminary results}

Our notation for finite groups agrees with that of \cite{ATLAS}. For groups $A$
and $B$ symbols $A\times B$ and $A\circ B$ denote direct and central products,
respectively. Recall that $A\circ B$ is a group possessing normal subgroups $A_1$ and $B_1$ isomorphic to $A$ and $B$
respectively such 
that $G=\langle A_1,B_1\rangle$ and $[A_1,B_1]=1$. By $A\splitext B$, $A\nonsplitext
B$, and
$A\arbitraryext B$ we denote a split, a nonsplit,
and an arbitrary extension of a group $A$ by a group $B$. For a group $G$ and a
subgroup $S$ of $\Sym_n$ by $G\wr S$ we always denote the permutation wreath
product. If $G$ is a finite group, then $O_\pi(G)$ denotes the
$\pi$-radical of $G$, i.e., the largest normal $\pi$-subgroup of $G$, while $O^{\pi'}(G)$ denotes the subgroup of $G$
generated by
all $\pi$-element. We write $m\le n$ if a real number $m$ is not greater than $n$, while we use notations $H\leq G$ and
$H\unlhd G$ instead
of ``$H$ is a subgroup of $G$'' and ``$H$ is a normal subgroup of $G$'', respectively. If $r$ is an odd prime and $q$ is
a positive integer
coprime to $r$, then $e(q,r)$  denotes
a~multi\-plicative order of $q$ modulo $r$, that is a minimal natural number $m$ with
$q^m\equiv1\pmod{r}$. If $q$ is an odd positive integer, then 
$$
e(q,2)=\left\{ 
\begin{array}{rl}
1&\text{  if }q\equiv1\pmod 4 ,\\
2&\text{  if }q\equiv-1\pmod 4 .\\
\end{array}
\right.
$$
For
$M\subseteq G$ we set $M^G=\{M^g\mid g\in G\}$.
For a group $G$ we denote by $\Aut(G)$, $\Inn(G)$, and $\Out(G)$ the groups of all, inner, and outer automorphisms of
$G$, respectively.

A finite group $G$ is called {\em $\pi$-separable}, if $G$
possesses a subnormal series such that each factor of the series is either a $\pi$- or a $\pi'$- group. It is clear that
every
$\pi$-separable group possesses a normal series $1=G_0\lhd G_1\lhd \ldots\lhd G_k=G$ such that each factor of the series
is either a $\pi$-
or a $\pi'$- group. Clearly we may assume that $G_i$ is invariant under $\Aut(G)$ for all~$i$.

In Lemma \ref{simpleepi} we collect some known facts about Hall subgroups in finite
groups. Most of the results mentioned in Lemma \ref{simpleepi} are known and we
just give hints of their proofs.

\begin{Lemma}\label{simpleepi}
Let $G$ be a finite group, $A$ a normal subgroup of $G$.

\begin{itemize}
\item[{\em (a)}]For $H\in \Hall_\pi(G)$, we have
$H\cap A\in\Hall_\pi(A)$ and $HA/A\in\Hall_\pi(G/A)$, in particular~${\Hall_\pi^G(A)\subseteq \Hall_\pi(A)}$.

\item[{\em (b)}]  If  $M/A$ is a $\pi$-subgroup of $G/A$, then  there exists a $\pi$-subgroup $H$ of $G$
with~${M/A=HA/A}$.

\item[{\em (c)}] If $G\in E_\pi$ (resp. $G\in D_\pi$), then  $G/A\in E_\pi$ (resp. $G/A\in D_\pi$).

\item[{\em (d)}] If $A$ is $\pi$-separable and $G/A\in E_\pi$ (resp. $G/A\in C_\pi$, $G/A\in D_\pi$), then $G\in E_\pi$
(resp. $G\in
C_\pi$, $G\in D_\pi$).

\item[{\em (e)}] Assume that $A\in E_\pi$ and $\pi(G/A)\subseteq \pi$. Then a $\pi$-Hall subgroup $H$ of $A$ lies in
$\Hall_\pi^G(A)$
if and
only if $H^A=H^G$.

\item[{\em (f)}] Assume that $A$ satisfies $C_\pi$, $G/A$ satisfies $E_\pi$, and $M\leq G$ is chosen so
that $A\leq M$ and $M/A\in\Hall_\pi(G/A)$. Then $\emptyset\not=\Hall_\pi(M)\subseteq\Hall_\pi(G)$, and every $H,K\in
\Hall_\pi(G)$
are conjugate in $G$ if and only if $HA/A$ and $KA/A$ are conjugate in $G/A$.  In particular, $G\in E_\pi$ and if
$G/A\in C_\pi$,
then~${G\in C_\pi}$ {\em (}this statement generalizes Chunikhin's results {\em \cite[Theorems~C1,C2]{HallProperties})}.
\end{itemize}
\end{Lemma}

\noindent{\slshape Proof.}\ \
\noindent {\slshape Proof.}\ \
(a) See \cite[Lemma~1]{HallProperties}.

(b)  Let $H$ be a minimal subgroup of $G$ such that $HA/A=M/A$. Then $H\cap A$ is contained
in the Frattini subgroup of $H$. Indeed, if there exists a maximal subgroup $L$ of $H$, not containing
$H\cap A$, then clearly $LA=M$, which contradicts the minimality of
$H$. Thus the group $H/(H\cap A)\simeq HA/A$ is a $\pi$-group, while  $H\cap A$ is a nilpotent normal subgroup of $H$.
Therefore a $\pi'$-Hall subgroup $K$ of $H\cap A$ is a normal $\pi'$-Hall subgroup of $H$. In view of the
Schur-Zassenhaus theorem \cite[Chapter
IV, Satz~27]{Zassenhaus} (see also \cite[Theorems~D6,~D7]{HallProperties}) we obtain that $H$ possesses a $\pi$-Hall
subgroup $H_1$. Clearly
$H_1A=HA$, hence $H_1=H$ and $H$ is a $\pi$-group.

(c) Follows from (a) and (b).

(d) Follows from the Schur-Zassenhaus theorem \cite[Chapter
IV, Satz~27]{Zassenhaus} (see also \cite[Theorems~D6,~D7]{HallProperties}) and induction on the order of $A$.

(e) The ``only if'' part is evident, since
$G=H_1A$ for a $\pi$-Hall subgroup $H_1$ of $G$, containing $H$. Now we prove the ``if'' part. Since $G$ leaves
the set $\{H^a\mid a\in A\}$ invariant, Frattini argument implies that
$G=AN_G(H)$. Now $N_G(H)$ possesses a normal series $1\lhd H\lhd N_A(H)\lhd N_G(H)$ and all
section in this series are either $\pi$- or $\pi'$- groups, so $N_G(H)$ is $\pi$-separable. Statement (d) of
the lemma implies that $N_G(H)\in D_\pi$. In particular, there exists a
$\pi$-Hall subgroup $H_1$ of $N_G(H)$ and $H\leq H_1$. From
$\pi(G/A)\subseteq\pi$, $\pi(\vert A:H\vert)\subseteq\pi'$, and $G=AN_G(H)$
we derive that $\vert G:N_G(A)\vert$ is a $\pi'$-number. Hence $H_1$ is a
$\pi$-Hall subgroup of~$G$.

(f) Assume that  $M/A\in \Hall_\pi(G)$ and $M$ is the complete preimage of $M/A$ in $G$ under the natural homomorphism
$G\rightarrow G/A$.
Point (e) of the lemma implies that $M\in E_\pi$, in particular there exists $L\in \Hall_\pi(M)$. Since $\vert
G:L\vert=\vert
G:M\vert\cdot \vert M:L\vert$ is a $\pi'$-num\-ber, we obtain that $L\in \Hall_\pi(G)$. Assume that $H,K\in
\Hall_\pi(G)$. If
$H$ and $K$ are conjugate in $G$, then it is clear that $HA/A$ and $KA/A$ are conjugate in $G/A$. Assume that $HA/A$ and
$KA/A$ are
conjugate in $G/A$. Then, up to conjugation in $G$, we have that $HA=KA$. Since $A\in C_\pi$ we may assume that $H\cap
A=K\cap A$, i.e.,
$H,K\leq N_{HA}(H\cap A)$. As in the proof of item (e) we obtain that $N_{HA}(H\cap A)\in D_\pi$, hence $H,K$ are
conjugate. \qed

If $\prec$ is a total ordering of $\pi(G)$, then $G$ has a Sylow tower of complexion $\prec$
provided $G$ has a normal series $G=G_0>G_1>\ldots>G_n=1$,
where $G_{i-1}/G_i$ is isomorphic to a Sylow $r_i$-subgroup of $G$ and $r_1\prec r_2\prec\ldots \prec r_n$. %If $\prec$
%is the natural ordering of $\pi(G)$, then we say simply that $G$ has a Sylow tower.

\begin{Lemma}\label{SylowTower}
{\em \cite[Theorem~A1]{HallProperties}} Assume that $H_1,H_2$ are two $\pi$-Hall subgroups of a finite group $G$
such that both $H_1,H_2$ have a Sylow tower of the same complexion. Then $H_1$ and $H_2$ are conjugate in~$G$.
\end{Lemma}

\begin{Lemma}\label{HallSymmetric} {\em (\cite[Theorem~A4 and the
note after it]{HallProperties}, \cite[Theorem 4.3 and
Corollary 4.4]{RevVdoContemp}, \cite{ThomHallSymmetric})}
Let $\pi$ be a set of primes. Then the following statements hold{\em:}
\begin{itemize}
\item[{\em (A)}] If $\Sym_n\in E_\pi$ and $H$ is a $\pi$-Hall subgroup of $\Sym_n$, then $n$, $\pi$, and $H$ satisfy
exactly one of
the following statements{\em:}
\begin{itemize}
\item[{\em (a)}] $\vert \pi\cap\pi(\Sym_n)\vert\le 1$. In this case a
$\pi$-Hall subgroup of $\Sym_n$ is either its Sylow $p$-subgroup (if
$\pi\cap\pi(\Sym_n)=\{p\}$) or trivial (if $\pi\cap\pi(\Sym_n)=\emptyset$).
\item[{\em (b)}] $n=p\ge7$ is a prime and $\pi\cap \pi(\Sym_p)=\pi((p-1)!)$. In
this case a $\pi$-Hall subgroup of $\Sym_p$ is isomorphic to $\Sym_{p-1}$.
\item[{\em (c)}] $\pi(\Sym_n)\subseteq \pi$ and $n\ge5$. In this case $\Sym_n$ is a
$\pi$-Hall subgroup of~$\Sym_n$.
\item[{\em (d)}] $\pi\cap\pi(\Sym_n)=\{2,3\}$ and $n\in\{3,4,5,7,8\}$. In this
case for a $\pi$-Hall subgroup $H$ of $\Sym_n$ we have that $H=\Sym_3$ if $n=3$,
$H\simeq \Sym_4$ if $n\in\{4,5\}$, $H\simeq \Sym_3\times \Sym_4$ if $n=7$, and
$H\simeq \Sym_4\wr\Sym_2$ if~${n=8}$.
\end{itemize}
\item[{\em (B)}] Conversely, if $n$ and $\pi$ satisfy one of statements {\em (a)--(d)}, then $\Sym_n\in E_\pi$.
\item[{\em (C)}] $\Hall_\pi^{\Sym_n}(\Alt_n)=\Hall_\pi(\Alt_n)$ and
$k_\pi(\Sym_n)=k_\pi^{\Sym_n}(\Alt_n)=k_\pi(\Alt_n)\in\{0,1\}$.
\item[{\em (D)}] If $2,3\in\pi$ and $\pi(n!)\not\subseteq\pi$, then both $\Alt_n$ and $\Sym_n$ do not satisfy~$D_\pi$.
\end{itemize}
\end{Lemma}

\begin{Cor}\label{NormalizersOfHallSubgroupsInSymmetric}
Assume that $\pi$ is a set of primes such that $2,3\in\pi$. Suppose that $\Sym_n\in E_\pi$ and $H\in\Hall_\pi(\Sym_n)$.
Then
$N_{\Sym_n}(H)=H$.
\end{Cor}

\begin{Lemma}\label{Numberkpi(G)Inextension}
Let $M$ be a $\pi$-subgroup of $\Sym_n$, $L$ be a finite group, and $$L\wr M=(L_1\times\ldots\times L_n)\splitext M$$ be
a permutation
wreath product and $L\simeq L_i$ for $i=1,\ldots,n$. Assume that $G$ is a normal subgroup of $L\wr M$ and $A=G\cap
(L_1\times\ldots\times L_n)$.
Suppose also that $A\simeq
L_1\circ \ldots \circ L_n$ and $G/A\simeq M$. Denote by $t$ the number of orbits of $M$, and by $k$ the number
$k_\pi(L)$.
Then~${k_\pi(G)=k_\pi(L\wr M)=k^t}$. Moreover, if $2,3\in\pi$ and $M$ is a $\pi$-Hall subgroup of $\Sym_n$, then
$k_\pi(G)=k_\pi(L\wr
M)=k_\pi(L\wr \Sym_n)$.
\end{Lemma}

\noindent{\slshape Proof.}\ \
By Lemmas \ref{simpleepi}(f) and \ref{HallSymmetric}(C) we may assume $k_\pi(L)\ge1$, i.e., $L\in E_\pi$. In view of
Lemma \ref{simpleepi}(f) we may assume that $Z(A)=1$, i.e.,
$A=L_1\times\ldots\times L_n$ and $G=L\wr M$. Assume that
$H,K\in\Hall_\pi(G)$. First we prove \begin{multline}\label{ConjHall}
H,K \text{ are conjugate in } G\text{  if and only if }\\ H\cap L_i, K\cap L_i\text{ are conjugate in } L_i\text{ for
every }i=1,\ldots,n.
\end{multline}
Assume that $H=K^g$ for some $g\in G$. Condition $\pi(M)\subseteq \pi$ implies $HA=KA=G$. So we may suppose $g\in A$ and
$g=g_1\ldots
g_n$, where $g_i\in L_i$. Since $[L_i,L_j]=1$ for $i\not=j$, we obtain $$H\cap L_i=K^g\cap L_i=(K\cap L_i)^g=(K\cap
L_i)^{g_i}$$ for
every $i=1,\ldots,n$. Conversely suppose that $H\cap L_i$ and $K\cap L_i$ are conjugate in $L_i$ for every
$i=1,\ldots,n$. It follows that
$H\cap A$ and $K\cap A$ are conjugate in $A$. Therefore we may assume that $H\cap A=K\cap A$, so $H,K\leq N_G(H\cap A)$.
As in the proof
of item (e) of Lemma \ref{simpleepi} we obtain $N_G(H\cap A)\in D_\pi$, hence $H,K$ are conjugate.

Let $\Omega_1,\ldots,\Omega_t$ be the orbits of $M$. Since, for every $H\in\Hall_\pi(G)$  we have $G=HA$ and $M\simeq
G/A$, it
follows that $\Omega_1,\ldots,\Omega_t$ are the orbits of $H$  on $\{L_1,\ldots,L_n\}$. Now for every $h\in H$ and
$i=1,\ldots,n$ the
identity $H\cap L_i^h=(H\cap L_i)^h$ holds, so for every
$j=1,\ldots,t$ the subgroup $\prod_{L_i\in\Omega_j} L_i$ possesses at most $k$  classes of conjugate $\pi$-Hall subgroup
invariant under $G$.
Thus Lemma \ref{simpleepi}(e) implies the inequality~${k_\pi(G)\le k^t}$.

Conversely assume that $H\in \Hall_\pi(L\wr M)$ (as we noted at the beginning of the proof, we may assume $G=L\wr M$).
If we take $L_{j_i}\in\Omega_i$ and $H_{j_i}\in \Hall_\pi(L_{j_i})$ for some
$i=1,\ldots ,t$, then for every $g_1,g_2\in M$ with $L_{j_i}^{g_1}=L_{j_i}^{g_2}$ we have
$(H_{j_i}^{L_{j_i}})^{g_1}=(H_{j_i}^{L_{j_i}})^{g_2}$. So for every
$i=1\ldots,t$ in the subgroup $\prod_{L_j\in\Omega_i} L_j$ we can take at least $k$  classes of conjugate
$\pi$-Hall subgroups invariant under $L\wr M$. Hence $A$ possesses at least $k^t$ classes of conjugate
$\pi$-Hall subgroups invariant under $L\wr M$. By Lemma \ref{simpleepi}(e) and statement \eqref{ConjHall} we obtain that
$L\wr M$ possesses at least $k^t$
classes of conjugate $\pi$-Hall subgroups, whence~${k^t\ge k_\pi(G)\ge  k^t}$.

Suppose that $2,3\in\pi$, $\Sym_n\in E_\pi$, and $M\in \Hall_\pi(\Sym_n)$. Suppose $H,K\in\Hall_\pi(L\wr\Sym_n)$. Since
$\Sym_n\in C_\pi$,
we may assume that $H,K\in\Hall_\pi(L\wr M)$. So we need to show that $H,K$ are conjugate in $L\wr M$ if and only if
$H,K$ are conjugate
in $L\wr \Sym_n$. Suppose there exists $x\in L\wr \Sym_n$ such that $H=K^x$. Then $L\wr M=(L_1\times\ldots\times L_n)H
=(L_1\times\ldots\times
L_n)K$, hence the image of $x$ under the natural homomorphism $L\wr M\rightarrow M$ is in $N_{\Sym_n}(M)$, which is
equal to $M$ in virtue
of Corollary \ref{NormalizersOfHallSubgroupsInSymmetric}. Therefore~${x\in L\wr M}$.
\qed

\begin{Lemma}\label{FWair} {\em (\cite[Lemma~2.2]{GrossNotDivisibleBy3}, \cite[Lemmas~2.4 and~2.5]{GrossOddOrder}, and
\cite{Wair})}
Let $q$ be a rational integer and $r$ a prime such that  $(q,r)=1$. Denote $e(q,r)$ by $e$ and
$$e^*=\left\{
\begin{array}{rl}
2e, & \text{ if } e\equiv 1\pmod 2,\\
e, & \text{ if } e\equiv 0\pmod 4,\\
e/2, & \text{ if } e\equiv 2\pmod 4.
\end{array}
\right.
$$
Then the following identities hold{\em:}
$$
(q^n-1)_r=
\left\{
\begin{array}{ll}
(q^e-1)_r (n/e)_{r}, & \text{ if } n \text{ is divisible by } e,\\
&\\
(r,2) & \text{ otherwise{\em;} }
\end{array}
\right.
$$

$$
(q^n-(-1)^n)_r=
\left\{
\begin{array}{ll}
(q^{e^*}-(-1)^{e^*})_r (n/{e^*})_{r}, & \text{ if } n \text{ is divisible by } e^*,\\
&\\
(r,2) & \text{ otherwise{\em;} }
\end{array}
\right.
$$

$$
\prod\limits_{i=1}^{n}(q^i-1)_r=(q^e-1)^{n/e}_r\big((n/e)!\big)_r.
$$
\end{Lemma}

If $r\in \{2,3\}$, then Lemma \ref{FWair} implies the following corollaries.

\begin{Cor}\label{3part} Let $q$ be a rational integer such that $(q,3)=1$. Then the following identities hold{\em:}
$$
(q^n-\eta^n)_\ti=
\left\{
\begin{array}{ll}
(q-\eta)_\ti n_\ti, & \text{ if } q\equiv\eta\pmod 3,\\
&\\
(q+\eta)_\ti\left(%\displaystyle\frac{n}{2}
n/2\right)_\ti, & \text{ if } q\equiv-\eta\pmod 3 \text{ and } n \text{ is even},\\
&\\
1, & \text{ if } q\equiv-\eta\pmod 3 \text{ and } n \text{ is odd};
\end{array}
\right.
$$

$$
\prod\limits_{i=1}^{n}(q^i-\eta^i)_\ti=\left\{
\begin{array}{ll}
(q-\eta)^n_\ti (n!)_\ti, & \text{ if } q\equiv\eta\pmod 3,\\
&\\
(q+\eta)^{[n/2]}_\ti\left(%\displaystyle\frac{n}{2}
[n/2]!\right)_\ti, & \text{ if } q\equiv-\eta\pmod 3.
\end{array}
\right.
$$
\end{Cor}

\begin{Cor}\label{2part}  Let $q$ be an odd rational integer. Then the following identities hold{\em:}
$$
(q^n-\eta^n)_\tw=(q-\eta)_\tw t, \text{ where }
t=\left\{
\begin{array}{ll}
(q+\eta)_\tw (n/2)_{\mbox{}_2}, & \text{ if }
n \text{ is even},\\

&\\
1, & \text{ if }
 n \text{ is odd{\em;}}

\end{array}
\right.
$$
$$
\prod\limits_{i=1}^{n}(q^i-\eta^i)_\tw=
(q-\eta)^n_\tw(q+\eta)^{[n/2]}_\tw \big([n/2]!\big)_{\mbox{}_2}.
$$
\end{Cor}

\begin{Lemma}\label{SymmetricDivizors}
Assume that $r,p$ are distinct odd primes, $q=p^\alpha$ and
$m\ge\frac{r+1}{2}$. Then the inequality $\big((q^2-1)(q^4-1)\ldots(q^{2(m-1)}-1)\big)_r>(m!)_r$ holds.
\end{Lemma}

\noindent{\slshape Proof.}\ \
Note that
\begin{equation*}
(m!)_r=\prod_{k=1}^{[m/r]}r\cdot k_r.
\end{equation*}
To every number $r\cdot k$ we put into correspondence the number $(q^{2k(r-1)}-1)$. Then
$$(q^{2k(r-1)}-1)_r=(q^{2(r-1)}-1)_r\cdot k_r\ge
r\cdot k_r.$$
So the following inequalities
\begin{multline*}
\big((q^2-1)(q^4-1)\ldots(q^{2(m-1)}-1)\big)_r\ge
(q^{2((r-1)/2)}-1)_r\cdot \prod_{k=1}^{[m/r]}(q^{2k(r-1)}-1)_r\ge
\\ r\cdot \prod_{k=1}^{[m/r]}r\cdot k_r>(m!)_r,
\end{multline*}
hold, whence the lemma follows.
\qed

For linear algebraic groups our notation agrees with that of
\cite{HuLinearAlgebraicGroups}. For finite groups of Lie type we use the notation
from \cite{CarSimpleGgroupsOfLieType}. If $\ov{G}$~ is a simple connected
linear algebraic group over the algebraic closure $\ov{\F}_p$ of a finite field
of characteristic $p$, then a surjective endomorphism $\sigma:\ov{G}\rightarrow
\ov{G}$ is called a {\em Frobenius map}, if the set of $\sigma$-stable points
$\ov{G}_\sigma$ is finite. Every group $G$ such that $O^{p'}(\ov{G}_\sigma)\leq
G\leq \ov{G}_\sigma$ is called a {\em finite group of Lie type}.
Notation for classical groups agrees with~\cite{KlLi}. In order to make uniform
statements and arguments we use the following notations $\GL_n^+(q)=\GL_n(q)$,
$\GL_n^-(q)=\GU_n(q)$, $\SL_n^+(q)=\SL_n(q)$, $\SL_n^-(q)=\SU_n(q)$. In this paper we consider groups of Lie type with
the base field $\F_q$
of odd
characteristic $p$ and order $q=p^\alpha$, and we fix the symbols $p$ and $q$ for this purposes. We always choose
$\varepsilon(q)\in\{+1,-1\}$ (usually we write just $\varepsilon$, since $q$ is fixed by the choice of $G$) so that
$q\equiv\varepsilon(q)\pmod4$, i.~e, $\varepsilon(q)=(-1)^{(q-1)/2}$. The same symbol $\varepsilon(q)$ is used to denote
the sign of $\varepsilon(q)$.
Following \cite{KlLi},  by $\Oo^\eta_n(q)$ we denote the general
orthogonal group of degree $n$ and of sign
$\eta\in\{\circ,+,-\}$ over $\F_q$, while the symbol $\GO^\eta_n(q)$ denotes the group of similarities. Here $\circ$ is
an empty symbol, and
we use it only if $n$ is odd. By $\eta$ we always mean an element  from the set $\{\circ,+,-\}$, and if
$\eta\in\{+,-\}$, then we use $\eta$
instead of $\eta1$ as well. In classical groups the symbol $\P$ will also denote the reduction modulo scalars. Thus for
every subgroup $H$ of
$\GL_n(q)$ the image of $H$ in $\PGL_n(q)$ is denoted by~${\P H}$.

Let $I\in\{\GL_n(q),\GU_n(q),\Sp_n(q),\Oo^\eta_n(q)\}$ be a general classical group over a field of   odd characteristic
$p$. Let $V$ be the
natural module for $I$ of dimension $n$
over a field either of order $q$ or of order $q^2$ if $I=\GU_n(q)$. In all cases we say that $\F_q$ is the base field
for $V$. Assume that
$V$ is equipped with the corresponding form
(trivial for $\GL_n(q)$, unitary for
$\GU_n(q)$, skew-symmetric for $\Sp_n(q)$, and symmetric for $\Oo_n^\eta(q)$). Then $I$ can be identified  in a natural
way with
the  {\it (general) group of isometries} of $V$. We set also
\begin{multline*}
S=S(V)= I(V)\cap \SL(V)\in\{\SL(V),\SU(V),\Sp(V),{\rm SO}^{\eta}(V)\}
\text{ and }\\
\Omega=\Omega(V)= O^{p'}(S(V)).
\end{multline*}
Given subspaces $U,W\leq V$, we write $U+ W=U\perp W$ if $U\cap W=\{0\}$ and $U,W$ are mutually orthogonal. Following
\cite{KlLi}, we say
that a subgroup $H$ of $G$, where $\Omega(V)\leq G\leq I(V)$, is of type $I(U)\perp I(W)$, if $H$ is the stabilizer in
$G$ of the
decomposition $U\perp W$ and $H$ stabilizes both $U,W$, while $H$ is isomorphic to $(I(U)\times I(W))\cap G$  as an
abstract group .
If  $V$ is an orthogonal space of even dimension, then we denote the sign of the corresponding quadratic form by
$\eta(V)$, and the
discriminant of the form  by $D(V)$. We write
$D(V)=\square$ if the discriminant of the form is a square in $\F_q$,  and $D(V)=\boxtimes$ if the discriminant of the
form is a non-square
in $\F_q$. If $\eta=\eta(V)$ and $\dim(V)=2m$ then  \cite[Proposition~2.5.10]{KlLi} implies that $D(V)=\square$ if and
only if
$\eta=\varepsilon(q)^m$.

In our arguments we use the classification of subgroups of odd index in finite
simple groups obtained by M.W.Liebeck and J.Saxl \cite{liesax} and
independently by W.M.Kantor \cite{Kant}. A more detailed description of subgroups
of odd index in the finite simple classical groups is obtained by N.V.Maslova
in~\cite[Theorem~1]{Maslova} and for classical groups we refer to this
description. Since \cite{Maslova} is published in Russian, we cite the main theorem from this paper here.

Assume that $n$ is a positive integer and $\alpha_0\cdot 2^0+\alpha_1\cdot 2^1+\ldots$, where $\alpha_i\in\{0,1\}$, is
the $2$-adic
expansion of
$n$ (for our purposes we assume that this expansion is infinite, but only finitely many coefficients are not equal to
$0$). Define
$\psi(n)=(\alpha_0,\alpha_1,\ldots )$. Let $\ge$ be a linear order on $\{0,1\}$ such that $1\ge 0$. We say that
$\psi(n)=(\alpha_0,\alpha_1,\ldots)\gg \psi(m)=(\beta_0,\beta_1,\ldots)$ if $\alpha_i\ge \beta_i$ for all $i$. Notice
that
$\gg$ is a partial order.

\begin{Theo}\label{Maslova0}{\em \cite[Theorem~1]{Maslova}}
Let $G$ be one of the finite classical groups{\em :} $\SL_n(q)$ with $n \ge 2$,
$\SU_n(q)$ with $n \ge 3$, $\Sp_n(q)$ with $n \ge 4$ and $n$ even, $\Omega_n(q)$ with $n \ge 7$ and $n$ odd, and
$\Omega_n^\pm(q)$ with $n \ge 8$ and $n$ even. Assume that the base field of $G$ has odd order $q$, and $V$ is the
natural module for
$G$. Then
$M$ is a maximal subgroup of odd index of $G$ if\footnote{In the original paper the theorem is proven in both
directions, i.e., it has ``if and only if'' form, however the proof of ``only if'' part does use the unpublished PhD
thesis of P.Kleidman, that known to have inaccuracies} one of the following statements holds{\em:}
\begin{itemize}
\item[{\rm (a)}] $M=N_G(C_G(\sigma))$, where $\sigma$ is a field automorphism of odd prime order of~$G$.
\item[{\rm (b)}] $G=\SL_n(q)$, $M$ is the stabilizer of a subspace of dimension $m$ of $V$ and $\psi(n) \gg
\psi(m)$.
\item[{\rm (c)}] $G=\SU_n(q)$ or $G=\Sp_n(q)$, $M$ is the stabilizer of a nondegenerate subspace of dimension~$m$ of $V$
and $\psi(n)
\gg \psi(m)$.
\item[{\rm (d)}] $G=\Omega_n(q)$, $n$ is odd, $M$ is the stabilizer of a nondegenerate subspace $U$ of even
dimension~$m$ of $V$, $D(U) =
\square$,
$\psi(n) \gg \psi(m)$ and $(q,m,\eta(U)) \not = (3,2,+)$.
\item[{\rm (e)}] $G=\Omega_n^\eta(q)$, $n$ is even, $M$ is the stabilizer of a nondegenerate subspace  $U$ of
dimension~$m$ of $V$, and one of the following holds:
    \begin{itemize}
    \item[{\rm (e.1)}] $m$ is odd, $D(V)=\boxtimes$, and $\psi(n-2) \gg \psi(m-1)$, except for the case  $m=n/2$ and
subspaces $U$ and  $U^\bot$ are nonisometric;
    \item[{\rm (e.2)}] $m$ is even, $(q,m,\eta(U)) \not = (3,2,+),(3,n-2,+)$,  $D(U)=D(V)=\boxtimes$, and $\psi(n-2) \gg
\psi(m-2)$;
    \item[{\rm (e.3)}] $m$ is even, $(q,m,\eta(U)) \not = (3,2,+),(3,n-2,+)$, $D(U)=D(V)=\square$, and $\psi(n) \gg
\psi(m)$.
    \end{itemize}
\item[{\rm (f)}] $G=\SL_n(q)$, $M$ is the stabilizer of a decomposition $V=\bigoplus V_i$ into a direct sum of subspaces
of the same
dimension $m$
and either $m=2^w \ge 2$ and $(n,m,q) \not = (4,2,3)$, or $m=1$, $q \equiv 1 \pmod 4$ and $q \ge 13$ if $n=2$.
\item[{\rm (g)}] $G=\SU_n(q)$, $M$ is the stabilizer of an orthogonal decomposition $V=\bigoplus V_i$ into a direct sum
of isometric
subspaces
$V_i$ of dimension $m$, and either $m=2^w \ge 2$, or $m=1$, $q \equiv 3 \pmod 4$ and $(n,q) \not = (4,3)$.
\item[{\rm (h)}] $G=\Sp_n(q)$, $H$ is a stabilizer of an orthogonal decomposition $V=\bigoplus V_i$ into a direct sum of
isometric subspaces
$V_i$
of dimension $m$ and $m=2^w \ge 2${\rm.}
\item[{\rm (i)}] $G=\Omega_n(q)$, $n$ is odd, $M$ is the stabilizer of an orthogonal decomposition $V=\bigoplus V_i$
into a direct sum of
isometric
subspaces
$V_i$ of dimension $1$,  $q$ is a prime, and $q \equiv \pm 3 \pmod 8$.
\item[{\rm (j)}] $G=\Omega_n^\eta(q)$, $n$ is even, $M$ is the stabilizer of an orthogonal decomposition $V=\bigoplus
V_i$ into a direct sum
of
isometric
subspaces   $V_i$ of dimension $m$, and either $m=1$, $q$ is a prime, $q \equiv \pm 3 \pmod 8$ and $(n,\eta) \not =
(8,+)$; or
$m=2^w \ge 2$, $D(V)=D(V_i)=\square$, and $(m,q,\eta(V_i)) \not = (2,3,\pm),(2,5,+)$.
\item[{\rm (k)}] $\P G=\PSL_2(q)$ and $\P M \simeq \PGL_2(q_0)$, where  $q=q_0^2$.
\item[{\rm (l)}] $\P G=\PSL_2(q)$ and $\P M \simeq \Alt_4$, where $q$ is a prime,  $q \equiv 3,5,13,27,37 \pmod {40}$.
\item[{\rm (m)}] $\P G=\PSL_2(q)$ and $\P M \simeq \Sym_4$, where $q$ is a prime, $q \equiv \pm 7 \pmod {16}$.
\item[{\rm (n)}] $\P G=\PSL_2(q)$ and $\P M \simeq \Alt_5$, where  $q$ is a prime, $q \equiv 11,19,21,29 \pmod {40}$.
\item[{\rm (o)}] $\P G=\PSL_2(q)$ and $\P M$ is the dihedral group of order $q+1$, where  $7 < q \equiv 3 \pmod 4$.
\item[{\rm (p)}] $\P G=\PSU_3(5)$ and $\P M \simeq M_{10}\simeq 3\arbitraryext \Alt_6\arbitraryext2$.
\item[{\rm (q)}] $\P G=\PSL_4(q)$ and $\P M \simeq 2^4\arbitraryext\Alt_6$, where $q$ is a prime, $q \equiv 5 \pmod 8$.
\item[{\rm (r)}] $\P G=\PSL_4(q)$ and $\P M \simeq \PSp_4(q)\arbitraryext 2$, where $q \equiv 3 \pmod 4$.
\item[{\rm (s)}] $\P G=\PSU_4(q)$ and $\P M \simeq 2^4\arbitraryext\Alt_6$, where $q$ is a prime, $q \equiv 3 \pmod 8$.
\item[{\rm (t)}] $\P G=\PSU_4(q)$ and $\P M \simeq \PSp_4(q)\arbitraryext 2$, where $q \equiv 1 \pmod 4$.
\item[{\rm (u)}] $\P G=\PSp_4(q)$ and $\P M \simeq 2^4\arbitraryext \Alt_5$, where $q$ is a prime, $q \equiv 3 \pmod 8$.
\item[{\rm (v)}] $\P G=\P \Omega_7(q)$ and $\P M \simeq \Omega_7(2)$, where $q$ is a prime, $q \equiv \pm 3 \pmod 8$.
\item[{\rm (w)}] $\P G=\P \Omega_8^+(q)$ and $\P M \simeq \Omega_8^+(2)$, where $q$ is a prime, $q \equiv \pm 3 \pmod
8$.
\end{itemize}
\end{Theo}

\begin{Lemma}\label{StructureInducedAutomorphisms} Let  $V$ be a vector space of dimension $n$ over a finite field
$\F_q$ of odd
characteristic $p$, equipped with a trivial, unitary, symmetric, or skew-symmetric form. Assume that
$V=V_1\perp\dots\perp V_k$ is a
decomposition of $V$ into a direct orthogonal sum of nondegenerate (arbitrary if the form is trivial) subspaces. Let
$L=\Omega(V)$,
$L_0=\left(
I(V_1)\times\dots\times I(V_k)\right)\cap L$, and
$\rho_i:I(V_1)\times\dots\times I(V_k)\rightarrow I(V_i)$ be the natural projection. Then
$L_0^{\rho_i}=I(V_i)$ for every
$i=1,\ldots,k$, except for the case, when   $V$ is an orthogonal space, $k=2$ and, up to renumbering,  $\dim(V_1)=n-1$,
$\dim(V_2)=1$, and $i=1$. In this exceptional case one of the following statements holds{\em:}
\begin{itemize}
\item[{\em (a)}]  $n$ is odd, $L_0^{\rho_1}$ is equal to $\Omega(V_1)$, extended by a graph automorphism.
\item[{\em (b)}]  $n$ is even, $D(V)=\square$,  $L_0^{\rho_1}=\Omega(V_1)$.
\item[{\em (c)}]  $n$ is even, $D(V)=\boxtimes$,  $L_0^{\rho_1}={\rm SO}(V_1)$.
\end{itemize}
\end{Lemma}

\noindent{\slshape Proof.}\ \  Assume that $L_0^{\rho_i}\ne I(V_i)$. Then  \cite[Lemma~4.1.1]{KlLi} implies that $V$ is
an orthogonal space,
$k=2$, up to renumbering,    $\dim( V_1)=n-1$, $\dim( V_2)=1$, and $i=1$. Now for $n$ even the statement of the lemma
follows from
\cite[Proposition~4.1.6]{KlLi} and the fact that  $\mathrm{O}_{n-1}(q)\times\mathrm{O}_1(q)$ does not induce  graph
automorphisms
on~$\Omega_{n-1}(q)$.

Assume that $n$ is odd. Let $K$ be an algebraic closure of $\F_q$, $\overline{L}=\mathrm{SO}_n(K)$ a simple linear
algebraic group
(the group of orthogonal matrices of determinant~$1$), $\mathrm{O}_n(K)$ a group of all orthogonal matrices. Assume that
$\sigma$ is a Frobenius map of $\overline{L}$, i.e., a surjective endomorphism such that the set of $\sigma$-stable
points
$\overline{L}_\sigma$ is finite.  Then $\overline{L}_\sigma=\mathrm{SO}_n(q)$ and
$L=O^{p'}(\overline{L}_\sigma)=\Omega_n(q)$, the index  $\vert \overline{L}_\sigma:G\vert$ is equal to $2$ (recall that
$q$ is odd)
and
$\overline{L}_\sigma=\mathrm{SO}_n(q)$ is generated by  $L$ and a diagonal  automorphism.

Let $V$ be the natural module for $L$. Then $\overline{V}=K\otimes_{\F_q}V$ is the natural module for~$\overline{L}$.
Moreover, if  $V_1\perp  V_2$ is a decomposition of  $V$ into an orthogonal direct sum, then
$(K\otimes_{\F_q}V_1)\perp (K\otimes_{\F_q}V_2)$ is a decomposition of  $\overline{V}$ into an orthogonal direct sum.
So, for every
subgroup $L_0$ of $L$, stabilizing the decomposition $V_1\perp V_2$,  there corresponds a unique subgroup
$\overline{L}_0$ of $\overline{L}$,
stabilizing the decomposition $(K\otimes_{\F_q}V_1)\perp  (K\otimes_{\F_q}V_2)$, and $L_0=\overline{L}_0\cap L$. The
subgroup $\overline{L}_0$
is a reductive subgroup of maximal rank of  $\overline{L}$. By \cite[Theorem~2]{vdogal} it follows that 
$\mathrm{Aut}_L(L_0^{\rho_1})$
does not contain diagonal automorphisms, hence by using \cite[Proposition~4.1.6]{KlLi}, we obtain the statement of the
lemma for  $n$ odd.
\qed

\begin{Lemma}\label{notinherit}
Assume that a simple classical group  $G$ and its subgroup $H$ of  satisfy  one of the following statements{\em:}
\begin{itemize}
 \item[{\em (a)}]\ \  $G\simeq \PSL_2(q)$, $(q,6)=1$, $H\simeq \Sym_4 ${\em;}
 \item[{\em (b)}]\ \  $G\simeq \PSL_2(q)$, $(q,30)=1$, $ H\simeq\SL_2(4)${\em;}
 \item[{\em (c)}]\ \  $G\simeq  {\rm P}\Omega_7(q)$, $(q,210)=1$, $ H\simeq\Omega_7(2)${\em;}
 \item[{\em (d)}]\ \  $G={\rm P}\Omega_8^+(q)$, $(q,210)=1$, $H\simeq \Omega_8^+(2)$.
\end{itemize}
Suppose $K$ is chosen so that $G<K\leq \widehat{G}$, where $\widehat{G}$ is the group of inner-diagonal automorphisms of
$G$.  Then there is no subgroup $H_1$ of  $K$ such that $H_1\cap G=H$ and~$|H_1:H|=|K:G|$.
\end{Lemma}

\noindent{\slshape Proof.}\ \
If $H$ satisfies either (a) or (b), then the lemma follows from \cite[Chapter~II, \S~8]{Hup}. If $H$ satisfies statement
(d), then the
lemma follows from  \cite[Proposition~2.3.8]{Kl}. Assume that $H$ satisfies (c). By using   \cite{ATLAS}, we obtain that
the minimal
nontrivial irreducible representation of $\Omega_7(2)$ has degree $7$.  Since
 $(\vert H\vert, q)=1$, all ordinary characters of this group have rational values,
  and in view of \cite[Theorem~9.14 and Corollary~15.12]{Isaacs}, the same property holds for the representations over
$\F_q$. In view of \cite{ATLAS} it also follows that  $\Out (\Omega_7(2))$ is trivial, while the universal central
extension
$2\arbitraryext\Omega_7(2)$ of $\Omega_7(2)$ has no faithful irreducible representations of degree 7 over $\F_q$.
Therefore
$N_{\widehat{G}}(\Omega_7(2))=N_G(\Omega_7(2))=\Omega_7(2)$.
\qed

The next lemma follows from Lemmas \ref{simpleepi}(e), \ref{StructureInducedAutomorphisms} and \ref{notinherit}.

\begin{Lemma}\label{HallAmostSimpleInduce}
Assume $\pi\cap\pi(G)=\{2,3,5,7\}$, $G$ is isomorphic to either  $\Omega_7(q)$, or $\Omega_8^+(q)$, or $\Omega_9(q)$,
and a $\pi$-Hall
subgroup  $H$ of $G$ is isomorphic to either $\Omega_7(2)$, or $2\arbitraryext\Omega_8^+(2)$, or $\left(2\arbitraryext
\Omega_8^+(2)\right)\arbitraryext 2$, respectively. Denote by $G_1$ either
$\mathrm{SO}_7(q)$, or $\mathrm{SO}_8^+(q)$, or $\mathrm{SO}_9(q)$, respectively. Then  $G_1$ does not possesses a
$\pi$-Hall subgroup
$H_1$ such that~${H_1\cap G=H}$.
\end{Lemma}

\section[Classical groups in small dimension]{Maximal subgroups of odd index in classical groups of small dimension}

In this section we classify $\pi$-Hall subgroups in groups $\SL_2^\eta(q)$ and $\GL_2^\eta(q)$, and give a
complete list of maximal subgroups of odd index in classical groups of small dimension.

\begin{Lemma}\label{SLU23dim2} Let $\pi$ be a set of primes with $2,3\in\pi$.  Assume that
$G\simeq\SL_2(q)\simeq\SL_2^\eta(q)\simeq\Sp_2(q)$,  where
$q$ is a power of an odd prime $p\not\in\pi$, and $\varepsilon=\varepsilon(q)$.  Recall that for a subgroup $A$ of $G$
we denote by
$\P A$ the reduction modulo scalars. Then the following statements hold{\em:}

\begin{itemize}
\item[{\em (A)}] If $G\in E_\pi$ and $H\in \Hall_\pi(G)$, then one of  the following statements holds{\em:}
\begin{itemize}
\item[{\em (a)}]  $\pi\cap\pi(G)\subseteq \pi(q-\varepsilon)$, $\P H$ is a
$\pi$-Hall subgroup in the dihedral subgroup $D_{q-\varepsilon}$ of order $q-\varepsilon$ of $\P G$.
All $\pi$-Hall subgroups of this type are conjugate in~$G$.

\item[{\em (b)}] $\pi\cap\pi(G)=\{2,3\}$, $(q^2-1)_{\mbox{}_{\{2,3\}}}=24$, $\P H\simeq \Alt_4$.
All  $\pi$-Hall subgroups of this type are conjugate in~$G$.

\item[{\em (c)}] $\pi\cap\pi(G)=\{2,3\}$, $(q^2-1)_{\mbox{}_{\{2,3\}}}=48$, $\P H\simeq \Sym_4$.
There exist exactly two classes of conjugate subgroups of this type, and
$\PGL_2^\eta(q)$ interchanges  these classes.

\item[{\em (d)}] $\pi\cap\pi(G)=\{2,3,5\}$, $(q^2-1)_{\mbox{}_{\{2,3,5\}}}=120$, $\P H\simeq
\Alt_5$. There exist exactly two classes of conjugate subgroups of this type, and
$\PGL_2^\eta(q)$ interchanges  these classes.
\end{itemize}
\item[{\em (B)}] Conversely, if $\pi$ and $(q^2-1)_\pi$ satisfy one of statements {\em (a)--(d)}, then~${G\in E_\pi}$.
\item[{\em (C)}] If $G\in E_\pi$, then  $k_\pi(G)\in\{1,2,3\}$.
\item[{\em (D)}] If $G\in E_\pi$ and $H\in\Hall_\pi(G)$, then  there exists a $\pi$-subgroup $K\leq G$, nonconjugate in
$G$ to a
subgroup of~$H$, in particular $G\not\in D_\pi$.
\item[{\em (E)}] Every $\pi$-Hall subgroup of $\P G$ can be obtained as $\P H$ for some~${H\in \Hall_\pi(G)}$.
Conversely, if $\P H\in
\Hall_\pi(\P G)$ and $H$ is a complete preimage of $\P H$ in $G$, then~$H\in \Hall_\pi(G)$.
\end{itemize}
\end{Lemma}

\noindent{\slshape Proof.}\ \
(E) Follows from Lemma \ref{simpleepi}(a), (d) and from $\vert Z(G)\vert=2$.

(A) Assume that $\P H$ is a $\pi$-Hall subgroup of $\PSL_2(q)$. Then $\P H$ is
included in a maximal subgroup $M$ of odd index of $\PSL_2(q)$. By using
Theorem \ref{Maslova0}, we obtain the following list of maximal subgroups of
odd index in $\PSL_2(q)$:
\begin{enumerate}
\item $M=N_G(C_G(\sigma))$, where $\sigma$ is a field automorphism of odd prime order of~$G$. In view of
\cite[Proposition~4.5.3]{KlLi} we
obtain that $M\simeq \PSL_2(q^{\frac{1}{r}})$, where $q=p^\alpha$ and $r$ is an odd prime divisor of~$\alpha$, and all
subgroups of this
type are conjugate in~$\PSL_2(q)$.
\item  $q=q_0^2$, $M\simeq \PGL_2(q_0)$.
\item $q\equiv1\pmod4$, $q\ge 13$, $M$ is the dihedral group $D_{q-1}$ of order $q-1$.
\item $q\equiv-1\pmod4$, $q\ge 11$, $M$ is the dihedral group $D_{q+1}$ of order $q+1$.
\item $q\equiv 3,5,13,27,37\pmod{40}$, $M\simeq \Alt_4$.
\item $q\equiv \pm7\pmod{16}$, $M\simeq \Sym_4$.
\item $q\equiv 11,19,21,29\pmod{40}$, $M\simeq \Alt_5$.
\end{enumerate}

Consider all these cases separately.

Assume that either $\P H\leq  \PSL_2(q^{\frac{1}{r}})=\PSL_2(q_0)$, or $\P H \leq
\PGL_2(q^{\frac{1}{2}})= \PGL_2(q_0)$.
The condition $p\not\in\pi$ implies that  $\PSL_2(q_0) \not \subseteq\P H$. By
Lemma
\ref{simpleepi}(a) we obtain that $\P H \cap \PSL_2(q_0)$ is a $\pi$-Hall
subgroup of $\PSL_2(q_0)$. Induction on $q$ implies that $\P H \cap
\PSL_2(q_0)$ satisfies  one of statements (a)--(d) of the lemma. So we obtain the statement of the lemma by induction if
$\P H
\subseteq \PSL_2(q_0)$. If $\P H \not\subseteq \PSL_2(q_0)$, then Lemma \ref{simpleepi}(e) implies  that $\P H\cap
\P\SL_2(q_0)$ satisfies
either (a) or (b) of the lemma, hence $\P H$
satisfies either (a) or (c) of the lemma.

If $\P H$ is included in a dihedral subgroup $D_{q-\varepsilon}$ of order $q-\varepsilon$, then
statement (a) of the lemma holds.

Assume that $\P H$ is included in $\Alt_4$. Since $2,3\in\pi$, we obtain that
$H=\Alt_4$. So statement (b) of the lemma holds in this case.

Assume that $\P H$ is included in $\Sym_4$. Since $2,3\in\pi$, we obtain that
$H=\Sym_4$. So statement (c) of the lemma holds
in this case.

Assume, finally,  that $\P H$ is included in $\Alt_5$. Since $2,3\in\pi$, by
Lemma \ref{HallSymmetric} we obtain that either $\P H=\Alt_5$, or
$\P H\simeq \Alt_4$. The case $\P H\simeq \Alt_4$ is considered above, so
we may assume that $\P H=\Alt_5$.  Thus
statement (d) of the lemma holds in this case.

(B) Now we prove that every subgroup $\P H$, satisfying one of statements (a)--(d) of the lemma, is a $\pi$-Hall
subgroup of $\P G$.
This fact is evident for statements (b)--(d), since $\vert \P G\vert_\pi=\frac{1}{2}(q^2-1)_\pi=\vert \P H\vert$. If
statement (a) holds,
i.e., $\P H$ is  a $\pi$-Hall subgroup  in a dihedral  subgroup $M$ of order $q-\varepsilon$ of $G$, then $\vert \P
G:M\vert=\frac{1}{2}q(q+\varepsilon)$ is a $\pi'$-number, whence~${\P H\in\Hall_\pi(\P G)}$.

(C) If a $\pi$-Hall subgroup $\P H$ satisfies statement (a) of
the
lemma, then $\P H$ has a Sylow tower of complexion $\prec$, where $\prec$ is the natural order, so all subgroups of this
type are conjugate by Lemma \ref{SylowTower}. If a $\pi$-Hall subgroup $\P
H$ satisfies statement (b) of the lemma, then $\P H$ has a Sylow tower of complexion $\prec$, where $3\prec 2$. Hence
all subgroups of this
type are conjugate by Lemma \ref{SylowTower}. In view of \cite[Chapter~XII]{Dick}, if either (c) or (d) holds, then $\P
G$ possesses
precisely two classes of conjugate subgroups isomorphic respectively to either $\Sym_4$ or $\Alt_5$, and these classes
are interchanged by
 $\PGL_2(q)$.  Notice that there can be more than one type of $\pi$-Hall subgroups in $G$, namely,  there can
exist subgroups, satisfying either (a) and (b), or (a) and (c), or (a) and (d). Hence, if $G\in E_\pi$, then
$k_\pi(G)\in\{1,2,3\}$.

(D) If $G$ possesses more than one class of conjugate $\pi$-Hall subgroups, then we have nothing to prove. So we may
assume that $G$
possesses
one class of conjugate $\pi$-Hall subgroups, i.e., either statement (a) or statement (b) of the lemma holds. Assume that
statement (a) of the lemma holds. Then $H$ is included in a dihedral subgroup $D_{2(q-\varepsilon)}$ of order
$2(q-\varepsilon)$, in particular $H$ possesses a normal abelian subgroup of index $2$, while its Sylow $2$-subgroup is
not normal. On the
other hand $G$ possesses a $\pi$-subgroup $K\simeq \SL_2(3)\simeq 2\nonsplitext 2^2\arbitraryext 3$. Clearly $K$ is not
isomorphic to a
subgroup of $H$. Assume that statement (c) of the lemma holds, i.e., $H\simeq \SL_2(3)$. Assume that $q\equiv \nu\pmod
3$, where $\nu=\pm1$,
then $G$ possesses a dihedral $\pi$-subgroup $K=D_{2(q-\nu)_{\mbox{}_{\{2,3\}}}}$ of order
$2(q-\nu)_{\mbox{}_{\{2,3\}}}$. Clearly $K$ is
not isomorphic to a
subgroup of~$H$.
\qed

As a corollary to Lemmas \ref{simpleepi}(e)  and \ref{SLU23dim2} we obtain the following lemma.

\begin{Lemma}\label{GLU23dim2}
Let $G=\GL^\eta_2(q)$, $\P G=G/Z(G)=\PGL_2^\eta(q)$, where $q$ is a power of a
prime $p$, and $\varepsilon=\varepsilon(q)$.  Let $\pi$ be a set of primes such that $2,3\in\pi$ and $p\not\in\pi$. A
subgroup  $H$ of $G$ is a $\pi$-Hall subgroup if and only if $H\cap \SL_2(q)$ is
a $\pi$-Hall subgroup of $\SL_2(q)$, $\vert H:H\cap\SL_2(q)\vert_\pi=(q-\eta)_\pi$, and either statement {\em (a)}, or
statement {\em (b)} of Lemma {\rm
\ref{SLU23dim2}} holds. More precisely, one of the following statements holds{\em:}
\begin{itemize}
\item[{\em (a)}]  $\pi\cap\pi(G)\subseteq \pi(q-\varepsilon)$, where
$\varepsilon=\varepsilon(q)$,
$\P H$ is a $\pi$-Hall subgroup in the dihedral group $D_{2(q-\varepsilon)}$ of order
$2(q-\varepsilon)$ of~$\P G$. All $\pi$-Hall subgroups of this type are
conjugate in~$G$.

\item[{\em (b)}] $\pi\cap\pi(G)=\{2,3\}$, $(q^2-1)_{\mbox{}_{\{2,3\}}}=24$, $\P H\simeq \Sym_4$.
All $\pi$-Hall subgroups of this type are conjugate in~$G$.
\end{itemize}
\end{Lemma}

\begin{Lemma}\label{OddIndexSp4}
Let $G=\mathrm{Sp}_4(q)$, $\varepsilon=\varepsilon(q)$, and $M$ be a maximal subgroup of odd index of $G$. Then one of
the following
statements hold{\em:}
\begin{itemize}
\item[{\em (a)}] $M\simeq \Sp_4(q^{\frac{1}{r}})$, where $q=p^\alpha$ and $r$ is an odd prime divisor of $\alpha$, and
all subgroups of this
type are conjugate in~$G${\em;}
\item[{\em (b)}] $M\simeq\mathrm{Sp}_2(q)\wr \Sym_2\simeq \mathrm{SL}_2(q)\wr \Sym_2$ and all subgroups of this type are
conjugate in~$G${\em;}
\item[{\em(c)}] $q\equiv\pm3\pmod8$,  $M\simeq2^{1+4}\arbitraryext \Omega_4^-(2)\simeq 2^{1+4}\arbitraryext
\mathrm{SL}_2(4)$ and all subgroups
of this type are conjugate in~$G$.
\end{itemize}
\end{Lemma}

\noindent{\slshape Proof.}\ \
In view of \cite{liesax} and Theorem~\ref{Maslova0}(a), (c), (h), (u) we obtain that either %\linebreak
$M=N_G(C_G(\sigma))$; or $M$
is the stabilizer of a nondegenerate subspace $U$ of dimension $m$ of $V$ and
$\psi(4) \gg \psi(m)$;   or $M$ is the stabilizer of an orthogonal decomposition $V=\bigoplus V_i$ into a direct sum of
isometric
subspaces $V_i$ of dimension $m=2^w\ge2$; or $q$ is a prime,  $q \equiv \pm 3 \pmod 8$, and $M\simeq
2^{1+4}\arbitraryext
\Omega_4^-(2)\simeq 2^{1+4}\arbitraryext \mathrm{SL}_2(4)$. In the first case \cite[Proposition~4.5.8]{KlLi} implies
statement (a) of the lemma. It is easy to see that the second case is impossible.
In the third case we obtain that $m=2$ and $M\simeq \Sp_2(q)\wr \Sym_2$. So \cite[Table~3.5.C]{KlLi} implies that  all
subgroups of this
type are conjugate in~$G$. In the
fourth case we obtain that $M$ satisfies (c) of the lemma and \cite[Table~3.5.C]{KlLi} implies that all
subgroups of this type are conjugate in~$G$.
\qed

\begin{Cor}\label{OddIndexO5}
Suppose $G=\Omega_5(q)$, $\varepsilon=\varepsilon(q)$, and let $M$ be a maximal odd-index subgroup of $G$.  Then one of
the following statements holds{\em:}

\begin{itemize}
\item[{\em (a)}] $M\simeq \Omega_5(q^{\frac{1}{r}})$, where $q=p^\alpha$ and $r$ is an odd prime divisor of $\alpha$,
and all subgroups of this
type are conjugate in~$G${\em;}
\item[{\em (b)}] $M\simeq\Omega_4^+(q)\arbitraryext 2$ and all subgroups of this type are conjugate in~$G${\em;}
\item[{\em (c)}] $q\equiv\pm3\pmod8$,  $M\simeq2^4\arbitraryext\Alt_5$ and all subgroups of this type are conjugate
in~$G$.
\end{itemize}
\end{Cor}

\noindent{\slshape Proof.}\ \
The corollary follows from known isomorphisms $\SL_2(q)\simeq \Sp_2(q)$,  $\mathrm{PSp}_4(q)\simeq \Omega_5(q)$, 
$\Omega_4^+(q)\simeq
\Sp_2(q)\circ \Sp_2(q)$, and Lemma~\ref{OddIndexSp4}.
\qed

\begin{Lemma}\label{OddIndexSL4}
Suppose $G\simeq \SL_4^\eta(q)$,  $\varepsilon=\varepsilon(q)$, and $M$ is a maximal subgroup of odd index of $G$. Then
one of the following statements holds{\em:}
\begin{itemize}
\item[{\em(a)}] $M\simeq \SL_4^\eta(q^{\frac{1}{r}})$, where $q=p^\alpha$ and $r$ is an odd prime divisor of $\alpha$,
and all subgroups
of this type are conjugate in~$G${\em;}
\item[{\em(b)}] $\eta=-\varepsilon$,\footnote{In this case $\vert Z(\SL_4^\eta(q))\vert=2$} $M\simeq
\Sp_4(q)\arbitraryext 2\simeq 2\arbitraryext
\SO_5(q)$, and there exist two classes of subgroup of this type, interchanged by~${\GL_4^\eta(q)}$;
\item[{\em(c)}] $M\simeq(\GL_2^\eta(q)\wr \Sym_2)\cap \SL_4^\eta(q)$, and all subgroups of this type are conjugate
in~$G${\em;}
\item[{\em(d)}] $\eta=\varepsilon$, $M\simeq (\GL_1^\eta(q)\wr\Sym_4)\cap \SL_4^\eta(q)$ and all subgroups of this type
are conjugate
in~$G${\em;}
\item[{\em(e)}] $\eta=\varepsilon$, $q\equiv 5\varepsilon\pmod8$, $M\simeq 4\arbitraryext 2^4\arbitraryext \Alt_6$, and 
there exist two
classes of subgroup of this type, interchanged by~${\GL_4^\eta(q)}$.
\end{itemize}
\end{Lemma}

\noindent{\slshape Proof.}\ \
In view of Theorem~\ref{Maslova0}(a), (b), (c), (f), (g), (q), (r), (s), (t) we obtain that $M$ satisfies one of the
following statements:
\begin{enumerate}
\item $M=N_G(C_G(\sigma))$, where $\sigma$ is a field automorphism of odd prime order of~$G$.
\item $M$ is the stabilizer of a nondegenerate (arbitrary if $\eta=+$) subspace  $U$ of dimension $m$ of $V$,
and
$\psi(4) \gg \psi(m)$.
\item $M$ is the stabilizer of an orthogonal (arbitrary if $\eta=+$)
decomposition $V=\bigoplus V_i$ into a direct sum of
isometric subspaces   $V_i$ of dimension $m=2^w\ge2$.
\item $M\simeq \Sp_4(q)\arbitraryext 2\simeq 2\arbitraryext
\SO_5(q)$, $q\equiv -\eta\pmod4$.
\item $M$ is the stabilizer of an orthogonal (arbitrary if $\eta=+$)
decomposition $V=\bigoplus V_i$ into a direct sum of
isometric subspaces   $V_i$ of dimension $1$, and $q\equiv\eta\pmod4$.
\item $M\simeq 4\arbitraryext 2^4\arbitraryext \Alt_6$, $q$ is  prime, and $q\equiv 5\eta\pmod8$.
\end{enumerate}
If $M$ satisfies the first statement, then \cite[Proposition~4.5.3]{KlLi} implies statement (a) of the lemma. The second
statement is
impossible. If $M$ satisfies
the third statement, then by using \cite[Proposition~4.2.9]{KlLi} we obtain statement (c) of the lemma. If $M$ satisfies
the fourth
statement, then  statement (b) of the lemma follows from \cite[Propositions~4.5.6 and~4.8.3]{KlLi}. If $M$
satisfies
the fifth statement, then by using \cite[Proposition~4.2.9]{KlLi} we obtain statement (d) of the lemma. Assume that $M$
satisfies
the sixth statement. Then \cite[Proposition~4.6.6]{KlLi} and the condition $q\equiv 5\eta \pmod8$ imply that
$\eta=\varepsilon$. Now
statement (e) of the lemma follows from~\cite[Proposition~4.6.6 and Tables~3.5.A and~3.5.B]{KlLi}.
\qed

\begin{Cor}\label{OddIndexO6}
Suppose $G=\mathrm{P}\Omega_6^\eta(q)$, $\varepsilon=\varepsilon(q)$, and $M$ is a maximal subgroup of odd index of $G$.
Then one of the
following statements holds{\em:}
\begin{itemize}
\item[{\em (a)}] $M\simeq \P\Omega_6^\eta(q^{\frac{1}{r}})$, where $q=p^\alpha$ and $r$ is an odd prime divisor of
$\alpha$, and all subgroups
of this type are conjugate in~$G${\em;}
\item[{\em (b)}] $\eta=-\varepsilon$, $M\simeq \Omega_5(q)\arbitraryext 2\simeq\mathrm{SO}_5(q)$, there exist two
classes of subgroups of this
type, and $\mathrm{PGO}_6^\eta(q)$ interchanges these classes{\em;}
\item[{\em (c)}] $\eta=-\varepsilon$, $M\simeq \left(\Omega_2^\eta(q)\times\Omega_4^+(q)\right)\arbitraryext
[4]$,\footnote{Here and
below, following
\cite{KlLi}, by $\mathrm{P}\Omega_2^\eta(q)$ we always mean a cyclic group of order $(q-\eta)/(4,q-\eta)$, while by
$\Omega_2^\eta(q)$ we always mean a cyclic group of order $(q-\eta)/(2,q-\eta)$.} and all subgroups of this type
are conjugate in~$G${\em;}
\item[{\em (d)}] $\eta=\varepsilon$, $M\simeq 2\arbitraryext
\left(\mathrm{P}\Omega_2^\eta(q)\times\mathrm{P}\Omega_4^+(q)\right)\arbitraryext
[4]$, and all subgroups of this type are conjugate in~$G${\em;}
\item[{\em (e)}] $\eta=\varepsilon$, $q\equiv5\varepsilon\pmod8$, $M\simeq 2^4\arbitraryext \Alt_6$, there exist two
classes of
subgroups of this type, these
classes are invariant under $\mathrm{PO}_6^\eta(q)$, and  $\mathrm{PGO}_6^\eta(q)$ interchanges these classes{\em;}
\item[{\em (f)}] $\eta=\varepsilon$, $M\simeq 2^2\arbitraryext \mathrm{P}\Omega_2^\eta(q)^3\arbitraryext
2^4\arbitraryext\Sym_3$,  and all
subgroups of this type are conjugate
in~$G$.
\end{itemize}
\end{Cor}

\noindent{\slshape Proof.}\ \
Follows from known isomorphism $\PSL_4^\eta(q)\simeq \P\Omega_6^\eta(q)$ and Lemma~\ref{OddIndexSL4}.
\qed

\section{Hall subgroups in linear, unitary, and symplectic groups}

\begin{Lemma}\label{mainstructureLinUniSym}
Let $\pi$ be a set of primes such that $2,3\in\pi$. Suppose $V$ is a linear, unitary, or symplectic space of dimension
$n$ with the
base field $\F_q$ of characteristic $p\not\in\pi$.  Assume that $G$ is chosen so that $\Omega(V)\leq G\leq I(V)$, and
$G$ possesses a
$\pi$-Hall subgroup $H$. Then one of the following statements holds{\em:}
\begin{itemize}
\item[{\em (a)}] $H$ stabilizes a decomposition $V=V_1\perp\ldots\perp V_k$ into a direct sum of pairwise orthogonal
nondegenerate
{\em(}arbitrary if $V$ is linear{\em)} subspaces $V_i$, and ${\dim(V_i)\le 2}$ for $i=1,\ldots,k$.
\item[{\em (b)}]  $V$ is a linear or a unitary space, $\dim(V)=4$, $I(V)=\GL^\eta(V)$,
${\vert \P G:\mathrm{PSL}^\eta(V)\vert}\le
2$,
$\pi\cap
\pi(G)=\{2,3,5\}$,
$q\equiv5\eta\pmod 8$ {\em(}in particular $\vert \PGL^\eta(V):\PSL^\eta(V)\vert=4$ and  $\P G\not=\PGL^\eta(V)${\em)},
$(q+\eta)_\ti=3$,
$(q^2+1)_{\mbox{}_5}=5$. Moreover $H\simeq 4\arbitraryext 2^4\arbitraryext \Sym_6$, if
$\vert \P G:\PSL^\eta(V)\vert=2$, and   $H\simeq 4\arbitraryext 2^4\arbitraryext \Alt_6$, if~${\P G=\PSL^\eta(V)}$.
\end{itemize}
\end{Lemma}

\noindent{\slshape Proof.}\ \
We proceed by induction on~$\dim(V)$. If $\dim(V)\le 2$ we have nothing to prove.

Assume that $\dim (V)>2$. Since $p\ne 2$, it follows that $G/Z(G)$ has a simple socle, and $G$ induces inner-diagonal
automorphisms on this socle. In
view of the main theorem from \cite{liesax} (see also Theorem \ref{Maslova0}) we obtain one of the following cases.
\begin{itemize}
\item[(1)] $V$ is unitary, $q=5$, $n=3$, and $H\cap \Omega(V)$ is included in $3\arbitraryext M_{10}\simeq
3\arbitraryext
\Alt_6 \arbitraryext 2$.
\item[(2)] $V$ is symplectic, $n=4$, and $H\leq M$, where $M$ is a maximal subgroup of odd index, satisfying
Lemma~\ref{OddIndexSp4}(c).
\item[(3)] $V$ is a linear or a unitary space, $n=4$, and   ${H\cap \Omega(V)\leq M}$, where $M$ is a maximal subgroup
of odd index in
$\Omega(V)$ satisfying Lemma~\ref{OddIndexSL4}(e).
\item[(4)] $H\leq M< I(V)$ for a group $M$ such that $\Omega(V_0)\leq M\leq I(V_0)$, where $I(V_0)$ is a group of the
same type as $I(V)$,
$\dim(V_0)=\dim(V)$ and the base field $\F_{q_{\mbox{}_0}}$  for~$V_0$ is a proper subfield of~$\F_q$.
\item[(5)] $V$ possesses a proper $H$-invariant nondegenerate (arbitrary if $V$ is linear) subspace~$U$.
\item[(6)] $H$ stabilizes a proper decomposition $V=U_1\perp\ldots\perp  U_m$
of $V$ into an orthogonal direct sum of pairwise isometric subspaces~$U_i$.
\end{itemize}

Now we proceed case by case.

(1) In this case $H\cap \Omega(V)$ is a $\pi$-Hall subgroup of $M\simeq 3\arbitraryext M_{10}\simeq3\arbitraryext \Alt_6
\arbitraryext 2$.
By Lemma \ref{simpleepi}(a) we obtain that $\Alt_6$ possesses a $\pi$-Hall subgroup. Lemma \ref{HallSymmetric} implies
that $\Alt_6$ does
not possess
a proper $\pi$-Hall subgroups with $2,3\in\pi$, hence $p=5\in \pi$, a contradiction  with~${p\not\in\pi}$.

(2) In this case $I(V)=\Omega(V)$ and $\vert \Omega(V):M\vert_\ti\ge 3$, a contradiction with $2,3\in\pi$.

(3) In this case $H\cap \Omega(V)$ is a $\pi$-Hall subgroup of $M\simeq 4\arbitraryext 2^4\arbitraryext \Alt_6$ and
$\eta=\varepsilon(q)$,
i.e., $(q-\eta)_\tw\ge4$. By using Lemmas \ref{simpleepi}(a) and \ref{HallSymmetric}, as in case (1) we obtain that
$\pi\cap \pi(G)
=\{2,3,5\}$ and $H\cap \Omega(V)=M$. Now $H$ is a $\pi$-Hall subgroup of $G$ if and only if $\pi(\vert
G:H\vert)\subseteq\pi'$. Since
$\vert\Omega(V):M\vert$ divides $\vert G:H\vert$, it follows that $H$ is a $\pi$-Hall subgroup of $G$ only if
$\pi(\vert\Omega(V):M\vert)$
is not divisible by $2$, $3$, and $5$, or, equivalently, only if $\vert \SL_4^\eta(q)\vert_{{\mbox{}_{\{2,3,5\}}}}=\vert
H\vert=2^9\cdot
3^2\cdot 5$.
The condition $p\not\in\pi$ implies that $p\not=2,3,5$. So $\vert
\SL_4^\eta(q)\vert_{{\mbox{}_{\{2,3,5\}}}}=\left(\left(q^2-1\right)\left(q^3-\eta\right)\left(q^4-1\right)\right)_{{
\mbox{}_{\{2,3,5\}}}}$.
By Lemma
\ref{OddIndexSL4}(e), we have $\vert \SL_4^\eta(q)\vert_\tw=\vert M\vert_\tw$ if and only if $q\equiv 5\eta\pmod8$.
Clearly,
$\left(\left(q^2-1\right)\left(q^3-\eta\right)\left(q^4-1\right)\right)_\ti=3^2$ if and only if $(q+\eta)_\ti=3$.
Finally,
$\left(\left(q^2-1\right)\left(q^3-\eta\right)\left(q^4-1\right)\right)_{{\mbox{}_5}}=5$ if and only if
$(q^2+1)_{{\mbox{}_5}}=5$.
Condition $\eta=\varepsilon(q)$ implies equality $\vert \P GL^\eta(V):\PSL^\eta(V)\vert=4$. By Lemma
\ref{OddIndexSL4}(e), it follows
that $\PSL^\eta(V)$ possesses  two classes of conjugate subgroups isomorphic to $L$ and  $\PGL^\eta(V)$ interchanges
these classes.   Since
$\PGL^\eta(V)/\PSL^\eta(V)$ is cyclic, Lemma
\ref{simpleepi}(e) implies that $\PGL^\eta(V)$ does not possesses a $\pi$-Hall subgroup $H$ such that $H\cap
\PSL^\eta(V)\simeq M$, while
every subgroup $G$ such that $\PSL^\eta(V)\leq G< \PGL^\eta(V)$ possesses a $\pi$-Hall subgroup $H$ such that
$H\cap\PSL^\eta(V)\simeq
M$. Thus statement (b) of the lemma holds in this case.

(4) We may assume that $q_0$ is the minimal possible number with $H\leq I(V_0)$. Since $\vert H\vert$ is coprime to $p$,
then $H$ is a
proper subgroup of $M \cap G$. Hence for $H$, $V_0$ either case (5) or case (6) holds. By using natural embeddings
$\F_{q_{\mbox{}_0}}\leq \F_q$ and
$V_0\leq V$ we obtain that for $H$, $V$ either case (5) or case (6) holds.

(5) In this case there exists a subspace $W\leq V$ such that $V=U\perp W$ and $W$ is $H$-in\-va\-ri\-ant (if
$I(V)=\GL(V)$, then the
existence
follows from Maschke Theorem, while in the remaining cases we can take $W=U^\perp=\{w\in W\mid (u,w)=0,\forall u\in
U\}$). Thus $H$ is
included in $G_0=G\cap (I(U)\times I(W))$ and $H$ is a $\pi$-Hall subgroup of~$G_0$. Denote by $\rho_U$ and $\rho_W$ the
projections from
$I(U)\times I(W)$ onto $I(U)$ and $I(W)$, respectively. Then $H\leq H^{\rho_U}\times H^{\rho_W}$, $\Omega(U)\leq
G_0^{\rho_U}\leq I(U)$, and
$\Omega(W)\leq G_0^{\rho_W}\leq I(W)$. Lemma \ref{simpleepi}(a) implies that $H^{\rho_U}$ and $H^{\rho_W}$ are
$\pi$-Hall subgroups of
$G_0^{\rho_U}$ and $G_0^{\rho_W}$, respectively. Furthermore, Lemma \ref{StructureInducedAutomorphisms} implies that
$G_0^{\rho_U}=I(U)$ and
$G_0^{\rho_W}=I(W)$. Hence $H^{\rho_U}$ in $G_0^{\rho_U}$ and $H^{\rho_W}$ in $G_0^{\rho_W}$ cannot satisfy statement
(b) of the lemma.
By induction
both
$U$ and $W$ have a decomposition into a direct sum of pairwise orthogonal nondegenerate
(arbitrary if $I(V)$ is linear) subspaces of dimensions at most $2$, whence statement (a) of the lemma holds.

(6) Since we have already considered case (5), we may assume that $H$ is irreducible. In this case $H$ is included in a
subgroup of type
$I(U_1)\wr \Sym_m$ of $I(V)$. In particular, $H$ normalizes the subgroup $G_0=G\cap (I(U_1)\times\ldots\times I(U_m))$
of $G$ and Lemma
\ref{simpleepi}(a) implies that $H_0=H\cap (I(U_1)\times\ldots\times I(U_m))$ is a $\pi$-Hall subgroup of $G_0$. Let
$N_1=\{x\in H\mid
U_1x=U_1\}$ be the stabilizer of $U_1$ in $H$. Clearly $H_0\leq N_1$. Denote by $\sigma$ the natural representation
$N_1\rightarrow I(U_1)$
of $N_1$. Assume also that $\rho:I(U_1)\times\ldots\times I(U_m)\rightarrow I(U_1)$ is the natural projection. We obtain
from definition
that the
restrictions of $\rho$ and $\sigma$ on $H_0$ coincide. Denote $G_0^\rho$ by $G_1$ and $H_0^\rho$ by $H_1$. Lemma
\ref{simpleepi}(a) implies that $H_1$ is a $\pi$-Hall subgroup of $G_1$. By using Lemma
\ref{StructureInducedAutomorphisms}, we also obtain
that $G_1=I(U_1)$. Thus $H_1$ in $G_1$ does not satisfy statement (b) of the lemma. Now $H_1=H_0^\rho=H_0^\sigma\leq
N_1^\sigma$. So
$N_1^\sigma$ is a $\pi$-Hall subgroup of $G_1$ and $N_1^\sigma=H_1$. By induction there exists an $N_1^\sigma$-invariant
decomposition
$U_1=W_{11}\perp\ldots\perp W_{1k}$ of $U_1$ into an orthogonal direct sum of nondegenerate (arbitrary if $V$ is linear)
subspaces of
dimensions at most $2$. By
definition this decomposition is $N_1$-invariant. Let $g_1=1,g_2,\ldots,g_m$ be a right transversal for the cosets of
$N_1$ in $H$. Since
$H$ is irreducible, without lost of generality  we may assume that $U_i=U_1g_i$. Now we set $W_{ij}=W_{1j}g_i$. Clearly
$H$ stabilizes the
decomposition $W_{11}\perp\ldots\perp W_{1k}\perp\ldots\perp W_{m1}\perp\ldots\perp W_{mk}$.
\qed

\begin{Lemma}\label{LU23.1}
Let $\pi$ be a set of primes with $2,3\in\pi$. Assume that $V$ is a linear or a unitary space with the base field $\F_q$
of characteristic
$p\not\in\pi$ and $G$ is chosen so that
$\SL^\eta(V)\leq G\leq \GL^\eta(V)$. Suppose also that
$H$ is a $\pi$-Hall subgroup of~$G$, and $H$ stabilizes a decomposition $V=V_1\perp\dots\perp V_m\perp
U_1\perp\dots\perp U_k$ into a
direct sum of pairwise orthogonal nondegenerate {\em(}arbitrary if $V$ is linear{\em)} subspaces such that  $\dim(
V_i)=2$ for $i=1,\dots, m$ and
$\dim (U_j)=1$ for
$j=1,\ldots,k$, and the decomposition cannot be refined. Then one of the following statements holds{\em:}
\begin{itemize}
\item[{\em (a)}]  $\dim (V)=2${\em;}
\item[{\em (b)}] $ q\equiv \eta \pmod {12}$ and $m=0${\em;}
\item[{\em (c)}] $ q\equiv -\eta \pmod {3}$ and $k\le 1${\em;}
\item[{\em (d)}] $ q\equiv -\eta \pmod {3}$, $ q\equiv \eta \pmod {4}$, $m$ is even, $k=3$, $(q+\eta)_\ti=3$, and
$m\not\equiv -1\pmod3${\em;}
\item[{\em (e)}] $q\equiv\eta\pmod{12}$, $m=1$, $k\equiv0\pmod3$, $k(k-1)\equiv 0\pmod4$.
\end{itemize}
\end{Lemma}

\noindent{\slshape Proof.}\ \
If $\dim(V)=2$ then statement (a) holds. Assume that $\dim(V)=2m+k>2$. Then $H$ is included in a subgroup $M$ of type
$(\GL_2^\eta(q)\wr
\Sym_m)\perp (\GL_1^\eta(q)\wr\Sym_k)$ of $\GL^\eta(V)$. Denote  by $L$ the intersection $M\cap G$. Since $H$ contains
a Sylow $2$-subgroup and a Sylow $3$-subgroup of $G$, the identities $\vert G:L\vert_\tw=\vert G:L\vert_\ti=1$ hold.

Since $H\cap \SL^\eta(V)$ is a $\pi$-Hall subgroup of $\SL^\eta(V)$, it is enough to prove statements (b)--(e) in case
$G=\SL^\eta(V)$.
In this case
$\vert M:L\vert= q-\eta$ and, by using Corollaries \ref{3part} and
\ref{2part}, we obtain the following identities
\begin{equation*}
|G|_\ti=\displaystyle\prod_{i=2}^{2m+k}(q^i-\eta^i)_\ti=
\left\{
\begin{array}{ll}
(q-\eta)^{2m+k-1}_\ti\big((2m+k)!\big)_\ti, & \text{ if } q\equiv\eta\pmod 3,\\
&\\
(q+\eta)_\ti^{m+[k/2]}\big(%\left[\displaystyle\frac{n}{2}\right]
(m+[k/2])!\big)_\ti, & \text{ if } q\equiv-\eta\pmod 3;
\end{array}
\right.
\end{equation*}
\begin{multline*}
|L|_\ti=\displaystyle\frac{1}{(q-\eta)_\ti}(q-\eta)_\ti^m(q^2-\eta^2)^m_\ti(m!)_\ti(q-\eta)^k_\ti
(k!)_\ti=\\
\left\{
\begin{array}{ll}
(q-\eta)^{2m+k-1}_\ti(m!)_\ti(k!)_\ti, & \text{ if } q\equiv\eta\pmod 3,\\
&\\
(q+\eta)_\ti^{m}(m!)_\ti(k!)_\ti, & \text{ if } q\equiv-\eta\pmod 3;
\end{array}
\right.
\end{multline*}
\begin{equation*}
|G|_\tw=\prod_{i=2}^{2m+k}(q^i-\eta^i)_\tw=(q-\eta)_\tw^{2m+k-1}(q+\eta)_\tw^{m+[k/2]}
\big((m+[k/2])!\big)_\tw;
\end{equation*}
\begin{multline*}
|L|_\tw=\displaystyle\frac{1}{(q-\eta)_\tw}(q-\eta)_\tw^m(q^2-\eta^2)^m_\tw(m!)_\tw(q-\eta)^k_\tw
(k!)_\tw\\ =(q-\eta)_\tw^{2m+k-1}(q+\eta)^m_\tw\big(m!k!\big)_\tw.
\end{multline*}

Assume first that $q\equiv\eta\pmod 3$. We have
$$|G:L|_\ti=\frac{\big((2m+k)!\big)_\ti}{(m!)_\ti (k!)_\ti}\ge \frac{\left(\left(2m\right)!\right)_\ti}{(m!)_\ti}=
(m+1)_\ti(m+2)_\ti\ldots (2m)_\ti.$$
Hence $m\le 1$ in this case. If $m=1$, then $$\vert G:L\vert_\ti=\frac{\left(\left(k+2\right)!\right)_\ti}{(k!)_\ti},$$
so
$\vert G:L\vert_\ti$ equals $1$ if and only if $(k+1)_\ti(k+2)_\ti$ equals $1$. Thus, if $q\equiv \eta\pmod3$, then the
equality $\vert
G:L\vert_\ti=1$ can be true only if either $m=0$, or $m=1$ and $k\equiv0\pmod3$. Furthermore
$$|G:L|_\tw=(q+\eta)_\tw^{[k/2]}\frac{\big(m+[k/2]!\big)_\tw}{(k!)_\tw}\ge(q+\eta)_\tw^{[k/2]}\frac{\big([k/2]!\big)_\tw
}{(k!)_\tw}.
$$
Now $$(k!)_\tw=2^{[k/2]+[k/2^2]+\dots}<2^{k/2+k/2^2+\dots}=2^k.$$ Since we have a strict inequality, it follows that
$(k!)_\tw\le 2^{k-1}$. Suppose $q\equiv-\eta\pmod 4$. Then  $(q+\eta)_\tw\ge 4$ and
$$1=|G:L|_\tw\ge 2^{2[k/2]}\frac{\big([k/2]!\big)_\tw}{2^{k-1}}%\ge \big([k/2]!\big)_\tw
,
$$
whence $k\le 1$. The case $m=1$ and $k=1$ is not possible, since $k\equiv0\pmod3$ for $m=1$. Therefore, if $q\equiv
\eta\pmod 3$, then
either
statement (a) of the lemma holds or $q\equiv\eta\pmod 4$ and $m\le 1$. If $q\equiv \eta\pmod 4$ and $m=0$, then we
obtain statement (b) of
the lemma. Assume that $q\equiv \eta\pmod4$ and $m=1$. Then Corollary \ref{2part} implies that $\vert
G\vert_\tw=(q-\eta)_\tw^{k+2}\cdot
((k+2)!)_\tw,$ and $\vert L\vert_\tw=(q-\eta)_\tw^{k+2}\cdot 2\cdot (k!)_\tw.$ Since  $\vert G:L\vert_\tw=1$, it follows
that
$\frac{(k+1)(k+2)}{2}\equiv1\pmod2$, hence $k(k-1)\equiv0\pmod4$,
and statement (e) of the lemma holds.

Now assume that $q\equiv-\eta\pmod 3$. Then
$$|G:L|_\ti=(q+\eta)_\ti^{[k/2]}\frac{\big((m+[k/2])!\big)_\ti}{(m!)_\ti (k!)_\ti}.$$
We have that
$$(k!)_\ti=3^{[k/3]+[k/3^2]+\dots}<3^{k/3+k/3^2+\dots}=3^{k/2},$$
hence  $(k!)_\ti\le 3^{[k/2]}$.
Now  $(q+\eta)^{[k/2]}_\ti\ge 3^{[k/2]}$, so
$$|G:L|_\ti\ge\frac{\big((m+[k/2])!\big)_\ti}{(m!)_\ti}.$$
Since  $|G:L|$ is not divisible by $3$, we have $[k/2]\le 2$, i.e., $k\le 5$. If
$k\le 1$, then statement (c) of the lemma holds. Suppose $k\in\{2,3,4,5\}.$  For $k=2$ we have
$$|G:L|_\ti=(q+\eta)_\ti\frac{\big((m+1)!\big)_\ti}{(m!)_\ti}\ge 3,$$
for $k=4,5$ we have
$$|G:L|_\ti=(q+\eta)_\ti^2\frac{\big((m+2)!\big)_\ti}{3(m!)_\ti}\ge 3,$$ which is impossible.
If  $k=3$, then $$\vert G:L\vert_\ti=\frac{(q+\eta)_\ti((m+1)!)_\ti}{3(m!)_\ti},$$ so $m\not\equiv -1\pmod 3$ and
$(q+\eta)_\ti=3$. Moreover
$$
|G:L|_\tw=(q+\eta)_\tw\frac{\big((m+1)!\big)_\tw}{2(m!)_\tw}=\frac{(q+\eta)_\tw}{2}(m+1)_\tw,
$$
so $m$ is even, $q\equiv\eta\pmod4$ and statement (d) holds.
\qed

\begin{Lemma} \label{HallSubgroupsOfLinearAndUnitaryGroups} Assume $G=\SL_n^\eta(q)$ is a special linear or unitary
group with the base
field  $\F_q$ of characteristic $p$ and $n\ge2$. Let $\pi$ be a set of primes such that $2,3\in\pi$ and $p\not\in\pi$.
Then the following statements
hold{\em:}
\begin{itemize}
\item[{\em (A)}] Suppose $G\in E_\pi$, and $H$ is a
$\pi$-Hall subgroup
of $G$. Then for $G$, $H$, and $\pi$ one of the following statements holds{\em:}
\begin{itemize}
\item[{\em (a)}] $n=2$ and one of the statements {\rm (a)--(d)}  of Lemma~{\rm\ref{SLU23dim2}} holds.

\item[{\em (b)}] either $q\equiv \eta \pmod {12}$, or $n=3$ and
 $q\equiv \eta \pmod 4$, $\Sym_n$ satisfies $E_\pi$, $\pi\cap\pi(G)\subseteq \pi(q-\eta)\cup \pi(n!)$ and if
$r\in(\pi\cap\pi(n!))\setminus
\pi(q-\eta)$, then $|G|_r=|\Sym_n|_r$. The subgroup $H$ is included in  $$M=L\cap G\simeq Z^{n-1}\arbitraryext
\Sym_n,$$ where $L=Z\wr
\Sym_n\leq \GL_n^\eta(q)$ and
$Z=\GL_1^\eta(q)$ is a cyclic group of order $q-\eta$. All $\pi$-Hall subgroups of this type are conjugate in~$G$.

\item[{\em (c)}] $n=2m+k$, where $k\in\{0,1\}$, $m\ge 1$,
$q\equiv -\eta \pmod {3}$, $\pi\cap\pi(G)\subseteq \pi(q^2-1)$, the groups $\Sym_m$ and $\GL_2^\eta(q)$ satisfy
$E_\pi$.\footnote{Notice
that,
in view of Lemma \ref{SLU23dim2}, the conditions $\GL_2^\eta(q)\in E_\pi$ and $q\equiv -\eta\pmod3$ imply that $q\equiv
-\eta\pmod{r}$ for
every odd $r\in\pi(q^2-1)\cap\pi$.}
The subgroup $H$ is included in $$M=L\cap
G\simeq(\underbrace{\GL_2^\eta(q)\circ\dots\circ\GL_2^\eta(q)}_{\displaystyle
m \mbox{ \rm times}})\arbitraryext\Sym_m\circ Z,$$ where $L=\GL_2^\eta(q)\wr \Sym_m\times Z\leq \GL_n^\eta(q)$ and $Z$
is a cyclic group of order $q-\eta$ for
$k=1$ and $Z$ is trivial for
$k=0$. The subgroup $H$ acting by conjugation on the set of factors in the central product 
\begin{equation}\label{prodgl2}
\underbrace{\GL_2^\eta(q)\circ\dots\circ\GL_2^\eta(q)}_{\displaystyle
m \mbox{ \rm times}}
\end{equation} has at most two orbits.
The intersection of $H$ with each factor $\GL_2^\eta(q)$ in \eqref{prodgl2} is a  $\pi$-Hall subgroup in
$\GL_2^\eta(q)$. The intersections with
the factors from the same orbit all satisfy the same statement {\em (a)} or {\em (b)} of Lemma~{\rm \ref{GLU23dim2}}.
Two $\pi$-Hall subgroups of $M$ are
conjugate in $G$ if
and only if they are conjugate in $M$. Moreover $M$ possesses one, two, or four classes of conjugate $\pi$-Hall
subgroups, while all
subgroups $M$
are conjugate in~$G$.

\item[{\em (d)}] $n=4$, $\pi\cap \pi(G)=\{2,3,5\}$,  $q\equiv5\eta\pmod 8$, $(q+\eta)_3=3$, and
$(q^2+1)_5=5$. The subgroup $H$ is isomorphic to $4\arbitraryext 2^4\arbitraryext \Alt_6$, $G$
possesses exactly two classes of conjugate $\pi$-Hall
subgroups of this
type and $\GL_4^\eta(q)$
interchanges these classes.

\item[{\em (e)}] $n=11$, $\pi\cap \pi(G)=\{2,3\}$, $(q^2-1)_{\mbox{}_{\{2,3\}}}=24$, $q\equiv -\eta\pmod3$, and
$q\equiv\eta\pmod4$. The subgroup  $H$ is
included in a
subgroup $M=L\cap G,$ where $L$ is a subgroup of $G$ of type $\left((\GL_2^\eta(q)\wr \Sym_4)\perp (\GL_1^\eta(q)\wr
\Sym_3)\right)$, and
$$H= \left(\left(\left(Z\circ
2\arbitraryext \Sym_4\right)\wr \Sym_4\right)\times
\left(Z\wr \Sym_3\right)\right)\cap G,$$ where $Z$ is a Sylow $2$-subgroup of a cyclic group of order $q-\eta$. All
$\pi$-Hall subgroups of
this type are conjugate in~$G$.
\end{itemize}
\item[{\em (B)}] Conversely, if the conditions on $\pi$, $n$, $\eta$, and $q$ in one of statements  {\em (a)--(e)} are
satisfied, then 
$G\in E_\pi$.
\item[{\em (C)}] If $G\in E_\pi$, then $k_\pi(G)\in\{1,2,3,4\}$.
\item[{\em (D)}] If $G\in E_\pi$ and $H\in\Hall_\pi(G)$, then $G$ possesses a $\pi$-subgroup  nonconjugate to a subgroup
of $H$, in
particular, $G$ does not satisfy~$D_\pi$.
\end{itemize}
\end{Lemma}

\noindent{\slshape Proof.}\ \
(A) Assume that $G\in E_\pi$. In view of Lemma \ref{SLU23dim2}, we may assume that $n>2$. By Lemmas
\ref{mainstructureLinUniSym} and
\ref{LU23.1}, it follows that either statement (d) of the
lemma holds, or every $\pi$-Hall subgroup $H$ of $G$ is included in $M$, where  $M$ is the stabilizer of an 
$H$-invariant decomposition
$V=V_1\perp\dots\perp V_m\perp U_1\perp\dots\perp U_k$ into a direct sum of pairwise orthogonal subspaces such that
$\dim (V_i)=2$ for
$i=1,\dots, m$ and $\dim (U_j)=1$ for $j=1,\dots, k$, and the decomposition cannot be refined. Below we denote by $T$
the intersection
$$(\GL^\eta(V_1)\times\ldots\times
\GL^\eta(V_m)\times \GL^\eta(U_1)\times\ldots \times \GL^\eta(U_k))\cap G$$ and by $\rho_i$ the natural projection
$T\rightarrow
\GL^\eta(V_i)$. Notice that Lemma \ref{StructureInducedAutomorphisms} implies the equality $T^{\rho_i}=\GL^\eta(V_i)$
for all $i$, since
$n>2$.
Moreover one of the statements (b)--(e) of Lemma \ref{LU23.1}
holds. Now we consider these statements separately.

Suppose statement (b) of Lemma \ref{LU23.1} holds, i.e., $m=0$ (so $k=n$) and $q\equiv \eta\pmod{12}$. In this case,
$M=T\arbitraryext
\Sym_n$, and $$T\simeq\underbrace{(q-\eta)\times\ldots\times(q-\eta)}_{(n-1)\text{ times}}.$$ Since $M$ includes a
$\pi$-Hall subgroup $H$
of $G$, we have $\vert G:M\vert_\pi=1$. If $r\in\pi(q-\eta)$, then Lemma \ref{FWair} (Corollary \ref{2part}, if $r=2$)
implies that
$\vert
G:M\vert_r=1$. If $r\not\in\pi(q-\eta)$, then $r\in \pi(\Sym_n)$ and the identity $\vert G:M\vert_r=1$ holds if and only
if $\vert
G\vert_r=\vert\Sym_n\vert_r$.  Thus for every $r\in \pi\cap \pi(G)\setminus\pi(q-\eta)$ the identity $\vert
G\vert_r=\vert\Sym_n\vert_r$
holds. Lemma \ref{simpleepi}(a) implies that $HT/T$ is a $\pi$-Hall subgroup of $M/T\simeq \Sym_n$, hence $\Sym_n\in
E_\pi$. Moreover by
\cite[Tables~3.5A and~3.5B]{KlLi} it follows that all subgroups of $G$ stabilizing a decomposition into direct
orthogonal sum of
$1$-dimension
subspaces are conjugate in $G$. By Lemma \ref{HallSymmetric}, we have $\Sym_n\in C_\pi$, hence, by Lemma
\ref{simpleepi}(f), $M\in C_\pi$
and all $\pi$-Hall subgroups of this type are conjugate in $G$. Thus statement  (b) of Lemma
\ref{HallSubgroupsOfLinearAndUnitaryGroups} holds.

Suppose statement (c) of Lemma \ref{LU23.1} holds, i.e., $k\le 1$ (so $m>0$) and $q\equiv -\eta\pmod3$. In this case,
$M=T\arbitraryext
\Sym_m$. Since $M$ includes a $\pi$-Hall subgroup $H$ of $G$, we have $\vert G:M\vert_\pi=1$, so $\pi\cap
\pi(G)\subseteq\pi(q^2-1)\cup
\pi(m!)$. Lemma \ref{SymmetricDivizors} implies that if $r\in \pi(m!)\setminus\pi(q^2-1)$, then $\vert G:M\vert_r>1$.
Hence $\pi\cap \pi(G)
\subseteq\pi(q^2-1)$. By Lemma \ref{simpleepi}(a), we obtain that $HT/T$ is a $\pi$-Hall subgroup of $M/T\simeq \Sym_m$,
while $H\cap T$ is
a $\pi$-Hall subgroup of $T$. Hence $(H\cap T)^{\rho_i}$ is a $\pi$-Hall subgroup of $\GL^\eta(V_i)\simeq \GL_2^\eta(q)$
and
$\GL_2^\eta(q)\in
E_\pi$. By Lemma \ref{simpleepi}(e), it follows that if $i,j$ are in the same $(HT/T)$-orbit, then $(H\cap T)^{\rho_i}$
in $\GL^\eta(V_i)$
and
$(H\cap T)^{\rho_j}$ in $\GL^\eta(V_j)$ satisfy  the same statement of Lemma \ref{GLU23dim2}. In view of Lemma
\ref{HallSymmetric}, $HT/T$
has at most two
orbits. Lemma \ref{GLU23dim2} implies that $\GL_2^\eta(q)$ possesses at most two classes of conjugate $\pi$-Hall
subgroups. By Lemma
\ref{Numberkpi(G)Inextension}, it follows that $M$
possesses one, two or four classes of conjugate $\pi$-Hall subgroups. Now we show that if, for some $g\in G$, the
$\pi$-Hall subgroups $H$
and
$H^g$ are in  $M$, then $g\in M$.  Assume  $\Delta$ is the set of subgroups $\SL^\eta(V_i)$. By \cite[Theorem~2]{Asc} it
follows that for every Sylow $2$-subgroup $S$ of $M$ (and, so of $G$) the identity
$\Delta={\rm Fun}(S)$ (in the notations of \cite{Asc}) holds. It follows that~${g\in N_G(\Delta)=M}$.

Suppose statement (d) of Lemma \ref{LU23.1} holds, i.e., $ q\equiv -\eta \pmod {3}$, $ q\equiv \eta \pmod {4}$, $m$ is
even, $k=3$,
$(q+\eta)_\ti=3$, and $m\not\equiv -1\pmod3$. Assume first that $m=0$. Then it is easy to see, that statement (b) of
Lemma \ref{HallSubgroupsOfLinearAndUnitaryGroups} holds in this case. Now assume that $m>0$. Then
$M=T\arbitraryext(\Sym_m\times \Sym_3)$,
and by Lemma \ref{simpleepi}(a), $HT/T$ is a $\pi$-Hall subgroup of $M/T$, while $H\cap T$ is a $\pi$-Hall subgroup of
$T$, and $(H\cap
T)^{\rho_i}$
is a $\pi$-Hall subgroup of $\GL^\eta(V_i)\simeq \GL_2^\eta(q)$. In particular, $\GL_2^\eta(q)\in E_\pi$. Since $q\equiv
\eta\pmod4$ and
$q\equiv -\eta\pmod 3$, Lemma \ref{GLU23dim2} implies that $\pi\cap \pi(\GL_2^\eta(q))=\{2,3\}$ and $\GL_2^\eta(q)\in
C_{\pi}$. By Lemma \ref{SymmetricDivizors}, it
follows that for every $r\in\pi(m!)\setminus\pi(q^2-1)$ the inequality $\vert G:M\vert_r>1$ holds. So $\pi\cap
\pi(G)=\{2,3\}$. Since
$HT/T$ is a $\pi$-Hall subgroup of $\Sym_m\times \Sym_3$, we obtain that $\Sym_m\in E_{\mbox{}_{\{2,3\}}}$. Lemma
\ref{HallSymmetric} and
the conditions
$m$ is even and $m\not\equiv -1\pmod3$ imply that $m=4$, in particular, $\Sym_m=\Sym_4$ is a $\{2,3\}$-gro\-up. Thus, by
Lemma
\ref{simpleepi}(f), we obtain that $M\in C_{\mbox{}_{\{2,3\}}}$, and statement (e) of the lemma holds.

Suppose statement (e) of Lemma \ref{LU23.1} holds, i.e., $q\equiv\eta\pmod{12}$, $m=1$, $k\equiv0\pmod3$, $k(k-1)\equiv
0\pmod4$. Again, by
Lemma \ref{simpleepi}(a), we obtain that $HT/T$ is a $\pi$-Hall subgroup of $M/T\simeq \Sym_k$, $H\cap T$ is a
$\pi$-Hall subgroup of $T$
and $(H\cap T)^{\rho_1}$ is a $\pi$-Hall subgroup of $\GL^\eta(V_1)\simeq \GL_2^\eta(q)$. So $(H\cap T)^{\rho_1}$
satisfies either statement
(a),
or statement (b) of Lemma \ref{GLU23dim2}. If $(H\cap T)^{\rho_1}$ satisfies statement (a) of Lemma \ref{GLU23dim2},
then $V_1$ possesses a
decomposition $V_1'\perp V_1''$ into an orthogonal sum of $1$-dimension subspaces and $(H\cap T)^{\rho_1}$ stabilizes
this decomposition.
Hence
$H$ stabilizes a decomposition $V_1'\perp V_1''\perp U_1\perp\ldots\perp U_k$ into an orthogonal sum of $1$-dimension
subspaces, i.e., $H$
satisfies statement (b) of Lemma \ref{LU23.1}. Since this case is already considered, we may assume that $(H\cap
T)^{\rho_1}$ satisfies
statement (b) of Lemma \ref{GLU23dim2}. In particular, $\pi\cap \pi(q^2-1)=\{2,3\}$ and
$(q^2-1)_{\mbox{}_{\{2,3\}}}=24$. Since
$\pi(M)=\pi(q^2-1)\cup\pi(k!)$ and $M$ contains a $\pi$-Hall subgroup of $G$, we obtain that for every $r\in \pi\cap
\pi(G)\setminus\{2,3\}$ the identity $\vert G\vert_r=\vert\Sym_k\vert_r$ holds. By Lemma \ref{simpleepi}(a), $HT/T$ is a
$\pi$-Hall
subgroup of $\Sym_k$. The conditions $k\equiv0\pmod3$ and $k(k-1)\equiv 0\pmod4$ imply that  $k\ge9$.  By Lemma
\ref{HallSymmetric}, it follows that $5\in\pi\cap\pi(G)$ (and $5\not\in\pi(q^2-1)$). By Lemma \ref{FWair}, we obtain
that $$\vert
G\vert_{\mbox{}_5}=(q^2+1)_{\mbox{}_5}^{\left[\frac{k+2}{4}\right]}\cdot
\left(\left(\frac{k+2}{4}\right)!\right)_{\mbox{}_5}.$$ Since
$$\vert
\Sym_k\vert_{\mbox{}_5}=5^{\left[\frac{k}{5}\right]+\left[\frac{k}{25}\right]+\ldots}<5^{\frac{k}{5}+\frac{k}{25}+\ldots
}=5^{\frac{k}{4}}\le
(q^2+1)_{\mbox{}_5}^{\left[\frac{k+2}{4}\right]}$$ (the last inequality is true, since $k(k-1)\equiv 0\pmod4$), we
obtain a contradiction
with $\vert
G\vert_r=\vert\Sym_k\vert_r$. Thus there does not exist a $\pi$-Hall subgroup satisfying statement (e) of
Lemma~\ref{LU23.1}.

(B) In view of Lemma \ref{SLU23dim2}, we may assume that $n>2$. If statement (d) holds, then it is easy to check that
$H$ is a $\pi$-Hall
subgroup of~$G$.

By \cite{KlLi}, it follows that there always exists a subgroup $M$ with the structure described in  statements (b), (c),
and (e) of the
lemma. Namely, $M$ is the stabilizer in $G$ of a decomposition
$$V=V_1\perp\dots\perp V_m\perp U_1\perp\dots\perp U_k$$ of the natural module  $V$ of $G$ into a direct orthogonal sum
of $m$
nondegenerate subspaces  $V_i$ of dimension 2 and  $k$ nondegenerate subspaces  $U_j$ of dimension  1. Moreover $m=0$ in
case (b), $k\le
1$ in case (c),  $m=4$ and
$k=3$ in case (e). By \cite[Tables~3.5A, 3.5B, and~3.5C]{KlLi}, it also follows that the stabilizers of decompositions
with the same numbers
$m$ and $k$ are conjugate in~$G$.

Assume that one of statements (b), (c), and (e) of the lemma holds. The subgroup $M$ of $G$ is obtained  as the
intersection
$L\cap G$, where $L$ is specified in the corresponding statements. In order to prove that $G\in E_\pi$, it is sufficient
to show that $\vert
\GL_n^\eta(q):L\vert_\pi=1$ and $L\in E_\pi$. Indeed, the identity $\vert \GL_n^\eta(q):L\vert_\pi=1$ implies that every
$H_1\in\Hall_\pi(L)$ is a $\pi$-Hall subgroup of $\GL_n^\eta(q)$, i.e., $\varnothing\not=\Hall_\pi(L)\subseteq
\Hall_\pi(\GL_n^\eta(q))$.
Since $G$ is normal in $\GL_n^\eta(q)$, Lemma \ref{simpleepi}(a) implies
that $H=H_1\cap G$ is a $\pi$-Hall subgroup of~$G$. Notice also that the identity $\vert\GL_n^\eta(q):L\vert_\pi=1$ and
the condition $L\in
E_\pi$ imply $\vert G:M\vert_\pi=1$ and~${M\in E_\pi}$. Now we consider  statements  (b), (c), and (e) separately.

Suppose statement (b) of the lemma holds. In this case, $L=\GL_1^\eta(q)\wr \Sym_n=(q-\eta)\wr \Sym_n$.  Assume first
that
$q\equiv\eta\pmod{12}$. Then, for every $r\in \pi(q-\eta)$, Lemma \ref{FWair} (Corollary \ref{2part} if $r=2$)
implies that $\vert \GL_n^\eta(q):L\vert_r=1$. If $r\in (\pi\cap \pi(G))\setminus \pi(q-\eta)$, then the condition
$\vert
G\vert_r=\vert\Sym_n\vert_r$ implies $\vert \GL_n^\eta(q):L\vert_r=1$. So, for every $r\in \pi$, the identity $\vert
\GL_n^\eta(q):L\vert_r=1$ holds. Hence $\vert \GL_n^\eta(q):L\vert_\pi=1$. Now assume that $n=3$, $q\equiv -\eta\pmod3$,
and $q\equiv
\eta\pmod4$. In this case, $3\not\in\pi(q-\eta)$ and the condition $\vert G\vert_\ti=\vert\Sym_3\vert_\ti$ implies
$\vert
\GL_3^\eta(q):L\vert_\ti=1$. For the remaining primes $r\in\pi\cap \pi(G)$, we have $r\in\pi(q-\eta)$ and Lemma
\ref{FWair} (Corollary
\ref{2part} if $r=2$) implies $\vert \GL_3^\eta(q):L\vert_r=1$. Hence $\vert \GL_3^\eta(q):L\vert_\pi=1$. Now we show
that $L\in E_\pi$.
Consider a $\pi$-Hall subgroup $\overline{H}_2$ of $\Sym_n$ (this subgroup exists, since $\Sym_n\in E_\pi$), let $H_2$
be its complete
preimage under the natural homomorphism $\varphi:L\rightarrow \Sym_n$. The kernel of $\varphi$ is an abelian subgroup
$$T=\underbrace{\GL_1^\eta(q)\times \ldots\times \GL_1^\eta(q)}_{n\text{ times}}=\underbrace{(q-\eta)\times \ldots\times
(q-\eta)}_{n\text{
times}},$$ and $H_2/T\simeq\overline{H}_2$ is a $\pi$-gro\-up. By Lemma \ref{simpleepi}(e) we obtain that $H_2\in
E_\pi$. Since $\vert
L:H_2\vert_\pi=\vert \Sym_n:\overline{H}_2\vert_\pi=1$, it follows that~${L\in E_\pi}$.

Suppose statement (c) of the lemma holds. In this case, $L=(\GL_2^\eta(q)\wr \Sym_m)\times Z$. By using Lemma
\ref{FWair} (Corollary
\ref{2part}, if $r=2$), we obtain that, for every $r\in\pi$, the identity $\vert \GL_n^\eta(q):L\vert_r=1$ holds. Hence
$\vert
\GL_n^\eta(q):L\vert_\pi=1$.
Now we show that $L\in E_\pi$. Consider a $\pi$-Hall subgroup $\overline{H}_2$ of $\Sym_m$ (this subgroup exists, since
$\Sym_m\in E_\pi$),
let $H_2$ be its complete preimage in $L$ under the natural homomorphism $\varphi:L\rightarrow \Sym_m$. The kernel of
$\varphi$ is the
subgroup $$T=\underbrace{\GL_2^\eta(q)\times\ldots\times \GL_2^\eta(q)}_{m\text{ times}}\times Z,$$ and
$H_2/T\simeq\overline{H}_2$ is a
$\pi$-group. Let $H_3$ be a $\pi$-Hall subgroup of $\GL_2^\eta(q)$ (this subgroup exists, since $\GL_2^\eta(q)\in
E_\pi$) and let $H_4$ be a
$\pi$-Hall subgroup of $Z$. Then $$H_5=\underbrace{H_3\times\ldots\times H_3}_{m\text{ times}}\times H_4$$ is a
$\pi$-Hall subgroup of $T$
and $H_5^T=H_5^{H_2}$. By Lemma \ref{simpleepi}(e), we obtain that $H_2\in E_\pi$. Since $\vert L:H_2\vert_\pi=\vert
\Sym_m:\overline{H}_2\vert_\pi=1$, it follows that~${L\in E_\pi}$.

Suppose finally that statement (e) of the lemma holds. Then $$L=(\GL_2^\eta(q)\wr\Sym_4)\times(\GL_1^\eta(q)\wr
\Sym_3),$$ $q\equiv
-\eta\pmod3$, $q\equiv \eta \pmod 4$, $\pi\cap \pi(G)=\{2,3\}$, and $n=11$. We have $$\vert
L\vert_{\mbox{}_{\{2,3\}}}=(q-\eta)^8_\tw\cdot
(q+\eta)_\ti^4
\cdot2^4\cdot 24\cdot (q-\eta)_\tw^3\cdot 2\cdot 3=(q-\eta)_\tw^{11}\cdot (q+\eta)_\ti^4\cdot 2^8\cdot 3^2.$$ By
Corollaries \ref{3part} and
\ref{2part}, we obtain that $$\vert \GL_{11}^\eta(q)\vert_{\mbox{}_{\{2,3\}}}=(q-\eta)_\tw^{11}\cdot (q+\eta)_\ti^5\cdot
 2^8\cdot 3.$$
The condition
$(q^2-1)_{\mbox{}_{\{2,3\}}}=24$ implies that $(q+\eta)_\ti=3$. Hence $\vert
\GL_{11}^\eta(q):L\vert_{\mbox{}_{\{2,3\}}}=1$. Now we show
that $L\in E_\pi$.
The group $L$ includes a subgroup $$H=\left(\left(\left(q-\eta\right)_{\mbox{}_{\{2,3\}}}\circ2\arbitraryext
\Sym_4\right)\wr \Sym_4\right)\times
\left(\left(q-\eta\right)_{\mbox{}_{\{2,3\}}}\wr \Sym_3\right),$$ which is a $\pi$-group. Since
$(q^2-1)_{\mbox{}_{\{2,3\}}}=24$ and $q\equiv-\eta\pmod3$,
we have $\vert
H\vert_{\mbox{}_{\{2,3\}}}=(4\cdot
24)^4\cdot 24\cdot 4^3\cdot 6=2^{30}\cdot 3^6,$ and $\vert L\vert_{\mbox{}_{\{2,3\}}}=(4\cdot 24)^4\cdot 24\cdot
4^3\cdot 6=2^{30}\cdot
3^6.$ Hence
$H$ is a $\pi$-Hall subgroup of~$L$.

Thus we proved that if one of the statements (a)--(e) of the lemma holds, then~${G\in E_\pi}$.

(C)
Note that if statement (e) of the lemma holds, then there also exist two classes of conjugate $\pi$-Hall subgroups of
$G$ satisfying
statement (c), while the remaining statements cannot be satisfied simultaneously. Therefore if $G\in E_\pi$, then
$k_\pi(G)\in\{1,2,3,4\}$.

(D) Given a $\pi$-Hall subgroup $H$ of
$G$, it remains to prove that there exists a $\pi$-sub\-gro\-up $K$ of $G$ such that $K$ is not conjugate to a subgroup
of $H$. If $n=2$,
then this statement is proven in Lemma \ref{SLU23dim2}(D). If statement (d) of the lemma holds, then $G$ possesses two
classes of
$\pi$-Hall subgroups. If statement (e) of the lemma holds, then, as we noted above, $k_\pi(G)=3$. So we may assume that
either statement
(b), or statement (c) of the lemma holds.
If
statement (b) holds then, in the proof of Lemma 6.1 from \cite{RevVdoContemp}, it is shown that $G$ possesses a
$\pi$-subgroup $K$ such that
$K$ is not isomorphic to a subgroup of $H$. Assume that statement (c) of the lemma holds.  $\Sym_4$ is known to have an
irreducible
representation of degree $3$ over a field of characteristic not equal to $2$ and $3$, while $\Alt_4$ is included in
$\SL^\eta_3(q)$ and
$\Alt_4$ does
not stabilize a decomposition into a sum of a $2$-dimensional and a $1$-dimensional subspaces.  Hence if statement (c)
holds, then we can
take $K$ to be equal to $\Alt_4$ acting  irreducibly on a $3$-dimensional subspace of the natural module of~$G$.
\qed

\begin{Lemma}\label{HallSubgroupsOfSymplecticGroups}
Let $G=\Sp_{2n}(q)$ be a symplectic group over a field $\F_q$ of characteristic $p$. Assume that  $\pi$ is a set of
primes
such that  $2,3\in\pi$ and $p\not\in\pi$. Then the following statements hold{\em:}
\begin{itemize}
\item[{\em (A)}] Suppose $G\in E_\pi$ and $H\in \Hall_\pi(G)$. Then both $\Sym_n$
and $\SL_2(q)$ satisfy $E_\pi$ and
$\pi\cap\pi(G)\subseteq \pi(q^2-1)$. Moreover, $H$ is a $\pi$-Hall subgroup of $$M=\Sp_2(q)\wr \Sym_n\simeq
\big(\underbrace{\SL_2(q)\times\dots\times\SL_2(q)}_{\displaystyle n
\mbox{ \rm times}}\big)\splitext \Sym_n\leq G.$$
\item[{\em (B)}] Conversely, if both $\Sym_n$ and $\SL_2(q)$ satisfy $E_\pi$ and $\pi\cap\pi(G)\subseteq
\pi(q^2-1)$, then $M\in E_\pi$ and every $\pi$-Hall subgroup $H$ of $M$ is a $\pi$-Hall subgroup of $G$.
\item[{\em (C)}] $\pi$-Hall
subgroups of $M$ are conjugate in $G$ if and only if they are conjugate in $M$, while
all such subgroups $M$ are conjugate in $G$. In particular, if $G\in E_\pi$, then $k_\pi(G)\in\{1,2,3,4,9\}$.
\item[{\em (D)}] If $G\in E_\pi$ and $H\in\Hall_\pi(G)$, then $G$ possesses a $\pi$-subgroup
nonconjugate to a subgroup of $H$, in  particular,~${G\not\in D_\pi}$.
\end{itemize}
\end{Lemma}

\noindent{\slshape Proof.}\ \
(A) Assume that $G\in E_\pi$ and $H$ is a $\pi$-Hall subgroup of $G$. The statement that $H$ is included in $M$ follows
from
Lemma~\ref{mainstructureLinUniSym}. By Lemma \ref{simpleepi}(a) we obtain that both $\Sym_n$ and $\SL_2(q)$ satisfy
$E_\pi$. Lemma
\ref{SymmetricDivizors} implies that for every $r\in\pi(n!)\setminus\pi(q^2-1)$ the inequality $\vert G:M\vert_r>1$
holds. So
$\pi\cap\pi(G)\subseteq \pi(q^2-1)$.

(B) Conversely, assume that  both $\Sym_n$ and $\SL_2(q)$ satisfy $E_\pi$ and $\pi\cap\pi(G)\subseteq \pi(q^2-1)$. In
view of
\cite[Table~3.5C]{KlLi} we obtain that $G$ possesses a subgroup $M\simeq \Sp_2(q)\wr\Sym_n$ and all such subgroups are
conjugate in $G$.
The condition $\pi\cap\pi(G)\subseteq\pi(q^2-1)$ and Corollaries \ref{FWair} and \ref{2part} imply that for every $r\in
\pi\cap\pi(G)$ the
equality
$\vert G:M\vert_r=1$ holds. Let $\overline{H}_1$ be a $\pi$-Hall subgroup of $\Sym_n$ and $H_1$ be its complete preimage
under the natural
homomorphism $M\rightarrow \Sym_n$. Let $H_2$ be a $\pi$-Hall subgroup of $\SL_2(q)$ and
$$H_3=\underbrace{H_2\times\ldots\times
H_2}_{n\text{ times}}\leq \underbrace{\SL_2(q)\times\ldots\times \SL_2(q)}_{n\text{ times}}=G_1\times \ldots\times
G_n=T\unlhd M.$$ Then
$H_3^T=H_3^{H_1}$, so Lemma \ref{simpleepi}(a) implies that $H_1$ possesses a $\pi$-Hall subgroup $H$. Since $\vert
M:H_1\vert_\pi=\vert\Sym_n:\overline{H}_1\vert_\pi=1$ and $\vert G:M\vert_\pi=1$, we obtain that $H$ is a $\pi$-Hall
subgroup of both $G$
and~$M$.

(C) Lemma \ref{HallSymmetric} (in the notations introduced in the proof of statement (B)) implies that $H$
has $t\le2$ orbits on the set of factors $\{G_1,\ldots,G_n\}$. By Lemma \ref{SLU23dim2} we obtain  that
$k_\pi(\SL_2(q))\in\{1,2,3\}$. Thus Lemma \ref{Numberkpi(G)Inextension}(b) implies that
$k_\pi(M)=k_\pi(\SL_2(q))^t\in\{1,2,3,4,9\}$. We
show
that if for some $g\in G$ the $\pi$-Hall subgroups $H$ and
$H^g$ are in  $M$, then $g\in M$, in particular, $k_\pi(G)=k_\pi(M)$.  Assume  $\Delta$ is the set of subgroups
$\SL_2(q)$ in $T$. By \cite[Theorem~2]{Asc}, it
follows that for every Sylow $2$-subgroup $S$ of $M$ (and, hence of $G$) the identity
$\Delta={\rm Fun}(S)$ (in the notations of \cite{Asc}) holds. It follows that~${g\in N_G(\Delta)=M}$.

(D) Given a $\pi$-Hall subgroup $H$ of
$G$ it remains to prove that there exists a $\pi$-sub\-gro\-up $K$ of $G$ such that $K$ is not conjugate to a subgroup
of $H$. If
$k_\pi(G)>1$ we have nothing to prove, so assume that $k_\pi(G)=1$. Lemma \ref{Numberkpi(G)Inextension} implies that
$k_\pi(\SL_2(q))=1$.
In the proof of Lemma \ref{SLU23dim2}(D), we have shown that if $k_\pi(\SL_2(q))=1$ then $\SL_2(q)$ possesses a
$\pi$-subgroup $X$ such
that
$X$ is not isomorphic to a subgroup of $H\cap G_1$. It is clear that $X\wr (HT/T)$ is a $\pi$-subgroup of $M$ (hence, of
$G$) and it is not
isomorphic to a subgroup of~$H$.
\qed

\section{Hall subgroups in orthogonal groups of dimension at most 6}

Recall that we denote  by $\Oo^\eta_n(q)$ the general
orthogonal group of degree $n$ and of sign
$\eta\in\{\circ,+,-\}$ over $\F_q$, while the symbol $\GO^\eta_n(q)$ denotes the group of similarities. Here $\circ$ is
an empty symbol, and
we use it only if $n$ is odd.

\begin{Lemma}\label{OrthogonalHallDim6}
Let $\pi$ be a set of primes with $2,3\in\pi$. Assume that $G=\Omega_n^\eta(q)$, where $q$ is a power of a prime
$p\not\in\pi$; and $H$ is a
$\pi$-Hall subgroup
of $G$. Suppose also that
$\varepsilon=\varepsilon(q)$. Then one of the following statements holds.
\begin{itemize}
\item[{\em (a)}] $n=2$, $G$ is cyclic of order $(q-\eta)/2$, $H$ is a unique $\pi$-Hall subgroup of~$G$.

\item[{\em (b)}]  $n=3$, $\pi\cap\pi(G)\subseteq \pi(q-\varepsilon)$, $H$ is a $\pi$-Hall subgroup in a dihedral
subgroup
$D_{q-\varepsilon}$ of order $q-\varepsilon$ of $G$. The subgroup $H$ stabilizes a decomposition $V_1\perp V_2$ of $V$
into an
orthogonal sum of subspaces of
dimension  $1$ and $2$, respectively, with $\eta(V_2)=\varepsilon$. All $\pi$-Hall subgroups of this type are conjugate
in~$G$.

\item[{\em (c)}] $n=3$, $\pi\cap\pi(G)=\{2,3\}$, $(q^2-1)_{\mbox{}_{\{2,3\}}}=24$, $H\simeq \Alt_4$.  The subgroup $H$
stabilizes the
decomposition
$V_1\perp V_2\perp V_3$ of $V$ into an orthogonal sum of $1$-dimension subspaces, and all $\pi$-Hall subgroups of this
type are conjugate
in~$G$.

\item[{\em (d)}] $n=3$, $\pi\cap\pi(G)=\{2,3\}$, $(q^2-1)_{\mbox{}_{\{2,3\}}}=48$, $H\simeq \Sym_4$.  The subgroup $H$
stabilizes the
decomposition
$V_1\perp V_2\perp
V_3$ of $V$ into an orthogonal sum of $1$-dimension subspaces,  there exist two classes of conjugate subgroups of this
type, and
$\mathrm{SO}_3(q)$ interchanges these classes.

\item[{\em (e)}] $n=3$, $\pi\cap\pi(G)=\{2,3,5\}$, $(q^2-1)_{\mbox{}_{\{2,3,5\}}}=120$, $H\simeq \Alt_5$.  The subgroup
$H$ is irreducible
and
primitive,
there exist two classes of conjugate subgroups of this type, and
$\mathrm{SO}_3(q)$ interchanges these classes.

\item[{\em (f)}]  $n=4$, $\eta=+$, $\pi\cap\pi(G)\subseteq \pi(q-\varepsilon)$, $H$ is a $\pi$-Hall subgroup in the
central product of two
subgroups isomorphic to $2\arbitraryext D_{q-\varepsilon}$.  The subgroup $H$ stabilizes  a decomposition
$V_1\perp V_2$ of $V$ into an orthogonal sum of $2$-dimensional subspaces with $\eta(V_1)=\eta(V_2)=\varepsilon$. All
$\pi$-Hall subgroups
of this type are conjugate in~$G$.

\item[{\em (g)}] $n=4$, $\eta=+$, $\pi\cap\pi(G)=\{2,3\}$, $(q^2-1)_{\mbox{}_{\{2,3\}}}=24$, $H\simeq
\SL_2(3)\circ \SL_2(3)$.
All $\pi$-Hall subgroups
of this type are conjugate in~$G$.

\item[{\em (h)}] $n=4$, $\eta=+$, $\pi\cap\pi(G)=\{2,3\}$, $(q-\varepsilon)_{\mbox{}_{\{2,3\}}}=12$,  $H\simeq
2\arbitraryext D_{12}\circ \SL_2(3)$. There exist two classes of conjugate
subgroups of this
type, $\mathrm{SO}_4^+(q)$ stabilizes these classes, while $\mathrm{O}_4^+(q)$ interchanges them.

\item[{\em (i)}] $n=4$, $\eta=+$, $\pi\cap\pi(G)=\{2,3\}$, $(q^2-1)_{\mbox{}_{\{2,3\}}}=48$, $H\simeq
2\arbitraryext\Sym_4\circ
2\arbitraryext\Sym_4$.
There exist four classes of conjugate subgroups of this type and every element of
$\mathrm{SO}_4^+(q)\setminus\Omega_4^+(q)$ induces an
involution with cyclic structure $(ij)(kl)$ on the set of these classes.

\item[{\em (j)}] $n=4$, $\eta=+$, $\pi\cap\pi(G)=\{2,3\}$, $(q-\varepsilon)_{\mbox{}_{\{2,3\}}}=24$,  $H\simeq
2\arbitraryext D_{24}\circ2\arbitraryext\Sym_4$. There exist four classes of conjugate subgroups of this type and each
element of
$\mathrm{O}_4^+(q)\setminus \Omega_4^+(q)$  induces an
involution with cyclic structure $(ij)(kl)$ on the set of these classes.

\item[{\em (k)}] $n=4$, $\eta=+$, $\pi\cap\pi(G)=\{2,3,5\}$, $(q^2-1)_{\mbox{}_{\{2,3,5\}}}=120$, $H\simeq
\mathrm{SL}_2(5)\circ\mathrm{SL}_2(5)$.
There exist four classes of conjugate subgroups of this type and every element of
$\mathrm{SO}_4^+(q)\setminus\Omega_4^+(q)$ induces an
involution with cyclic structure $(ij)(kl)$ on the set
of these classes.

\item[{\em (l)}] $n=4$, $\eta=+$, $\pi\cap\pi(G)=\{2,3,5\}$, $(q-\varepsilon)_{\mbox{}_{\{2,3,5\}}}=60$,  $H\simeq
2\arbitraryext D_{60}\circ\mathrm{SL}_2(5)$. There exist four classes of conjugate subgroups of this type and each
element of
$\mathrm{O}_4^+(q)\setminus \Omega_4^+(q)$  induces an
involution with cyclic structure $(ij)(kl)$ on the set of these classes.

\item[{\em (m)}]  $n=4$, $\eta=-$, $\pi\cap\pi(G)\subseteq \pi(q^2-1)$, $H$ is a $\pi$-Hall subgroup in a dihedral group
$D_{q^2-1}$. The
subgroup  $H$ stabilizes a decomposition
$V_1\perp V_2$ of $V$ into an orthogonal sum of $2$-dimensional subspaces with $\eta(V_1)=+$ and $\eta(V_2)=-$. All
$\pi$-Hall subgroups of this type are conjugate in~$G$.

\item[{\em (n)}] $n=4$, $\eta=-$, $\pi\cap\pi(G)=\{2,3\}$, $(q^2-1)_{\mbox{}_{\{2,3\}}}=24$, $H\simeq \Sym_4$. The
subgroup $H$ stabilizes a
decomposition $V_1\perp V_2\perp V_3\perp V_4$ of $V$ into an orthogonal sum of $1$-dimension subspaces. There exist two
classes
of conjugate $\pi$-Hall subgroups of this type,  invariant under~$\mathrm{O}_4^-(q)$, and
$\GO_4^-(q)$ interchanges these classes.

\item[{\em (o)}] $n=5$, $\pi\cap \pi(G)\subseteq \pi(q-\varepsilon)$, $H$ is isomorphic to an extension of a
$\pi$-Hall subgroup in the central product of two groups  isomorphic to  $2\arbitraryext
D_{q-\varepsilon}$ by a group of order $2$. The subgroup $H$ stabilizes a decomposition $V_1\perp V_2\perp V_3$ of $V$
with
$\dim(V_1)=1$, $\dim(V_2)=\dim(V_3)=2$, and $\eta(V_2)=\eta(V_3)=\varepsilon$.  All
$\pi$-Hall subgroups of this type are conjugate in~$G$.

\item[{\em (p)}] $n=5$, $\pi\cap\pi(G)=\{2,3\}$, $(q^2-1)_{\mbox{}_{\{2,3\}}}=24$, $H\simeq
(2\arbitraryext\Alt_4\circ2\arbitraryext\Alt_4)\arbitraryext 2$. The subgroup $H$  stabilizes a decomposition $V_1\perp
V_2$ of $V$ with
$\dim(V_1)=1$, $\dim(V_2)=4$, $\eta(V_2)=+$, and all $\pi$-Hall
subgroups of this type are conjugate in~$G$.

\item[{\em (q)}] $n=5$, $\pi\cap\pi(G)=\{2,3\}$, $(q^2-1)_{\mbox{}_{\{2,3\}}}=48$, $H\simeq (2\arbitraryext\Sym_4\circ
2\arbitraryext\Sym_4)\arbitraryext 2$. The subgroup $H$  stabilizes a decomposition $V_1\perp V_2$ of $V$ with 
$\dim(V_1)=1$,
$\dim(V_2)=4$, $\eta(V_2)=+$,  there exist two
classes of conjugate subgroups of this type, and  $\mathrm{SO}_5(q)$ interchanges these classes.

\item[{\em (r)}] $n=5$, $\pi\cap\pi(G)=\{2,3,5\}$, $(q^2-1)_{\mbox{}_{\{2,3,5\}}}=120$, $H\simeq (\mathrm{SL}_2(5)\circ
\mathrm{SL}_2(5))\arbitraryext
2$. The subgroup $H$  stabilizes a decomposition $V_1\perp V_2$ of $V$ with  $\dim(V_1)=1$, $\dim(V_2)=4$,
$\eta(V_2)=+$,  there exist two
classes of conjugate subgroups of this type, and  $\mathrm{SO}_5(q)$ interchanges these classes.

\item[{\em (s)}] $n=6$, $\eta=\varepsilon$,
$\pi\cap\pi(G)\subseteq \pi(q-\varepsilon)\cup\{3\}$ and, if $3\not\in\pi(q-\varepsilon)$, then $(q+\varepsilon)_\ti=3$,
$H$ is a $\pi$-Hall subgroup of a
solvable group  $\left((q-\varepsilon)^2\times (q-\varepsilon)/2\right)\arbitraryext \Sym_3$. The subgroup $H$
stabilizes a decomposition  $V_1\perp
V_2\perp V_3$ of $V$ into a sum of isometric $2$-dimensional subspaces with  $\eta(V_i)=\eta$ for $i=1,2,3$. All
$\pi$-Hall subgroups of this type are conjugate in~$G$.

\item[{\em (t)}] $n=6$, $\eta=-\varepsilon$,
$\pi\cap\pi(G)\subseteq \pi(q-\varepsilon)$, $H$ is a $\pi$-Hall subgroup of a solvable group
$\left((q^2-1)\times (q-\varepsilon)/2\right)\arbitraryext 2$. The subgroup $H$ stabilizes a decomposition  $V_1\perp
V_2\perp V_3$ of $V$ into a sum of $2$-dimensional subspaces with   $\eta(V_1)=\eta(V_2)=\varepsilon$,
$\eta(V_3)=-\varepsilon$. All
$\pi$-Hall subgroups of this type are conjugate in~$G$.

\item[{\em (u)}] $n=6$, $q\equiv -\eta\pmod{3}$, $\pi\cap\pi(G)=\{2,3\}$, $(q^2-1)_{\mbox{}_{\{2,3\}}}=24$,  $H$ is
isomorphic to
$\left((q-\eta)_2\circ 2\arbitraryext\Sym_4\circ2\arbitraryext\Sym_4\right)\arbitraryext 2$. The subgroup $H$ stabilizes
a decomposition
$V_1\perp V_2$ of $V$ with $\dim(V_1)=2$, $\dim(V_2)=4$, $\eta(V_1)=\eta$, $\eta(V_2)=+$. All
$\pi$-Hall subgroups of this type are conjugate in~$G$.

\item[{\em (v)}] $n=6$, $\eta=\varepsilon$,
$q\equiv\pm3\pmod8$, $\pi\cap\pi(G)=\{2,3,5\}$, $(q^2-1)_{\mbox{}_{\{2,3\}}}=24$, $q^2\equiv-1\pmod 5$,
$q\equiv-\eta\pmod 3$, $H\simeq
2^{5}\arbitraryext \Alt_6$. The subgroup $H$ stabilizes decomposition $V=V_1\perp V_2\perp V_3\perp V_4\perp V_5\perp
V_6$ of $V$ into a sum
of $1$-dimensional subspaces. There exist two classes of $\pi$-Hall subgroups of this type, invariant under 
$\mathrm{O}_6^\eta(q)$, and
$\mathrm{GO}_6^\eta(q)$ interchanges these classes.
\end{itemize}
\end{Lemma}

\noindent{\slshape Proof.}\ \
In the proof we use the known isomorphisms of orthogonal groups of small dimension and other groups  (see
\cite[Proposition~2.9.1]{KlLi},
for
example). Without further references we also use the fact that $q$ is odd.

If $n=2$, then
the isomorphism
$\Omega_2^\eta(q)\simeq (q-\eta)/2$ implies statement (a) of the lemma.

If $n=3$, then the isomorphisms  $\Omega_3(q)\simeq \mathrm{PSL}_2(q)$ and
$\mathrm{SO}_3(q)\simeq
\mathrm{PGL}_2(q)$ together with  Lemma~\ref{SLU23dim2} imply statements (b)--(e) (except for the decomposition part).
Indeed, if $H$ is a
$\pi$-Hall subgroup of $\Omega_3(q)\simeq \PSL_2(q)$, then one of the following holds:
\begin{enumerate}
\item $\pi$ and $H$ satisfy Lemma \ref{SLU23dim2}(a), so $\pi\cap\pi(G)\subseteq \pi(q-\varepsilon)$ and $H$ is a
$\pi$-Hall subgroup in
the dihedral subgroup $D_{q-\varepsilon}$ of $G$. Moreover all $\pi$-Hall subgroups of this type are conjugate in $G$.
So statement (b) of
the lemma holds (except for the decomposition part).
\item $\pi$ and $H$ satisfy Lemma \ref{SLU23dim2}(b), so $\pi\cap\pi(G)=\{2,3\}$, $(q^2-1)_{\mbox{}_{\{2,3\}}}=24$, and
$H\simeq \Alt_4$.
Moreover all $\pi$-Hall subgroups of this type are conjugate in $G$. So statement (c) of
the lemma holds (except for the decomposition part). 
\item $\pi$ and $H$ satisfy Lemma \ref{SLU23dim2}(c), so $\pi\cap\pi(G)=\{2,3\}$, $(q^2-1)_{\mbox{}_{\{2,3\}}}=48$, and
$H\simeq \Sym_4$.
Moreover there exist exactly two classes of conjugate subgroups of this type, and $\SO_3(q)\simeq \PGL_2(q)$
interchanges these classes.  So
statement (d) of
the lemma holds (except for the decomposition part). 
\item $\pi$ and $H$ satisfy Lemma \ref{SLU23dim2}(d), so $\pi\cap\pi(G)=\{2,3,5\}$, $(q^2-1)_{\mbox{}_{\{2,3\}}}=120$,
and $H\simeq \Alt_5$.
Moreover there exist exactly two classes of conjugate subgroups of this type, and $\SO_3(q)\simeq \PGL_2(q)$
interchanges these classes.  So
statement (e) of
the lemma holds (except for the decomposition part). 
\end{enumerate}
Now we show the statements about decomposition in statements  (b)--(e) of the lemma.

There exists a subgroup
$L=\left(\mathrm{O}_1(q)\times \mathrm{O}^\varepsilon_2(q)\right)\cap \Omega_3(q)$.  By Lemma
\ref{StructureInducedAutomorphisms},  the
subgroup $L$ is isomorphic to a dihedral group of order  $q-\varepsilon$.  Since all $\pi$-Hall subgroups satisfying
statement (b) of the
lemma are conjugate, it follows that  $H\leq L$, whence the decomposition in statement (b) of the lemma.

In view of \cite[Proposition~4.2.15]{KlLi}, there exists a subgroup of type $\Oo_1(q)\wr S_3$ of $\Omega_3(q)$, and it
is isomorphic to
$\Alt_4$ under the conditions of statement (c) of the lemma (in this case all subgroups of this type are conjugate), and
to $\Sym_4$ under
the conditions of statement (d) of the lemma (in this case there exists two classes of subgroups of this type).  It
follows that there exist
decompositions of $V$ as in statements (c) and~(d).

Suppose that $H\simeq\Alt_5$. The minimal degree of faithful complex representation of  $\Alt_5$ equals $3$. Since
$(q,\vert H\vert)=1$,
we obtain that the minimal
degree of a nontrivial representation of $\Alt_5$ over an algebraically closed field of characteristic $p$ equals $3$ as
well. Hence
the subgroup $\Alt_5$ in statement (e) is absolutely irreducible and primitive.

Suppose $n=4$ and $\eta=+$. Then   $\Omega_4^+(q)$ is isomorphic to  $\mathrm{SL}_2(q)\circ\mathrm{SL}_2(q)$,
$\mathrm{SO}_4^+(q)$ normalizes each
factor and induces on it the group of all in\-ner-di\-a\-go\-nal automorphisms, while each element of 
$\mathrm{O}_4^+(q)\setminus
\mathrm{SO}_4^+(q)$ interchanges the factors. Thus statements (f)--(l) of the lemma (except for the decomposition of $V$
in statement (f))
follows from Lemma~\ref{SLU23dim2}. A $\pi$-Hall subgroup  $H$ in statement (f), up to conjugation, is included in
$\left(\mathrm{O}_2^\varepsilon(q)\times\mathrm{O}_2^\varepsilon(q)\right)\cap \Omega_4^+(q)$, whence the decomposition
in statement~(f).

Suppose $n=4$ and $\eta=-$. We have that $\Omega_4^-(q)$ is isomorphic to   $\mathrm{PSL}_2(q^2)$, while  
$\mathrm{O}_4^-(q)$ is isomorphic
to
$2\times \mathrm{PSL}_2(q^2)\arbitraryext 2$ and induces on  $\mathrm{PSL}_2(q^2)$ the subgroup of $\Aut(\PSL_2(q^2))$
generated by the 
inner
automorphisms and a field automorphism of order $2$ \cite[Proposition~2.9.1(v) and proof]{KlLi}. Since
$((q^2)^2-1)_{\mbox{}_{\{2,3\}}}$ divides $48$, it follows that $H\cap \PSL_2(q^2)$ is a $\pi$-Hall subgroup of
$\PSL_2(q^2)$ satisfying
either (a) or (c) of Lemma~\ref{SLU23dim2}. If   $H\cap \PSL_2(q^2)$ satisfies Lemma~\ref{SLU23dim2}(a), i.e., $H\cap
\PSL_2(q^2)$ is a
$\pi$-Hall subgroup in $D_{q^2-1}$, then $(H\cap \PSL_2(q^2))^{\PSL_2(q^2)}=(H\cap \PSL_2(q^2))^{\Aut(\PSL_2(q^2))}$, so
by Lemma
\ref{simpleepi}(e) we obtain that statement (m) of the lemma (except the decomposition) holds. If $H\cap \PSL_2(q^2)$
satisfies
Lemma~\ref{SLU23dim2}(c), i.e., $H\cap \PSL_2(q^2)\simeq \Sym_4$, then $\PSL_2(q^2)$ possesses two classes of conjugate
subgroups of this
type. Moreover $H\cap \PSL_2(q^2)$ is a subgroup of symplectic type (type $\mathcal{C}_6$ in terms of \cite{KlLi}) of
$\PSL_2(q^2)$ and
\cite[Tables~3.5A and~3.5.G]{KlLi} imply that a field automorphism leaves these classes invariant. So Lemma
\ref{simpleepi}(e) implies
statement (n) of the lemma (except the decomposition). A subgroup $H$ in statement
(m) is included in
$\left(\mathrm{O}_2^+(q)\times\mathrm{O}_2^-(q)\right)\cap \Omega_4^-(q)$, whence the decomposition in this statement. A
subgroup $H$ in
statement (n) is included in $\left(\mathrm{O}_1(q)\times\left(\mathrm{O}_1(q)\wr \Sym_3\right)\right)\cap
\Omega_4^-(q)$, whence the
decomposition in statement~(n).

Suppose $n=5$, then   $\Omega_5(q)\simeq \mathrm{PSp}_4(q)$. By Lemma \ref{HallSubgroupsOfSymplecticGroups}, it follows
that  $M$
is included in
\begin{multline*}L=(\mathrm{O}_1(q)\times
\mathrm{O}_4^+(q))\cap\Omega_5(q)\simeq
\left(\mathrm{Sp}_2(q)\circ\mathrm{Sp}_2(q) \right)\arbitraryext \Sym_2\\ \simeq
\left(\mathrm{SL}_2(q)\circ\mathrm{SL}_2(q) \right)\arbitraryext \Sym_2
\end{multline*} and statements (o)--(r) (except the statements about classes in  (q) and (r)) of the lemma follows from
statements
(f)--(l) of the lemma. It is easy to see, that two classes of conjugate $\pi$-Hall subgroups in statements (i) and (k)
are invariant
under a graph automorphism, while the remaining two classes are interchanged by this automorphism. Hence the statement
on the number of
classes follows from Lemmas~\ref{simpleepi}(e) and~\ref{StructureInducedAutomorphisms}.

Suppose $n=6$. Then  $\mathrm{P}\Omega_6^\eta(q)\simeq \mathrm{PSL}_4^\eta(q)$. Statements (s) and (t) (except
the
decomposition of $V$) follow from Lemma~\ref{HallSubgroupsOfLinearAndUnitaryGroups}(b). Statement (u) (except the
decomposition of $V$)
follows from Lemma~\ref{HallSubgroupsOfLinearAndUnitaryGroups}(c). Statement (v) (except the decomposition of $V$)
follows from
Lemma~\ref{HallSubgroupsOfLinearAndUnitaryGroups}(d). A subgroup  $H$, satisfying statement (s), is included in
$\left(\mathrm{O}_2^\eta(q)\wr \Sym_3\right)\cap \Omega_6^\eta(q)$, whence the decomposition of $V$ in statement (s). A
subgroup
$H$, satisfying statement (t), is included in $\left(\mathrm{O}_2^\varepsilon(q)\times\mathrm{O}_2^\varepsilon(q)\times
\mathrm{O}_2^{-\varepsilon}(q)\right)\cap \Omega_6^\eta(q)$, whence the decomposition of $V$ in statement  (t). A
subgroup $H$,
satisfying  (u), is included in $\left(\mathrm{O}_2^\eta(q)\times
\mathrm{O}_4^{+}(q)\right)\cap \Omega_6^\eta(q)$, whence the decomposition of $V$ in statement (u). Finally, by
\cite[Proposition~4.2.15]{KlLi} it follows that $H$ satisfying  (v) exists only if  $\eta=\varepsilon$, and $H$ is
included in
$\left(\mathrm{O}_1(q)\wr S_6\right)\cap \Omega_6^\eta(q)$, whence the decomposition of $V$ in statement~(v).
\qed

\section{Hall subgroups in orthogonal groups}

In this section we say that a hypothesis $(\star)$ is true, if the following statements hold:

\begin{enumerate}
 \item $V$ is a vector space with a nondegenerate symmetric bilinear form over a field  $\F_q$ of odd characteristic
$p$;
\item $\varepsilon=\varepsilon(q)$;
\item $\pi$ is a set of primes such that  $2,3\in\pi$, and $p\not\in\pi$;
\item $\Omega(V)\leq G\leq I(V)$;
\item $H$ is a $\pi$-Hall subgroup of $G$;
\item $V$, $G$, and $H$ do not satisfy  the following statements:
\begin{enumerate}
\item $\dim (V)=7$, $G=\Omega(V)$, $H\simeq \Omega_7(2)$;
\item $\dim (V)=8$, $\eta(V)=+$, $G=\Omega(V)$, $H\simeq 2\arbitraryext \Omega^+_8(2)$;
\item $\dim (V)=9$, $G=\Omega(V)$, $H$ stabilizes a decomposition $V=U\perp W$, $\dim(U)=8$, $\eta(U)=+$, $\dim(W)=1$,
$H\leq (I(U)\times
I(W))\cap G$ and the projection of $H$ on $I(U)$ is isomorphic to $2\arbitraryext{\Omega_8^+(2)\arbitraryext 2}$.
\end{enumerate}
\end{enumerate}

\begin{Lemma} \label{orthored} Assume that $(\star)$ holds. Then there exists an $H$-in\-va\-ri\-ant decomposition
$$V=V_1\perp\ldots\perp V_l$$ of $V$ into an orthogonal sum of nondegenerate subspaces such that $\dim(V_i)\le 4$ for
all~${i=1,\ldots,l}$.
\end{Lemma}

\noindent{\slshape Proof.}\ \
If $n\le 4$, we have nothing to prove. If $n=5$ or $n=6$, then the claim follows from Lemma
\ref{OrthogonalHallDim6}. So we may assume $n\ge 7$. Since $p\not\in\pi$, we obtain that $H\not=G$. Hence $H$ is
included in a
maximal subgroup $M$ of odd index of $G$. In view of \cite{liesax}  we obtain that one of the following statements
holds:
\begin{itemize}
\item[(1)] $\Omega(V_0)\leq M\leq I(V_0)$, where $I(V_0)$ is a group of the same type as $I(V)$, $\dim(V_0)=\dim(V)$,
and the base field
$\F_{q_0}$ of $V_0$ is a proper subfield of $\F_q$.
\item[(2)] $\dim(V)=7$,  $M\cap\Omega(V)=\Omega_7(2)$, $q$ is prime and $q\equiv\pm3\pmod8$.
\item[(3)] $\dim(V)=8$, $\eta=+$, $M\cap\Omega(V)=2\arbitraryext \Omega_8^+(2)$, $q$ is prime and $q\equiv\pm3\pmod8$.
\item[(4)] $M$ is the stabilizer of a non-singular subspace.
\item[(5)] $M$ is the stabilizer of a decomposition $V=V_1\perp\ldots\perp V_k$ with all $V_i$ isometric.
\end{itemize}
Now we consider all statements separately.

(1) As in the proof of Lemma \ref{mainstructureLinUniSym} we obtain the statement of the lemma by induction on~$q$.

(2), (3) By \cite[Theorem~1.2]{RevinSIBAM} it follows that both $\Omega_7(2)$ and $\Omega_8^+(2)$ do not possess proper
$\pi$-Hall subgroups
with $2,3\in\pi$. Hence $H\cap \Omega(V)=M\cap \Omega(V)$, i.e.,  $(\star)$ is not satisfied, a contradiction.

(4) In this case $M$ stabilizes a decomposition $U\perp W$ of $V$, hence $M=(I(U)\times I(W))\cap G$. Without lost of
generality we may
assume that $\dim(U)=k\ge\dim(W)=m$. Denote by $\rho_U$ and $\rho_W$ the natural projections of $M$ into $I(U)$ and
$I(W)$, respectively. By
Lemma \ref{simpleepi}(a) it follows that $H^{\rho_U}\in\Hall_\pi(M^{\rho_U})$ and  $H^{\rho_W}\in\Hall_\pi(M^{\rho_W})$.
If $m$ is greater
than $1$, then Lemma \ref{StructureInducedAutomorphisms} implies that $M^{\rho_U}=I(U)$ and $M^{\rho_W}=I(W)$. Hence
Lemma
\ref{HallAmostSimpleInduce}
implies that for $U$, $M^{\rho_U}$, $H^{\rho_U}$ and $W$, $M^{\rho_W}$, $H^{\rho_W}$ condition  $(\star)$ holds and we
obtain the statement
of the
lemma by induction on dimension. If $m=1$ and $k\not\in\{7,8,9\}$, then again for $U$,  $M^{\rho_U}$, $H^{\rho_U}$ and
$W$, $M^{\rho_W}$,
$H^{\rho_W}$
hypothesis $(\star)$
holds and the statement of the lemma follows by induction. If $m=1$ and $k\in\{7,8,9\}$, then Theorem \ref{Maslova0}(e)
implies that $k=8$
and $n=9$. In this case either condition  $(\star)$ is not true, or we obtain the statement of the lemma by induction.

(5) Since we already considered statement (4), we may assume that $M$ is irreducible. Then $$M=\left(\left(
I(V_1)\times\ldots\times
I(V_k)\right)\splitext \Sym_k\right)\cap G$$ and, by Lemma \ref{StructureInducedAutomorphisms}, we have
$M^{\rho_{V_i}}=I(V_i)$ for
$i=1,\ldots,k$. By using Lemma \ref{HallAmostSimpleInduce}, we obtain the statement of the lemma as in the proof of
Lemma~\ref{mainstructureLinUniSym}.
\qed

Assume that   $(\star)$ holds. By Lemma \ref{orthored}, there exists an $H$-in\-va\-ri\-ant decomposition
\begin{equation}\label{decomorth}
V=V_1\perp\dots\perp V_l,
\end{equation}
which cannot be refined. All subspaces $V_i$ in decomposition \eqref{decomorth} are nondegenerate and $\dim(V_i)\le 4$.
 For $k\in\{1,2,3,4\}$ and $\delta\in\{+,-,\circ\}$, we denote by  $V(k,\delta)$ the sum of
$V_i$ with  $\dim (V_i)=k$ and $\eta(V_i)=\delta$, and by $d(k,\delta)$ the number of such subspaces $V_i$. It is clear
that
$\dim(V(k,\delta))=kd(k,\delta)$, subspace $V(k,\delta)$ is $H$-in\-va\-ri\-ant for every $k$ and $\delta$, and
\begin{multline}\label{orthodec}
V=\sum_{(k,\delta)}V(k,\delta)=\\ V(1,\circ)\perp V(2,\varepsilon)\perp V(2,-\varepsilon)\perp V(3,\circ)\perp
V(4,+)\perp
V(4,-).\end{multline}
For brevity we write
$\eta(k,\delta)$ and $D(k,\delta)$ instead of $\eta(V(k,\delta))$ and $D(V(k,\delta))$, respectively. If $k$ is even
then
${\eta(k,\delta)=\delta^{d(k,\delta)}}$. Set also $I(k,\delta)=I(V(k,\delta))$, $\Omega(k,\delta)=\Omega(V(k,\delta))$.
Thus $H$ is
included in a subgroup $G_0$ of  $G$ of type
$$G_0=(I(1,\circ)\times I(2,\varepsilon)\times I(2,-\varepsilon)\times I(3,\circ)\times I(4,+)\times I(4,-))\cap G.
$$
Denote by $G(k,\delta)$ and $H(k,\delta)$ the projection of $G_0$ and $H$ on $I(k,\delta)$, respectively. Then
$\Omega(k,\delta)\leq
G(k,\delta)$ and by Lemma \ref{simpleepi}(a) it follows that $H(k,\delta)$ is a $\pi$-Hall subgroup of~$G(k,\delta)$.

We also set $$U=\sum_{(k,\delta)\ne (2,\varepsilon)}V(k,\delta),$$ and $W=V(2,\varepsilon).$ So $V=U\perp W$. Now we
show that the dimension
of  $U$ is bounded. Clearly it is enough to prove that $d(k,\delta)$ is bounded for~${(k,\delta)\ne(2,\varepsilon)}$.

Since the decomposition \eqref{decomorth} cannot be refined, Lemma \ref{OrthogonalHallDim6} implies that $d(4,-)=0$. By
Lemma
\ref{OrthogonalHallDim6}, we
obtain that a $\pi$-Hall subgroup $H$ of $\Omega_3(q)$ is irreducible if and only if $H\simeq \Alt_5$. By using Lemmas
\ref{simpleepi}(e),
\ref{StructureInducedAutomorphisms}, and \ref{OrthogonalHallDim6}(e), we obtain that  $d(3,\circ)=0$. Thus we obtain

\begin{Lemma}\label{d(4,-)=d(3,0)=0}
In the above notation we have $d(4,-)=d(3,\circ)=0$.
\end{Lemma}

Now we need to consider the possible values for   $d(1,\circ)$,
$d(2,-\varepsilon)$, and $d(4,+)$.
First assume that $d(1,\circ)=t$. Note that \cite[Tables 3.5D, 3.5E, and~3.5F]{KlLi} imply that if $t$ is odd or
$D(1,\circ)=\square$ then $H(1,\circ)$ is included in  $\Oo_1(q)\wr \Sym_{t}$ (the stabilizer of the decomposition of 
$V(1,\circ)$ into an
orthogonal direct sum of isometric one-dimensional subspaces), while  if $t$ is even and
$D(1,\circ)=\boxtimes$, then $H(1,\circ)$ is included in a subgroup of type $\Oo_1(q)\perp \Oo_{t-1}(q)$, and the
projection of
$H(1,\circ)$ on $\Oo_{t-1}(q)$ is included in $\Oo_1(q)\wr \Sym_{t-1}(q)$. In order to distinguish between these cases
for $t$ even, we
write
respectively
$d(1,\circ)=t_{\square}$ and $d(1,\circ)=t_{\boxtimes}$.

\begin{Lemma} \label{val_d} Suppose $(\star)$ holds and the decomposition {\em \eqref{decomorth}} of $V$ cannot be
refined. In the
above-introduced notation the following statements hold{\em:}
\begin{itemize}
\item[{\em (a)}] $d(1,\circ)\in\{0,1,2_\square,2_\boxtimes,3,4_\boxtimes,6_\square\}${\em;}
\item[{\em (b)}] $d(2,-\varepsilon)\in\{0,1\}${\em;}
\item[{\em (c)}] $d(4,+)\in\{0,1\}${\em;}
\item[{\em (d)}] if $d(1,\circ)=6_\square$, then $\dim(V)=6$,  and $\eta(V)=\varepsilon${\em;}
\item[{\em (e)}] if $d(1,\circ)=4_\boxtimes$, then $\eta(1,\circ)=-$ and $d(2,-\varepsilon)=d(4,+)=0${\em;}
\item[{\em (f)}] if $d(1,\circ)=3$, then  $d(2,-\varepsilon)= d(4,+)=0${\em;}
\item[{\em (g)}] if $d(1,\circ)=2_\boxtimes$, then $\eta(1,\circ)=-\varepsilon$, and $d(2,-\varepsilon)=0${\em;}
\item[{\em (h)}] if $d(1,\circ)=1$, then $d(2,-\varepsilon)=0${\em;}
\item[{\em (i)}] if $d(1,\circ)= 2_\square$, then either $\dim(V)=2$ and $\eta(V)=\varepsilon$, or $\dim(V)=4$ and
$\eta(V)=-$, or
$\dim(V)=6$ and
$\eta(V)=\varepsilon$.
\end{itemize}
\end{Lemma}

\noindent{\slshape Proof.}\ \  Let $M$ be the stabilizer in $G(k,\delta)$ of the decomposition of $V(k,\delta)$ into the
$H$-in\-va\-ri\-ant sum of
$k$-dimensional subspaces. Set  $\Omega=\Omega(k,\delta)$, $M_\Omega=\Omega\cap M$, $L=\Omega/Z(\Omega)$, and
$M_L=M_\Omega/Z(\Omega)$.
Since $M$
includes a  $\pi$-Hall subgroup $H(k,\delta)$ of $G(k,\delta)$, Lemma \ref{simpleepi}(a) implies that  $M_L$  includes a
$\pi$-Hall subgroup
of~$L$. In particular, $|L:M_L|$ is not divisible by 2 and 3.

(a) Assume $(k,\delta)=(1,\circ)$. Denote $d(1,\circ)$ by~$t$.

Consider the case $t$ odd first.  We set $m=(t-1)/2$ for brevity.
Then
$$|L|_\ti=\prod\limits_{i=1}^m(q^{2i}-1)_\ti=(q^2-1)^m_\ti(m!)_\ti\ge 3^m(m!)_\ti,$$
$$|M_L|_\ti=(t!)_\ti=3^{[t/3]+[t/3^2]+\dots}<3^{t/3+t/3^2+\dots}=3^{t/2}.$$
Since  $t$ is odd, we obtain
$|M_L|_\ti\le 3^{(t-1)/2}=3^m$. So $|L:M_L|_\ti\ge (m!)_\ti$, whence $m<3$ and $t<7$. If $t=5$, then $|L|_\ti\ge 3^2$,
while
$|M_L|_\ti=3$, a contradiction. Thus, if $t$ is odd, then~${t\in\{1,3\}}$.

Now assume that $t$ is even. If $D(1,\circ)=\boxtimes$, then case $t$ odd and Lemma~\ref{OrthogonalHallDim6} imply
that $t\in\{2_\boxtimes,4_\boxtimes\}$. Thus it remains to consider the case~${D(1,\circ)=\square}$. Notice that
$\eta(1,\circ)=\varepsilon^{t/2}$ in this case. Set $m=t/2-1$. We have
\begin{equation*}|L|_\ti=(q^{m+1}-\varepsilon^{m+1})_\ti\prod\limits_{i=1}^m(q^{2i}-1)_\ti=(q^{m+1}-\varepsilon^{m+1}
)_\ti(q^2-1)^m_\ti(m!)_\ti\ge
3^m(m!)_\ti,\end{equation*}
$$|M_L|_\ti=(t!)_\ti=3^{[t/3]+[t/3^2]+\dots}<3^{t/3+t/3^2+\dots}=3^{t/2}=3^{m+1}.$$ Since the last inequality is strict,
it follows
that $|M_L|_\ti\le 3^m$. Therefore, $|L:M_L|_\ti\ge (m!)_\ti$, whence $m<3$ and $t<8$. For $t=4_\square$, we get
$|L|_\ti=(q^2-1)_\ti\ge 3^2$, while $|M_L|_\ti=3$ that is impossible. Hence, if $t$ is even then
$t\in\{2_\square,2_\boxtimes,4_\boxtimes,6_\square\}$.

(b)  Assume  $(k,\delta)=(2,-\varepsilon)$ and set $d(2,-\varepsilon)=t$. By Theorem \ref{Maslova0}(j), we obtain that
$|L:M_L|$ is even
if~${t>1}$.

(c) Assume $(k,\delta)=(4,+)$ and set $d(4,+)=t$. By using \cite[Proposition~4.2.11]{KlLi}, we obtain that $$\vert
L\vert_\ti=(q^2-1)_\ti^{2t}((2t)!)_\ti,\ \ \ \vert M_L\vert_\ti=(q^2-1)_\ti^{2t}(t!)_\ti,$$ hence
$|L:M_L|$ is divisible by $3$ if~${t>1}$.

(d) Assume $d(1,\circ)=6_\square$ and $n=\dim (V)>6$. Since $d(3,\circ)=0$, it follows that $\dim(V)$ is even, i.e.,
$\Omega(V)=\Omega_{2m}^\eta(q)$ for some integer $m$. Since $H(1,\circ)$ is a $\pi$-Hall subgroup of $G(1,\circ)$, then
Lemma~\ref{OrthogonalHallDim6}(v) implies that $q\equiv -\varepsilon\pmod3$. Moreover, $H$ is included in a subgroup $F$
of type
$\Oo_6^\varepsilon(q)\perp \Oo_{2m-6}^{\eta\varepsilon}(q)$. Since $Z(G)$ is a $2$-group, the inclusions $Z(G)\leq H\leq
F$ hold. So the
index of  $F/Z(G)$  in $G/Z(G)$ is a $\pi'$-num\-ber.  By \cite[Proposition~4.1.6]{KlLi}, we obtain that either
$D=\boxtimes$
(i.e., $\eta=-\varepsilon^m$) and $$(F\cap\Omega(V))/Z(G)=\left(\Omega_6^\varepsilon(q)\times
\Omega_{2m-6}^{-\varepsilon^{m-3}}(q)\right)\arbitraryext [4],$$ or $D=\square$ (i.e., $\eta=\varepsilon^m$) and
$$(F\cap\Omega(V))/Z(G)=
2\arbitraryext\left(\mathrm{P}\Omega_6^\varepsilon(q)\times
\mathrm{P}\Omega_{2m-6}^{\varepsilon^{m-3}}(q)\right)\arbitraryext [4].$$

In the first case, we have
$$\vert
G:F\vert_{p'}=\frac{(q^m+\varepsilon^m)(q^{2(m-1)}-1)(q^{2(m-2)}-1)(q^{2(m-3)}-1)}{2(q^2-1)(q^4-1)(q^3-\varepsilon)(q^{
m-3}
+\varepsilon^ { m-3 }) }.$$
Consider the $3$-part of this index:
$$\vert
G:F\vert_\ti=\frac{(q^m+\varepsilon^m)_\ti(q^2-1)^3_\ti(m-1)_\ti(m-2)_\ti(m-3)_\ti}{(q^2-1)_\ti^2(q^{m-3}+\varepsilon^{
m-3})_\ti}\ge
(q^m+\varepsilon^m)_\ti.$$ Since $q\equiv-\varepsilon\pmod3$, it follows that $(q^m+\varepsilon^m)_3>1$ if and only if
$m$ is odd. Hence if
$m$ is odd, then   $F$ does not include $\pi$-Hall subgroups of $G$ (recall that  $3\in\pi$). If $m$ is even, then
$$\vert G:F\vert_\tw=\frac{(q^2-1)_\tw^3\cdot
2\cdot(m-2)_\tw}{(q^2-1)_\tw^2\cdot(q-\varepsilon)_\tw\cdot 2\cdot
2}=\frac{(m-2)_\tw}{2},$$ i.e., we must have $m\equiv 2\pmod4$ and $\eta=-\varepsilon$ in this case. Since
$\eta=-\varepsilon^m$, it
follows that $d(2,-\varepsilon)\not=0$, so $d(2,-\varepsilon)=1$ by statement (b) of the lemma. Now Lemma
\ref{simpleepi}(a) implies that
$\left(G(1,0)\times G(2,-\varepsilon)\right)\cap\Omega_8^-(q)$ contains a $\pi$-Hall subgroup of
$\Omega_8^-(q)$. Since $D(V(1,\circ)\perp V(2,-\varepsilon))=\boxtimes$, for the subspace
$V(1,\circ)+V(2,-\varepsilon)$, we can apply the
same arguments as
for $V$. In particular, we have
$4=\frac{1}{2}\dim(V(1,\circ)+V(2,-\varepsilon))\equiv 2\pmod4$, a contradiction.

If $D=\square$, then
 \begin{equation*}
\vert
G:F\vert_\ti=\frac{(q^m-\varepsilon^m)_\ti(q^2-1)^3_\ti(m-1)_\ti(m-2)_\ti(m-3)_\ti}{(q^2-1)_\ti^2(q^{m-3}-\varepsilon^{
m-3})_\ti}\ge
(q^2-1)_\ti\ge3,
\end{equation*}
so $F$ does not include $\pi$-Hall subgroups of $G$ in this case.

(e), (f) Assume $d(1,\circ)=4_\boxtimes$, then $\eta(1,\circ)=-$. Suppose $d(2,-\varepsilon)=1$ (resp. $d(4,+)=1$). Then
$\Oo(V(1,\circ)\perp V(2,-\varepsilon))\simeq \Oo_6^\varepsilon(q)$ (resp. $\Oo(V(1,\circ)\perp V(4,+))\simeq
\Oo_8^-(q)$). Set
$\Omega=\Omega(V(1,\circ)\perp
V(2,-\varepsilon))\simeq\Omega_6^\varepsilon(q)$ (resp. $\Omega=\Omega(V(1,\circ)\perp V(4,+))\simeq \Omega_8^-(q)$) and
define groups
$M_\Omega$, $L$, and $M_L$ as above. Then $M_L$ includes a $\pi$-Hall subgroups of $L$, therefore $\vert L:M_L\vert$ is
not divisible by
$2$ and $3$. If $L\simeq\P\Omega_6^\varepsilon(q)$, then Corollary \ref{OddIndexO6} implies that $\vert L:M_L\vert$ is
even, while if
$L\simeq \P\Omega_8^-(q)$, then \cite[Proposition~4.1.6]{KlLi} implies that $\vert L:M_L\vert_\ti=3$. Statement (f) can
be obtained arguing
in the same way.

(g) The identity $\eta(1,\circ)=-\varepsilon$ is evident. The identity $d(2,-\varepsilon)=0$ follows from the fact that
the decomposition
\eqref{decomorth} cannot be
refined.

(h) can be obtained arguing as in (e) and (f) by using the fact that a subgroup of type $\Oo_1(q)\perp
\Oo_2^{-\varepsilon}(q)$ in ${\rm P}\Omega_3(q)$ has even index.

(i) If $d(1,\circ)=2_\square$, then  $\eta(1,\circ)=\varepsilon$ and, since the decomposition \eqref{decomorth} cannot
be refined, $d(2,\varepsilon)=0$.
Hence, the previous
statements imply that one of the following cases occurs:

(1) $d(2,-\varepsilon)=d(4,+)=0$, $\dim (V)=2$, $\eta(V)=\varepsilon$;

 (2) $d(2,-\varepsilon)=1$, $d(4,+)=0$, $\dim (V)=4$, $\eta(V)=-$;

(3) $d(2,-\varepsilon)=0$, $d(4,+)=1$, $\dim (V)=6$, $\eta(V)=\varepsilon$;

(4)  $d(2,-\varepsilon)=d(4,+)=1$, $\dim (V)=8$, $\eta(V)=-$.

As in the proof of statement (e) one can  show that the last case is impossible.
\qed

\begin{Lemma} \label{qequivpm3} Suppose $(\star)$ holds, $\dim(V)\ge 7$, and the decomposition {\em \eqref{decomorth}}
of $V$
cannot be refined. Assume also that  either $d(4,+)=1$ or $d(1,\circ)\ge 3$. Then $q\equiv \pm3 \pmod 8$.
\end{Lemma}

\noindent{\slshape Proof.}\ \
Assume that $(k,\delta)\in\{(1,\circ),\,(4,+)\}.$ Since $\dim(V)\ge 7$, Lemmas \ref{StructureInducedAutomorphisms} and
\ref{val_d} imply
that
$G(k,\delta)=I(V(k,\delta))\geq \SO(V(k,\delta))$. The statement follows from Lemma \ref{OrthogonalHallDim6}(c), (d),
(g), (h), (i), (j), (k), (l), (v) for a
$\pi$-Hall subgroup $H(k,\delta)\cap\SO(V(k,\delta))$ of $\SO(V(k,\delta))$.
\qed

Lemmas \ref{val_d} and \ref{qequivpm3} give us the following:

\begin{Lemma} \label{dimU} Suppose $(\star)$ holds, $\dim(V)\ge 7$, and the decomposition {\em \eqref{decomorth}} of $V$
cannot be refined.  Then $d(3,\circ)=d(4,-)=0$, and one of the following statements holds{\em:}
\begin{itemize}
\item[{\em(a)}] $\dim(U)=0$ and $d({k,\delta})=0$ for $(k,\delta)\ne (2,\varepsilon)${\em;}
\item[{\em(b)}] $\dim(U)=1$ and  $d(1,\circ)=1$ and $d(2,-\varepsilon)=d(4,+)=0${\em;}
\item[{\em(c)}] $\dim( U)=2$, $\eta(U)=-\varepsilon$,  $d(1,\circ)=2_\boxtimes$, and
$d(2,-\varepsilon)=d(4,+)=0$, or
$d(2,-\varepsilon)=1$ and $d(1,\circ)=d(4,+)=0${\em;}
\item[{\em(d)}] $\dim(U)=3$, $q\equiv \pm3 \pmod 8$,  $d(1,\circ)=3$, and $d(2,-\varepsilon)=d(4,+)=0${\em;}
\item[{\em(e)}] $\dim(U)=4,$ $\eta(U)=+$, $q\equiv \pm3 \pmod 8$,  $d(4,+)=1$, and
$d(1,\circ)=d(2,-\varepsilon)=0${\em;}
\item[{\em(f)}] $\dim(U)=4,$ $\eta(U)=-$, $q\equiv \pm3 \pmod 8$,  $d(1,\circ)=4_\boxtimes$, and
$d(2,-\varepsilon)=d(4,+)=0${\em;}
\item[{\em(g)}] $\dim(U)=5$, $q\equiv \pm3 \pmod 8$,  $d(4,+)=d(1,\circ)=1$, and $d(2,-\varepsilon)=0${\em;}
\item[{\em(h)}] $\dim(U)=6$, $\eta(U)=-\varepsilon$, $q\equiv \pm3 \pmod 8$,  and either $d(1,\circ)=2_\boxtimes$,
$d(4,+)=1$ and
$d(2,-\varepsilon)=0$, or $d(2,-\varepsilon)=d(4,+)=1$ and $d(1,\circ)=0$.
\end{itemize}
In particular, $\dim(U)\le 6$.
\end{Lemma}

\begin{Lemma}  \label{qequivvarepsilonpmod12} Suppose $(\star)$ holds and  $\dim( V)\ge 7$. Then
$s\in\pi(q-\varepsilon)$ whenever  $s\in\pi\cap \pi(q^2-1)$. In particular, $q\equiv \varepsilon \pmod
{12}$.
\end{Lemma}

\noindent{\slshape Proof.}\ \
Assume to the contrary that  $q\equiv-\varepsilon\pmod s$, in particular $s$ is odd.  Lemma \ref{dimU} implies that
$\dim(W)\ge 1$, i.e.,
$t=d(2,\varepsilon)\ge 1$.

Consider the case  $t\ge 2$ first. Let $\Omega=\Omega(2,\varepsilon)$. Then a $\pi$-Hall subgroup $H(2,\varepsilon)\cap
\Omega$ of
$\Omega$ is included in a subgroup $M_\Omega$ of type $\Oo_2^\varepsilon(q)\wr S_t$.  We also have that
$$|\Omega|_s\ge
\prod\limits_{i=1}^{t-1}(q^{2i}-1)_s=(q^2-1)_s^{t-1}\big((t-1)!\big)_s\ge s^{t-1}\big((t-1)!\big)_s,$$
  $$|M_\Omega|_s=(t!)_s<s^{t-1}\big((t-1)!\big)_s\le|\Omega|_s,$$ since $t\ge 2$. Hence $|\Omega:M_\Omega|$ is divisible
by $s$, a
contradiction.

Consider the case  $t=1$. Then $\dim(W)=2$. By Lemma \ref{dimU} it follows that one of the following statements holds.

(1) $\dim(U)=5$, $\dim( W)=2$, and $\dim( V)=7$.

(2) $\dim( U)=6$, $\eta(U)=-\varepsilon$, $\dim( W)=2$, $\eta(W)=\varepsilon$, $\dim( V)=8$,  and ${\eta(V)=-}$. In this
case we also have that
$d(4,+)=1$ and either  $d(1,\circ)=2_\boxtimes$, or $d(2,-\varepsilon)=1$.

Statement (1) is impossible, since by \cite[Proposition~4.1.6]{KlLi} the index of a subgroup of type $\Oo_5(q)\perp
\Oo_2^\varepsilon(q)$ in
${\rm P}\Omega_7(q)$ is divisible by 3. Statement (2) can be eliminated arguing as in the proof of statements (e) and
(i) of
Lemma~\ref{val_d}.
\qed

By Lemma  \ref{qequivvarepsilonpmod12} it follows that if   $\dim(V)\ge 7$, hypothesis $(\star)$ holds, and the
decomposition \eqref{decomorth} of $V$
cannot be refined, then a $\pi$-Hall subgroup of $\Oo_2^{-\varepsilon}(q)$ has order $4$ and is equal to the
stabilizer of a decomposition of the natural $2$-dimensional module into a sum of nondegenerate $1$-dimensional
subspaces.
Hence we also have
\begin{equation}\label{d(2,-epsilon)=0}
{d(2,-\varepsilon)=0}.\end{equation}

\begin{Lemma}
\label{HallSubgroupsOfOrthogonalGroupsOfEvenDimension}\label{HallSubgroupsOfOrthogonalGroupsOfOddDimension} Assume that
$G=\Omega_n^\eta(q)$, $\eta\in\{+,-,\circ\}$, $q$ is a power of a prime $p$, $n\ge 7$, $\varepsilon=\varepsilon(q)$. Let
$\pi$ be a set of
primes such that  $2,3\in\pi$, $p\not\in\pi$. Then the following statements hold{\em:}

\begin{itemize}
\item[{\em(A)}] If  $G$ possesses a $\pi$-Hall subgroup $H$, then
one of the following statements holds{\em:}
\begin{itemize}
\item [{\em (a)}] $n=2m+1$, $\pi\cap\pi(G)\subseteq\pi(q-\varepsilon)$, $q\equiv \varepsilon\pmod {12}$, and $\Sym_m\in
E_\pi$. The subgroup $H$ is a
$\pi$-Hall subgroup in $M=\big(\Oo_2^\varepsilon(q)\wr \Sym_m\times \Oo_1(q)\big)\cap G.$
All $\pi$-Hall subgroup of this type are conjugate.

\item [{\em (b)}] $n=2m$, $\eta=\varepsilon^m$, $\pi\cap\pi(G)\subseteq\pi(q-\varepsilon)$, $q\equiv \varepsilon\pmod
{12}$, and $\Sym_m\in E_\pi$. The subgroup $H$ is a $\pi$-Hall
subgroup in $M=\big(\Oo_2^\varepsilon(q)\wr \Sym_m\big)\cap G.$ All $\pi$-Hall subgroup of this type are conjugate.

\item [{\em (c)}] $n=2m$, $\eta=-\varepsilon^m$, $\pi\cap\pi(G)\subseteq\pi(q-\varepsilon)$, $q\equiv \varepsilon\pmod
{12}$, and $\Sym_{m-1}\in E_\pi$. The subgroup $H$ is a
$\pi$-Hall subgroup of $M=\big(\Oo_2^\varepsilon(q)\wr \Sym_{m-1}\times \Oo^{-\varepsilon}_2(q)\big)\cap G.$ All
$\pi$-Hall subgroup of this type are conjugate.

\item [{\em (d)}] $n=11$,  $\pi\cap\pi(G)=\{2,3\}$, $q\equiv \varepsilon\pmod {12}$, and  $(q^2-1)_\pi=24$. The
subgroup $H$ is a $\pi$-Hall subgroup
of $M=\big(\Oo_2^\varepsilon(q)\wr \Sym_{4}\times \Oo_1(q)\wr \Sym_3\big)\cap G.$ All
$\pi$-Hall subgroup of this type are conjugate.

\item [{\em (e)}] $n=12$, $\eta=-$, $\pi\cap\pi(G)=\{2,3\}$, $q\equiv \varepsilon\pmod {12}$, and  $(q^2-1)_\pi=24$.
The subgroup $H$ is a $\pi$-Hall
subgroup of $M=\big(\Oo_2^\varepsilon(q)\wr \Sym_{4}\times  \Oo_1(q)\wr \Sym_3\times \Oo_1(q)\big)\cap G$.
There exist precisely two classes of conjugate subgroups of this type in $G$,
and the
automorphism of order $2$ induced by the group of similarities of the natural module interchanges these classes.

\item[{\em (f)}] $n=7$, $\pi\cap\pi(G)=\{2,3,5,7\}$, and
$\vert G\vert_\pi=2^9\cdot 3^4\cdot5\cdot7$. The subgroup  $H$ is isomorphic to $\Omega_7(2)$. There exist
precisely two classes of conjugate subgroups of
this type in $G$, and  ${\rm SO}_7(q)$ interchanges these classes.

\item[{\em (g)}] $n=8$, $\eta=+$, $\pi\cap\pi(G)=\{2,3,5,7\}$, and $\vert G\vert_\pi=2^{13}\cdot3^5\cdot5^2\cdot7$. The
subgroup $H$ is isomorphic to $\Omega_8^+(2)$. There exist precisely four classes of conjugate
subgroups of
this type in $G$. The subgroup of $\Out(G)$ generated by diagonal and graph automorphisms is isomorphic to $\Sym_4$ and
acts on the set of
these classes as $\Sym_4$ in its natural permutation representation, and every diagonal automorphism acts without fixed
points.

\item[{\em (h)}]\ $n=9$, $\pi\cap\pi(G)=\{2,3,5,7\}$, and
$\vert G\vert_\pi=2^{14}\cdot3^5\cdot5^2\cdot7$. The subgroup $H\simeq 2\arbitraryext\Omega_8^+(2)\arbitraryext
2$. There exist precisely two classes
of conjugate subgroups of
this type in $G$, and  ${\rm SO}_9(q)$ interchanges these classes.
\end{itemize}
\item[{\em (B)}] Conversely, if arithmetic conditions in one of the statements {\rm (a)--(h)} holds, then $G$ possesses
a $\pi$-Hall subgroup with the given
structure.
\item[{\em (C)}] If $G\in E_\pi$, then  $k_\pi(G)\in\{1,2,3,4\}$.
\item[{\em (D)}] If $H\in \Hall_\pi(G)\not=\emptyset$, then $G$ possesses a $\pi$-subgroup not conjugate to a subgroup
of $H$, in
particular $G\not\in D_\pi$.
\item[{\em (E)}] All $\pi$-Hall subgroups of $\P G$ have the form~${\P H}$.
\end{itemize}
\end{Lemma}

\noindent{\slshape Proof.}\ \  (A) Suppose $H$ is a $\pi$-Hall subgroup of $G$. Assume first that  $(\star)$ does not
hold. Then $n=7$, $8$, or $9$ and
$H$ is isomorphic to $\Omega_7(2)$, $2\arbitraryext \Omega_8^+(2)$, or $2\arbitraryext\Omega_8^+(2)\splitext 2$,
respectively. Thus we
obtain statements (f), (g), (h) of the lemma and we only need to prove the statements about the classes of conjugate
$\pi$-Hall subgroups
in this case.

(f)  By \cite[Lemma~1.7.1]{Kl} the number of
classes of conjugate subgroups isomorphic to $\Omega_7(2)$ in
${\Oo}_7(q)$ is not greater than the number of nonequivalent
irreducible representation of $\Omega_7(2)$  of degree $7$ over a field of order $q$. By
using ordinary character tables of $\Omega_7(2)$ given in \cite{ATLAS} or
\cite{GAP} we obtain that it has precisely one irreducible complex
representation of degree
$7$, while $2.\Omega_7(2)$ has no
complex representations of degree $7$. Since $p$ does not divide the order of $\Omega_7(2)$, by
\cite[Theorem~15.3 and Corollary~9.7]{Isaacs}, this statement also holds for the
representations over $\mathbb{F}_q$. Thus there exist at most two classes of conjugate $\pi$-Hall subgroups isomorphic
to
$\Omega_7(2)$. By Lemma \ref{HallAmostSimpleInduce}, it follows that $\SO_7(q)$ does not possesses a $\pi$-Hall subgroup
$H_1$ such that
$H_1\cap\Omega_7(q)\simeq \Omega_7(2)$. Thus by Lemma \ref{simpleepi}(e) there exists at least two classes of conjugate
$\pi$-Hall subgroups
isomorphic to
$\Omega_7(2)$, and $\SO_7(q)$ interchanges these classes.

(g) By
\cite[Proposition~2.3.8]{Kl}, we obtain that there exist $4$ classes of
conjugate subgroups isomorphic to $\Omega_8^+(2)$ in $\mathrm{P}\Omega_8^+(q)$.  Moreover $\Aut(S)$
interchanges these four classes in the natural way. Therefore there exists a homomorphism
$\Aut(S)\rightarrow \mathrm{Sym}_4$ and, by \cite[Proposition~2.3.8]{Kl},
this homomorphism is surjective.

(h) In this case a $\pi$-Hall subgroup $H\simeq2\arbitraryext\Omega_8^+(2)\arbitraryext 2$ is included in a subgroup of
type
$\Oo_1(q)\perp\Oo_8^+(q)$ which, by Lemma \ref{StructureInducedAutomorphisms}(a), is isomorphic to
$\Omega_8^+(q)\splitext\langle\tau\rangle$, where $\tau$ is a graph
automorphism of order $2$, and by
\cite[Table~3.5D]{KlLi} all such subgroups are conjugate. By
\cite[Proposition~2.3.8]{Kl}, there exists $4$ classes of conjugate
subgroups isomorphic to $2.\Omega_8^+(2)$ in $\Omega_8^+(q)$, and the graph
automorphism $\tau$ of order $2$ stabilizes $2$ of these classes and
interchanges the remaining $2$ classes. By Lemma \ref{simpleepi}(e),
$\Omega_8^+(q)\splitext\langle\tau\rangle$ has $2$ classes  of conjugate
$\pi$-Hall subgroups isomorphic to
$2\arbitraryext\Omega_8^+(2)\arbitraryext 2$. By Lemma \ref{HallAmostSimpleInduce}, it follows that $\SO_9(q)$ does not
possesses a
$\pi$-Hall
subgroup $H_1$ such that
$H_1\cap\Omega_9(q)\simeq 2\arbitraryext\Omega_8^+(2)\arbitraryext 2$. Thus $\SO_9(q)$ interchanges these classes.

Now assume that $(\star)$ holds. Consider a decomposition of the natural module $V$ of  $G$ into $H$-in\-va\-ri\-ant
direct orthogonal sum
of subspaces of dimension at most $4$ such that the decomposition cannot be refined (the existence of this decomposition
follows from Lemma
\ref{orthored}). We preserve the above notations, in particular, the symbols $V(k,\delta),$ $U$, and $W$ are as defined
above. Then one of
statements (a)--(h) of Lemma  \ref{dimU} holds. Moreover, as we noted before Lemma
\ref{HallSubgroupsOfOrthogonalGroupsOfEvenDimension},
the condition that the decomposition cannot be refined implies that $d(2,-\varepsilon)=0$. Consider statements (a)--(h)
of Lemma
\ref{dimU} separately.

Suppose statement (a) of Lemma  \ref{dimU} holds. Then $V=W=V(2,\varepsilon)$. So $n=\dim (V)=2m$, where
$m=d(2,\varepsilon)$,
$\eta =\varepsilon^m$, and
$H\leq (\Oo_2^\varepsilon(q)\wr \Sym_m)\cap G$. Lemmas \ref{SymmetricDivizors} and \ref{qequivvarepsilonpmod12} imply
statement (b) of the
lemma. By
\cite[Tables~3.5.E and~3.5.F]{KlLi}, we obtain that all subgroups of type $\Oo_2^\varepsilon(q)\wr S_m$ are conjugate in
$G$. By Lemma
\ref{HallSymmetric}, we obtain $\Sym_m\in C_\pi$.  Since
$\Oo_2^\varepsilon(q)$ is solvable, Lemma \ref{simpleepi}(f) implies that
$(\Oo_2^\varepsilon(q)\wr \Sym_m)\cap G \in C_\pi$, hence all $\pi$-Hall subgroups of this type are conjugate.

Suppose statement (b) of Lemma  \ref{dimU} holds. Then  $n=\dim (V)=2m+1$, where $m=d(2,\varepsilon)$ and $H\leq
((\Oo_2^\varepsilon(q)\wr \Sym_m)\perp \Oo_1(q))\cap G$. Using \cite[Table~3.5.D]{KlLi} instead of \cite[Tables~3.5.E
and~3.5.F]{KlLi} we
obtain
statement (a) of the lemma as in the previous case.

Suppose statement (c) of Lemma  \ref{dimU} holds. In this case $n=\dim (V)=2m$,  $d(2,\varepsilon)=m-1$,
$d(1,\circ)=2_{\boxtimes}$ (in view of \eqref{d(2,-epsilon)=0}, the case $d(2,-\varepsilon)=1$ is impossible), $\eta
=-\varepsilon\cdot\varepsilon^{m-1}=-\varepsilon^m$ and $H\leq ((\Oo_2^\varepsilon(q)\wr \Sym_{m-1})\perp
\Oo^{-\varepsilon}_2(q))\cap G$. As in statement (a) of Lemma \ref{dimU} we obtain that statement~(c) of the lemma holds
in this case.

Suppose statement (d) of Lemma  \ref{dimU} holds, i.e., $\dim (U)=3$, $U=V(1,\circ)$. Set
$d(2,\varepsilon)=m-1$. Then $\eta(W)=\varepsilon^{m-1}$, $n=\dim(V)=2m+1$. Lemma \ref{StructureInducedAutomorphisms}
implies that
$G(1,\circ)\simeq\Oo_3(q)$. Now $H(1,\circ)$ is a $\pi$-Hall subgroup of $G(1,\circ)$, stabilizing a decomposition of
$U$ into a sum of
$1$-di\-men\-si\-o\-nal subspaces. Lemmas \ref{OrthogonalHallDim6}(b)--(e) and \ref{simpleepi}(e) imply that $\vert
H(1,\circ)\vert=24$ and
$H(1,\circ)\simeq
\Sym_4$, hence  $(q^2-1)_\pi=|\Oo_3(q)|_\pi=24$. By \cite[Proposition~4.1.6]{KlLi}, we obtain that the $p'$-part of
index of a subgroup $M$
of type $\Oo_{2(m-1)}^{\varepsilon^{m-1}}(q)\perp \Oo_3(q)$ in $G$ is equal to
$$\vert L:M\vert_{p'}=\frac{1}{2}\cdot
\frac{q^{2m}-1}{q^2-1}(q^{m-1}+\varepsilon^{m-1}).$$  Lemma \ref{qequivvarepsilonpmod12} implies
that $q\equiv\varepsilon\pmod3$, so we obtain
$(q^{m-1}+\varepsilon^{m-1})_\ti=1$. Since $M$ contains $H$, it follows that $1=|L:M|_\tw=m_\tw$ and
$1=|L:M|_{\ti}=m_\ti$, whence    $m$ is not divisible by  2 and by 3, i.e., $m\equiv \pm 1\pmod 6$. By Lemma
\ref{SymmetricDivizors}, we obtain that every prime from  $\pi\cap \pi(\Sym_{m-1})$ is contained in $\pi(q^2-1)$, i.e.,
$\Sym_{m-1}\in
E_{\{2,3\}}$. Now Lemma \ref{HallSymmetric} implies the inequalities   $m\le 9$ and $m\not=7$. All restrictions on $m$
give us the only
possible value $m=5$, whence statement (d) of the lemma follows. The statement about the classes of conjugate $\pi$-Hall
subgroups we
obtain as in statement (a) of Lemma \ref{dimU}, using \cite[Table~3.5.D]{KlLi} instead of \cite[Tables~3.5.E
and~3.5.F]{KlLi}.

Suppose statement (e) of Lemma  \ref{dimU} holds, i.e., $d(4,+)=1$, $U=V(4,+)$, $\dim (U)=4$,
$\eta(U)=+$. Set $d(2,\varepsilon)=m-2$. Then  $\dim(V)=2m$, $\eta (V)=\varepsilon^{m-2}=\varepsilon^m$. Lemma
\ref{StructureInducedAutomorphisms} implies that $G(4,+)\simeq \Oo_4^+(q)$. Now Lemma \ref{OrthogonalHallDim6}(f)--(l)
and the fact that
the decomposition \eqref{decomorth} cannot be refined imply that $H(4,+)\cap \Omega(4,+)$ is isomorphic to
$2\arbitraryext\Alt_4\circ2\arbitraryext\Alt_4$, whence $(q^2-1)_\pi=24$. By \cite[Proposition~4.1.6]{KlLi}, the
$p'$-part of
the index of a
subgroup $M$ of type $\Oo_{2(m-2)}^{\varepsilon^{m-2}}(q)\perp \Oo_4^+(q)$ in $G$ is equal to
$$\vert G:M\vert_{p'}=\frac{q^{m-2}+\varepsilon^{m-2}}{2}\cdot
\frac{(q^{2(m-1)}-1)(q^m-\varepsilon^m)}{(q^2-1)^2}.$$ Since, by Lemma \ref{qequivvarepsilonpmod12}, we have
$q\equiv\varepsilon\pmod3$, it
follows that $$|G:M|_\tw=\frac{m_\tw(m-1)_\tw}{2} ,\,\,|G:M|_{\ti}=m_\ti(m-1)_\ti.$$ Moreover $H$ is included in $M$,
whence  $(m-2)(m-3)$
is
divisible by 4, and $m-2$ is divisible by 3. Thus $m\equiv
2,11\pmod{12}$. Lemma \ref{SymmetricDivizors} implies that each prime in $\pi\cap \pi(\Sym_{m-2})$ is contained in
$\pi(q^2-1)$, i.e.,
$\Sym_{m-2}\in E_{\{2,3\}}$. Lemma \ref{HallSymmetric} implies that $m\le 10$ and $m\not=8$.
All restrictions on $m$ and the condition  $n\ge7$ imply that there does not exist a $\pi$-Hall subgroup satisfying
statement (e) of
Lemma~\ref{dimU}.

Suppose statement (f) of Lemma  \ref{dimU} holds. By definition, the identity $d(1,\circ)=4_\boxtimes$ implies that this
case follows from
the already considered statement (d) of Lemma  \ref{dimU}, and statement (e) holds.

Suppose statement (g) or (h) of Lemma  \ref{dimU} holds. As in the already considered statement (e) of Lemma 
\ref{dimU}, we obtain that
there
does not exist a $\pi$-Hall subgroup satisfying statement (g) or (h) of
Lemma~\ref{dimU}, respectively.

(B) Clearly a subgroup satisfying statement (f), (g) or (h) of the lemma is a $\pi$-Hall subgroup in the corresponding
group $G$. If one of
statements (a)--(e) of the lemma holds, then a direct calculation  using Lemma \ref{FWair} and Corollaries \ref{3part},
\ref{2part}
implies that $\vert
G:M\vert_\pi=1$, hence $\emptyset\not=\Hall_\pi(M)\subseteq \Hall_\pi(G)$.

(C) If one of statements (f), (g), (h) of the lemma holds, then $\vert G\vert_7=7$, whence $7\not\in\pi(q-\varepsilon)$
and none of
statements (a)--(e) can be fulfilled. Since either statements (a) and (d), or statements (c) and (e) of the lemma can be
fulfilled
simultaneously, and the remaining cannot,  we obtain that in this
case $k_\pi(G)=3$ and (C) follows.

(D) If one of the statements (e)--(h) of the lemma holds, then $G$ possesses more than one class of conjugate $\pi$-Hall
subgroups and we
have nothing to prove. Assume that statement (d) of the lemma holds. Then $d(2,\varepsilon)=4$, and $n=\dim( V)=11$. 
Then  $L=\big((\Oo_2^\varepsilon(q)\wr \Sym_4)\times \Oo_2^\varepsilon(q)\times O_1(q)\big)\cap G$  possesses a
$\pi$-Hall subgroup $F$, and $F$ is not isomorphic to a subgroup of $H$. Assume that one of statements (a) and (b) of
the lemma holds. Then
$G$ possesses a subgroup of type $\Oo_1(q)\wr S_n$, which contains a subgroup  $2^{n-1}\arbitraryext \Alt_n$. 
Statements (a) and (b) imply
that $\Sym_m\in E_\pi$, where $m=[n/2]$. Consider a subgroup $\Alt_m\times \Alt_m\leq \Alt_n$, let  $R$ be its
$\pi$-Hall subgroup and
$R_1$ the complete preimage in $2^{n-1}\arbitraryext \Alt_n$. It is clear that $R_1$ is not isomorphic to any subgroup
of  $H$. If statement
(c) holds, we proceed in the same way considering a subgroup of type~$\Oo_1(q)\perp (\Oo_1(q)\wr \Sym_{n-1})$.

(E) Follows from Lemma \ref{simpleepi}(a), (c), (f).
\qed

\section{Hall subgroups in exceptional groups of Lie type}

Recall that $q=p^\alpha$ is a power of a prime $p$,  $\varepsilon=\varepsilon(q)=(-1)^{(q-1)/2}$, $\pi$ is a set of
primes such that $2,3\in\pi$,  and
$p\not\in\pi$.

\begin{Lemma}\label{HallG2}
Suppose $G\simeq G_2(q)$, $2,3\in\pi$, $p\not\in\pi$, and $\varepsilon=\varepsilon(q)$. Then the following statements
hold{\em:}
\begin{itemize}
\item[{\em (A)}] If $G$ possesses  a $\pi$-Hall subgroup  $H$, then one of the following statements holds{\em:}
\begin{itemize}
\item[{\em (a)}] $\pi\cap \pi(G)=\{2,3,7\}$, $(q^2-1)_{\{2,3,7\}}=24$, and $(q^4+q^2+1)_7=7$. The subgroup  $H$ is
isomorphic to $G_2(2)$, and all $\pi$-Hall  subgroups of
this type are conjugate in~$G$;
\item[{\em (b)}] $\pi\cap \pi(G)\subseteq \pi(q-\varepsilon)$. The subgroup $H$ is a $\pi$-Hall subgroup of a solvable
group\linebreak
${({q-\varepsilon})}^2\arbitraryext W(G_2)$, and all $\pi$-Hall subgroups of this type are conjugate in~$G$.
\end{itemize}
\item[{\em (B)}] Conversely, if arithmetic conditions in one of statements {\em (a), (b)} holds, then $G\in E_\pi$.
\item[{\em (C)}] If $G\in E_\pi$, then  $k_\pi(G)=1$, i.e., $G\in C_\pi$.
\item[{\em (D)}] If $G\in E_\pi$ and $H\in \Hall_\pi(G)$, then $G$ possesses a $\pi$-subgroup that is not isomorphic to
a subgroup of~$H$,
in particular,~${G\not\in~D_\pi}$.
\end{itemize}
\end{Lemma}

\noindent{\slshape Proof.}\ \  (A) Suppose $H$ is a $\pi$-Hall subgroup of $G$. Then it is included in a maximal
subgroup $M$ of odd index of $G$.  By
\cite{liesax} one of the following statements holds (recall that $p\not\in\pi$, while $2\in\pi$, whence $p$ is odd):
\begin{itemize}
\item[(1)] $M\simeq G_2(q_0)$, where $\F_{q_{\mbox{}_0}}$ is a subfield of~$\F_q$.
\item[(2)] $M\simeq G_2(2)$, $\vert M\vert=2^6\cdot 3^3\cdot7$,
\item[(3)]  $M=N_G(2\arbitraryext(\mathrm{PSL}_2(q)\times\mathrm{PSL}_2(q)))$, $\vert M\vert=q^2(q^2-1)^2$,
\item[(4)] $M\simeq 2^3\arbitraryext\mathrm{SL}_3(2)$, $\vert M\vert=2^6\cdot 3\cdot 7$,
\item[(5)] $M\simeq \mathrm{SL}_3^\varepsilon(q)\arbitraryext 2$, $\vert M\vert=q^3(q^3-\varepsilon)(q^2-1)\cdot 2$,
\item[(6)] $M\simeq (q-\varepsilon)^2\arbitraryext W(G_2)$, $\vert M\vert=(q-\varepsilon)^2\cdot2^2\cdot3$.
\end{itemize}
Suppose statement (1) holds. Since $p\not\in\pi$ we obtain that $H$  is a proper subgroup of $M$. So $H$ is included in
a maximal subgroup
$M_1$ of odd index of $M\simeq G_2(q_0)$. In particular, $M_1$ satisfies one of statements (1)--(6). If $M_1$ satisfies
statement (1), we
again obtain that $H$ is a proper subgroup of $M_1$ and repeat the arguments. So we may assume that $M$ satisfies one of
statements
(2)--(6).\footnote{The case, when a maximal subgroup $M$ is a subgroup of the same type, as $G$, defined over a proper
subfield in proofs of
Lemmas \ref{HallF4}--\ref{Hall3D4} can be excluded by arguing in the same way, and we do not consider this case.}
Assume statement (2) holds. Then   \cite[Theorem~1.2]{RevinSIBAM} implies that   $G_2(2)$ does not possesses proper 
$\pi$-Hall subgroups
with $2,3\in\pi$. Therefore  $H=G_2(2)$. Moreover all subgroups isomorphic to $G_2(2)$ are conjugate in
view of
\cite[Theorem~A]{Kl0}, whence statement (a) of the lemma follows.
If either (3) or (4) holds, then $\vert G:M\vert_\ti>1$, whence $M$ cannot include a $\pi$-Hall subgroup of $G$ with
$3\in\pi$.  Assume
statement (5) holds. Then $\vert G:M\vert_\ti=1$ if and only if  $q\equiv\varepsilon\pmod3$. Moreover we obtain that a
$\pi$-Hall subgroup
$H_1$ of $\SL_3^\varepsilon(q)$ satisfies either statement (b), or statement (c)  of Lemma
\ref{HallSubgroupsOfLinearAndUnitaryGroups}. If $H_1$
satisfies statement (b) of Lemma \ref{HallSubgroupsOfLinearAndUnitaryGroups}, then $H_1$ normalizes a maximal torus of
order
$(q-\varepsilon)^2$, so $H\leq M$, where $M$ satisfies (6). If $H_1$ satisfies statement (c) of Lemma
\ref{HallSubgroupsOfLinearAndUnitaryGroups}, then $H_1$ is
included in a subgroup $\GL_2^\varepsilon(q)$ of $\SL_3^\varepsilon(q)$. Condition $q\equiv\varepsilon\pmod3$ implies
that $\vert
\SL_3^\varepsilon(q):\GL_2^\varepsilon(q)\vert_\ti=3$, hence $H_1$ cannot satisfy statement  (c) of Lemma
\ref{HallSubgroupsOfLinearAndUnitaryGroups}.  Finally, if statement (6) holds, then direct calculations show that the
condition $\vert
G:(q-\varepsilon)^2\arbitraryext W(G_2)\vert_\ti=1$ implies $q\equiv\varepsilon\pmod3$. The conjugacy of maximal tori of
order
$(q-\varepsilon)^2$ is obtained in \cite[Lemma~3.10]{RevVdoContemp}, whence  we obtain statement (b) of the lemma.

(B) If statement (a) of the lemma holds, then the index $\vert G:H\vert$ is a $\pi'$-number, so $G_2(2)$  is a
$\pi$-Hall subgroup of $G$. If  statement (b) of the lemma holds, then the index $\vert
G:(q-\varepsilon)^2.W(G_2)\vert$ is a $\pi'$-num\-ber, so $\emptyset\not=\Hall_\pi((q-\varepsilon)^2\arbitraryext
W(G_2))\subseteq
\Hall_\pi(G)$.

(C) Notice that  statement (a) of the lemma implies that $7$ does not divide $(q^2-1)$, while statement (b) implies
that $7$ divides $q^2-1$, whence these statements cannot hold simultaneously. Hence $k_\pi(G)=1$ if~${G\in E_\pi}$.

(D) In view of \cite{CLSS}, there exists a subgroup $L$ of $G$ isomorphic to  $2^3\arbitraryext\mathrm{SL}_3(2)$.  If
statement (a) holds,
then $L$ is a $\pi$-group and $L$ clearly is not isomorphic to a subgroup of $G_2(2)$. In the second case of the proof
of Lemma~6.9 in
\cite{RevVdoContemp}, we have constructed a $\{2,3\}$-subgroup of $L$ which cannot be embedded in a maximal torus,
whence if statement (b)
holds, then $G$ possesses a $\pi$-subgroup not isomorphic to a subgroup of a $\pi$-Hall subgroup of~$G$.
\qed

\begin{Lemma}\label{HallF4}
Suppose that  $G\simeq F_4(q)$, $2,3\in\pi$, $p\not\in\pi$, and $\varepsilon=\varepsilon(q)$. Then the following
statements hold{\em:}
\begin{itemize}
\item[{\em (A)}] If $G\in E_\pi$ and  $H\in\Hall_\pi(G)$, then  $\pi\cap \pi(G)\subseteq\pi(q-\varepsilon)$ and $H$ is
included in a
solvable group $(q-\varepsilon)^4\arbitraryext W(F_4)$.
\item[{\em (B)}] Conversely, if  $\pi\cap \pi(G)\subseteq\pi(q-\varepsilon)$, then  $G\in E_\pi$.
\item[{\em (C)}] If $G\in E_\pi$, then $k_\pi(G)=1$, i.e., $G\in C_\pi$.
\item[{\em (D)}] If $H\in\Hall_\pi(G)$, then $G$ possesses a $\pi$-subgroup  not isomorphic to a subgroup of~$H$. In
particular, $G\not\in
D_\pi$.
\end{itemize}
\end{Lemma}

\noindent{\slshape Proof.}\ \  (A) Assume that $G$ possesses a $\pi$-Hall subgroup  $H$. Then $H$ is included in a
maximal subgroup  $M$ of odd index  of
$G$. As in the proof of Lemma \ref{HallG2}, by using \cite{liesax} we may assume that one of the following statements
holds:
\begin{enumerate}
\item[(1)] $M=N_G(2\arbitraryext \mathrm{P}\Omega_9(q))$, $\vert M\vert=q^{16}(q^2-1)(q^4-1)(q^6-1)(q^8-1)$,
\item[(2)]  $M=N_G(2^2\arbitraryext\mathrm{P}\Omega_8^+(q))$, $\vert M\vert=q^{12}(q^2-1)(q^4-1)^2(q^6-1)\cdot2\cdot3$.
\end{enumerate}
If statement (1) holds, then  $\vert G:M\vert_\ti\ge 3$, hence $M$ cannot contain a $\pi$-Hall subgroup of $G$. If
statement  (2) holds,
then
by Lemma \ref{simpleepi}(a) the image of $H$ in $\mathrm{P}\Omega_8^+(q)$ is a $\pi$-Hall subgroup of
$\mathrm{P}\Omega_8^+(q)$. By Lemma \ref{HallSubgroupsOfOrthogonalGroupsOfEvenDimension}, we obtain that this image is
either  equal to
$\Omega_8^+(2)$, or included in the normalizer of a maximal torus (i.e., statement (b) of Lemma
\ref{HallSubgroupsOfOrthogonalGroupsOfEvenDimension} holds). If $H$ is included in the normalizer of a maximal torus, we
obtain statement
(A) of the lemma. If the image of  $H$ in $\mathrm{P}\Omega_8^+(q)$ equals
$\Omega_8^+(2)$, then $\vert N_G(2^2\arbitraryext\mathrm{P}\Omega_8^+(q))\vert_{\mbox{}_7}=7$, and  $\vert
G\vert_{\mbox{}_7}\ge 7^2$. Hence,
this case is impossible.

(B), (C) and (D) are proven in~\cite[Lemma~6.8]{RevVdoContemp}.
\qed

\begin{Lemma}\label{HallE6}
Suppose that $G\simeq E_6^\eta(q)$, $2,3\in\pi$, $p\not\in\pi$, and $\varepsilon=\varepsilon(q)$. Then the following
statements hold{\em:}
\begin{itemize}
\item[{\em (A)}] If $G\in E_\pi$ and $H\in\Hall_\pi(G)$, then $\pi\cap \pi(G)\subseteq \pi(q-\varepsilon)$ and one of
the following
statements holds:
\begin{itemize}
\item[{\em(a)}] $\eta=\varepsilon$, $5\in\pi$, and  $H$ is a $\pi$-Hall subgroup in $T\arbitraryext W(E_6)$, where $T$
is a maximal split
torus of order ${(q-\eta)^6/(3,q-\eta)}$;
\item[{\em(b)}] $\eta=-\varepsilon$ and $H$ is a  $\pi$-Hall subgroup in ${(q^2-1)^2(q+\eta)}^2\arbitraryext W(F_4)$.
\end{itemize}
\item[{\em (B)}] Conversely, if $\pi\cap \pi(G)\subseteq \pi(q-\varepsilon)$, and $5\in\pi$ for $\eta=\varepsilon$, then
$G\in E_\pi$.
\item[{\em (C)}] If $G\in E_\pi$, then $k_\pi(G)=1$, i.e., $G\in C_\pi$.
\item[{\em (D)}] If $H\in\Hall_\pi(G)$, then $G$ possesses a $\pi$-subgroup  not isomorphic to a subgroup of~$H$. In
particular, $G\not\in
D_\pi$.
\end{itemize}
\end{Lemma}

\noindent{\slshape Proof.}\ \
If we show that a $\pi$-Hall subgroup $H$ normalizes a maximal torus, then the claim follows from
\cite[Lemma~6.5]{RevVdoContemp}.  Suppose
that
$H$ is a $\pi$-Hall subgroup of $G$. Then it is included in a maximal subgroup $M$ of odd index of $G$. In view of
\cite{liesax} and arguing
as in the proof of Lemma \ref{HallG2}, we may assume that one of the
following statements holds:
\begin{enumerate}
\item[(1)] $\eta=+$, %$H$ is a parabolic subgroup with Levi factor of type~$D_5$,
$\vert M \vert=\frac{q^{36}(q-1)(q^2-1)(q^4-1)(q^5-1)(q^6-1)(q^8-1)}{(3,q-1)}$;
\item[(2)]  $M=N_G((4,q-\eta)\arbitraryext \mathrm{P}\Omega_{10}^\eta(q))$, $\vert
M\vert=\frac{q^{20}{(q-\eta)(q^2-1)(q^4-1)}\cdot{(q^5-\eta)(q^6-1)(q^8-1)}}{(3,q-\eta)}$;
\item[(3)] $M=N_G(2^2\arbitraryext\mathrm{P}\Omega_8^+(q))$, $\vert
M\vert=\frac{q^{12}(q-\eta)^2(q^2-1)(q^4-1)^2(q^6-1)\cdot2\cdot3}{(3,q-\eta)}$;
\item[(4)] $M=T\arbitraryext W(E_6)$,  where $T$ is a maximal split
torus of order $\frac{(q-\eta)^6}{(3,q-\eta)}$, and $\vert M\vert=\frac{(q-\eta)^6\cdot 2^7\cdot 3^4\cdot
5}{(3,q-\eta)}$.
\end{enumerate}
If either statement (1) or statement (2) holds, then $\vert G:M\vert_\ti\ge3$, whence  $M$ cannot include a $\pi$-Hall
subgroup of $G$.
If
statement (3) holds, then by Lemma \ref{simpleepi}(a) we obtain that the image of $H$ in $\mathrm{P}\Omega_8^+(q)$ is a
$\pi$-Hall subgroup
of $\P \Omega_8^+(q)$. Lemma \ref{HallSubgroupsOfOrthogonalGroupsOfEvenDimension} implies that either this image
satisfies statement (b) of
Lemma \ref{HallSubgroupsOfOrthogonalGroupsOfOddDimension}, or this image is isomorphic to $\Omega_8^+(2)$. In the first
case it follows that
$H$ normalizes a maximal torus, in the second case we obtain that $\vert H\vert_{{\mbox{}_7}}=7$, while
$\vert G\vert_{{\mbox{}_7}}\ge7^2$, a contradiction. If statement (4) holds, then $M$ normalizes a maximal torus.
\qed

\begin{Lemma}\label{HallE7}
Suppose that $G\simeq E_7(q)$, $2,3\in\pi$, $p\not\in\pi$, $\varepsilon=\varepsilon(q)$ and  $T$ is a maximal torus  of
order $(q-\varepsilon)^7/2$ of $G$.  Then the
following  statements hold{\em:}
\begin{itemize}
\item[{\em (A)}] If $H$ is a $\pi$-Hall subgroup of $G$, then $\pi\cap \pi(G)\subseteq \pi(q-\varepsilon)$, $5,7\in\pi$,
and $H$,  up to
conjugation, is included in $T\arbitraryext W(E_7)$.
\item[{\em (B)}] Conversely, if  $\pi\cap \pi(G)\subseteq \pi(q-\varepsilon)$ and $5,7\in\pi$, then $G\in E_\pi$.
\item[{\em (C)}] If $G\in E_\pi$, then $k_\pi(G)=1$, i.e., $G\in C_\pi$.
\item[{\em (D)}] If $H\in\Hall_\pi(G)$, then $G$ possesses a $\pi$-subgroup  not isomorphic to a subgroup of~$H$. In
particular, $G\not\in
D_\pi$.
\end{itemize}
\end{Lemma}

\noindent{\slshape Proof.}\ \
As in the previous lemma, if we show that a $\pi$-Hall subgroup $H$ normalizes a maximal torus, then the claim follows
from
\cite[Lemma~6.6]{RevVdoContemp}. Assume that
$H$ is a $\pi$-Hall subgroup of $G$. Then $H$ is included in a maximal subgroup $M$ of odd index of $G$. In view of
\cite{liesax} and
arguing as in the proof of Lemma \ref{HallG2}, we may assume that one of the
following statements holds.
\begin{itemize}
\item[(1)] $\vert M\vert=q^7(q^2-1)^7\cdot 2^2\cdot 3\cdot 7$,
\item[(2)] $\vert M\vert=1/2\cdot q^{31}(q^2-1)^2(q^4-1)(q^6-1)^2(q^8-1)(q^{10}-1)$,
\item[(3)] $M=N_G(2^2\arbitraryext (\mathrm{PSL}_2(q)^3\times\mathrm{P}\Omega_8^+(q)))$, $\vert M\vert=
q^{15}(q^2-1)^4(q^4-1)^2(q^6-1)\cdot
3$,
\item[(4)] $M= N(G,T)$, where $T$ is a maximal torus defined in the lemma, and $N(G,T)$ is its algebraic normalizer
(defined
in~\cite{RevVdoContemp}).
\end{itemize}
If one of statements (1)--(3) holds, then $\vert G:M\vert_\ti\ge3$, and if statement (4) holds, then $H$ is
in~${N(G,T)=T\arbitraryext
W(E_7)}$.
\qed

\begin{Lemma}\label{HallE8}
Suppose that $G\simeq E_8(q)$, $2,3\in\pi$, $p\not\in\pi$, $\varepsilon=\varepsilon(q)$ and  $T$ is a maximal torus  of
order $(q-\varepsilon)^8$ of $G$.  Then the
following  statements hold{\em:}
\begin{itemize}
\item[{\em (A)}] If $G\in E_\pi$ and $H\in \Hall_\pi(G)$, then $\pi\cap \pi(G)\subseteq \pi(q-\varepsilon)$,
$5,7\in\pi$,
and $H$, up to
conjugation, is included in $T\arbitraryext W(E_8)$.
\item[{\em (B)}] Conversely, if  $\pi\cap \pi(G)\subseteq \pi(q-\varepsilon)$ and $5,7\in\pi$, then $G\in E_\pi$.
\item[{\em (C)}] If $G\in E_\pi$, then $k_\pi(G)=1$, i.e., $G\in C_\pi$.
\item[{\em (D)}] If $H\in\Hall_\pi(G)$, then $G$ possesses a $\pi$-subgroup  not isomorphic to a subgroup of~$H$. In
particular, $G\not\in
D_\pi$.
\end{itemize}
\end{Lemma}

\noindent{\slshape Proof.}\ \
If we show that a $\pi$-Hall subgroup $H$ normalizes a maximal torus, then the lemma will follow from
\cite[Lemma~6.7]{RevVdoContemp}. In
view
of
\cite{liesax} and arguing as in the proof of Lemma \ref{HallG2}, we may assume that  a maximal subgroup $M$ of $G$ that
includes $H$
satisfies to one of the following statements.
\begin{enumerate}
\item[(1)] $\vert M\vert=q^8(q^2-1)^8\cdot 2^6\cdot 3\cdot 7$,
\item[(2)] $\vert M\vert=q^{56}(q^2-1)(q^4-1)(q^6-1)^2(q^8-1)^2(q^{10}-1)(q^{12}-1)(q^{14}-1)$,
\item[(3)] $M=N_G(2^2\arbitraryext (\mathrm{P}\Omega_8^+(q)\times\mathrm{P}\Omega_8^+(q)))$,  $\vert
M\vert=q^{24}(q^2-1)^2(q^4-1)^4(q^6-1)^2\cdot 2^2\cdot 3$,
\item[(4)] $M= N(G,T)$, where $T$ is a maximal torus defined in the lemma, and $N(G,T)$ is its algebraic normalizer.
\end{enumerate}
If one of statements (1)--(3) holds, then $\vert G:M\vert_\ti\ge3$, and if statement (4) holds, then $H$ is
in~${N(G,T)=T\arbitraryext
W(E_8)}$.
\qed

\begin{Lemma}\label{Hall3D4}
Suppose that $G\simeq {}^3D_4(q)$, $2,3\in\pi$, $p\not\in\pi$, $\varepsilon=\varepsilon(q)$ and  $T$ is a maximal torus 
of order $(q-\varepsilon)(q^3-\varepsilon)$  of
$G$. Then the following  statements hold{\em:}
\begin{itemize}
\item[{\em (A)}] If $G\in E_\pi$ and $H\in \Hall_\pi(G)$, then $\pi\cap \pi(G)\subseteq
\pi(q-\varepsilon)$, and $H$, up to conjugation,  is
included in $T\arbitraryext W(G_2)$.
\item[{\em (B)}] Conversely, if  $\pi\cap \pi(G)\subseteq \pi(q-\varepsilon)$, then $G\in E_\pi$.
\item[{\em (C)}] If $G\in E_\pi$, then $k_\pi(G)=1$, i.e., $G\in C_\pi$.
\item[{\em (D)}] If $H\in\Hall_\pi(G)$, then $G$ possesses a $\pi$-subgroup  not isomorphic to a subgroup of~$H$. In
particular, $G\not\in
D_\pi$.
\end{itemize}
\end{Lemma}

\noindent{\slshape Proof.}\ \
If we show that a $\pi$-Hall subgroup $H$ normalizes a maximal torus, then the claim follows from
\cite[Lemma~6.10]{RevVdoContemp}. In view
of
\cite{liesax} and arguing as in the proof of Lemma \ref{HallG2}, we may assume that  a maximal subgroup $M$ of $G$ that
includes $H$,
satisfies to one of the following statements.
\begin{enumerate}
\item[(1)] $M\simeq G_2(q)$, $\vert M\vert=q^{6}(q^2-1)(q^6-1)$,
\item[(2)]  $M=N_G(2\arbitraryext (\mathrm{PSL}_2(q)\times \mathrm{PSL}_2(q)))$, $\vert
M\vert=q^{4}(q^2-1)(q^6-1)$,
\item[(3)] $M=N_G(\mathrm{SL}_3^\varepsilon(q))$, $\vert M\vert=q^{3}(q^3-\varepsilon)(q^2-1)(q^2+\varepsilon
q +1)/2$.
\end{enumerate}
If one of statements (1), (2) holds, then $\vert G:M\vert_\ti\ge3$. Suppose statement (3) holds. Then the index $\vert
G:M\vert$  is not
divisible by $3$ if and only if $q\equiv \varepsilon\pmod3$. By Lemma \ref{HallSubgroupsOfLinearAndUnitaryGroups}, we
obtain that either a
$\pi$-Hall subgroup normalizes a maximal torus, or is included in  $\mathrm{GL}_2^\varepsilon(q)$. In the second case we
have $\vert M:
\mathrm{GL}_2^\varepsilon(q)\vert_\ti>1$.
\qed

\section{Proof of the Theorem 1.1}

We first prove Theorem \ref{IndClassNumb} in the particular case $2,3\in\pi$ and~${G=S}$.

\begin{Lemma}\label{ClassNumb} Let $\pi$ be a set of primes such that $2,3\in\pi$ and let $S$ be a finite simple group.
Then one of the
following statements holds:
\begin{itemize}
 \item[{\em (A)}]  $k_\pi(S)\in\{0,1,2,3,4\}$;
 \item[{\em (B)}]  $k_\pi(S)=9$, $S\simeq {\rm PSp}_{2n}(q)$, $q$ is a power of a prime $p\not\in\pi$ and
\begin{itemize}
\item[{\em (a)}]  either $n\in\{5,7\}$, $\pi\cap\pi(S)=\{2,3\}\subseteq\pi(q^2-1)$ and $(q^2-1)_\pi=48$;
\item[{\em (b)}] or $n=7$, $\pi\cap\pi(S)=\{2,3,5\}\subseteq\pi(q^2-1)$ and $(q^2-1)_\pi=120$.
\end{itemize}
\end{itemize}
\end{Lemma}

\noindent{\slshape Proof.}\ \  We may assume that $S$ is nonabelian and $S\in E_\pi$. Consider all nonabelian finite
simple $E_\pi$-groups with
$2,3\in\pi$.

\begin{center}
\begin{table}\caption{Proper $\pi$-Hall subgroups in sporadic groups, $2,3\in\pi$}\label{tb0}
\begin{center}
\begin{tabular}{|l|l|r|} \hline  $S$ & $\pi$ &
Structure $H$ \\ \hline $M_{11}$  & $\{2,3\}$ & $3^2\splitext Q_8\arbitraryext 2$ \\ &
$\{2,3,5\}$       &  $\Alt_6\arbitraryext 2$\\ \hline $M_{22}$  & $\{2,3,5\}$ &
$2^4\splitext\Alt_6$\\ \hline $M_{23}$  & $\{2,3\}$         & $2^4\splitext(3\times
\Alt_4)\splitext2$\\ & $\{2,3,5\}$ &  $2^4\splitext\Alt_6$\\ & $\{2,3,5\}$       &
$2^4\splitext(3\times \Alt_5)\splitext2$\\ & $\{2,3,5,7\}$     &  ${\rm
          L}_3(4)\splitext 2_2$\\ & $\{2,3,5,7\}$     &  $2^4\splitext \Alt_7$\\ &
          $\{2,3,5,7,11\}$  & $M_{22}$\\ \hline $M_{24}$  & $\{2,3,5\}$
          & $2^6\splitext 3\nonsplitext \Sym_6$\\ \hline $J_1$     & $\{2,3\}$         &
$2\times \Alt_4$\\ & $\{2,7\}$         &  $2^3\splitext7$\\ & $\{2,3,5\}$ &
          $2\times \Alt_5$\\ & $\{2,3,7\}$       &  $2^3\splitext7\splitext3$\\ \hline
$J_4$     & $\{2,3,5\}$       &  $2^{11}\splitext(2^6\splitext3\nonsplitext\Sym_6)$\\
\hline
\end{tabular}
\end{center}
\end{table}
\end{center}

If $S\simeq {\Alt}_n$, then the claim follows from Lemma~\ref{HallSymmetric}. In sporadic groups all proper $\pi$-Hall
subgroups with
$2,3\in\pi$ are found in \cite[Theorem~4.1]{rev2} and are given in Table~\ref{tb0}. By using \cite{ATLAS}, it is easy to
check
that~${k_\pi(S)\le2}$ in this case.

If $S$ is a finite group of Lie type with the base field of characteristic $p\not\in\pi$, then the claim follows from
Lemmas \ref{HallSubgroupsOfLinearAndUnitaryGroups}(C), \ref{HallSubgroupsOfSymplecticGroups}(C),
\ref{HallSubgroupsOfOrthogonalGroupsOfEvenDimension}(C),  \ref{HallG2}(C), \ref{HallF4}(C), \ref{HallE6}(C),
\ref{HallE7}(C),
\ref{HallE8}(C), \ref{Hall3D4}(C), and the proof of Lemma~\ref{HallSubgroupsOfSymplecticGroups}.

Now assume that $S$ is a finite group of Lie type with the base field $\F_q$ of characteristic $p\in\pi$. Then
\cite[Theorem~3.2]{GrossConjectureHall}, \cite[Theorem~3.1]{GrossNotDivisibleBy3}, and~\cite[Theorem~1.2]{RevinSIBAM}
imply that one of
the following cases 1--4 holds:

%\begin{enumerate}
%\item
Case 1. $\pi(S)\subseteq\pi$. Clearly~${k_\pi(S)=1}$.

%\item
Case 2. $\pi\cap\pi(S)\subseteq \pi(q-1)\cup\{p\}$ and each $\pi$-Hall subgroup $H$ of $S$ is included in a Borel
subgroup of $S$. In this
case the solvability of Borel subgroups and the fact that all Borel subgroups are conjugate imply that~${k_\pi(S)=1}$.

Case 3.
%\item
$S= D_n^\eta(q)$, $p=2$, $\pi\cap\pi(S)=\pi(D_{n-1}^\eta(q))$ and each $\pi$-Hall subgroup $H$ of $S$ is a parabolic
subgroup with the Levi
factor isomorphic to $D_{n-1}^\eta(q)$. Since all parabolic subgroups of this type are conjugate, we
have~${k_\pi(S)=1}$.

Case 4.
%\item
$S=A_{n-1}(q)\simeq {\rm PSL}(V)$, where $V$ is a vector space of dimension $n$ over $\F_q$;
$\pi\cap\pi(S)=\pi(H)$, where $H$ is the image in  ${\rm PSL}(V)$ of the stabilizer in~$\SL(V)$  of a series
$$0=V_0<V_1<\ldots<V_s=V$$ such that
 $\dim V_i/V_{i-1}=n_i$,  $n=\sum n_i$ and one of the following statements holds:
\begin{itemize}
 \item[(4.1)] $n$ is an odd prime, $s=2$, $\{n_1,n_2\}=\{1,n-1\}$;
\item[(4.2)]
$n=4$, $s=2$, $n_1=n_2=2$;
\item[(4.3)] $n=5$, $s=2$ $\{n_1,n_2\}=\{2,3\}$;
\item[(4.4)] $n=5$, $s=3$, $\{n_1,n_2,n_3\}=\{1,2\}$;
 \item[(4.5)] $n=7$, $s=2$, $\{n_1,n_2\}=\{3,4\}$;
\item[(4.6)] $n=8$, $s=2$, $n_1=n_2=4$;
\item[(4.7)] $n=11$, $s=2$, $\{n_1,n_2\}=\{5,6\}$.
\end{itemize}

Every $\pi$-Hall subgroup of $S$ is equal to one of such subgroups~$H$. It is easy to see that in distinct cases either
the dimensions, or
the sets
$\pi(H)$ are distinct, so distinct statements cannot be fulfilled simultaneously. The number of classes of conjugate
$\pi$-Hall subgroups
in each statement can be calculated directly. They are as given below.
\begin{itemize}
 \item[(4.1)]
$k_\pi(S)=2$;
\item[(4.2)]
 $k_\pi(S)=1$;
\item[(4.3)] $k_\pi(S)=2$;
\item[(4.4)]  $k_\pi(S)=3$;
 \item[(4.5)]  $k_\pi(S)=2$;
\item[(4.6)]  $k_\pi(S)=1$;
\item[(4.7)]  $k_\pi(S)=2$.
\end{itemize}
Thus the claim follows.
\qed

\begin{Lemma}\label{InducedSubgroupsInSymplectic}
Suppose $S\simeq \PSp_{2n}(q)\in E_\pi$, where $2,3\in\pi$ and $p\not\in\pi$. Assume also that $k_\pi(S)=9$. Choose
$G\in E_\pi$ so that $S\leq G\leq \Aut(S)$. Then $k^G_\pi(S)\in\{1,9\}$.
\end{Lemma}

\noindent{\slshape Proof.}\ \
By Lemma \ref{HallSubgroupsOfSymplecticGroups}, each $\pi$-Hall subgroup $H$ of $S$ is included in a subgroup
$$M=(\Sp_2(q)\wr
\Sym_n)/Z(\Sp_{2n}(q))=L\arbitraryext \Sym_n,$$ where $L=L_1\circ\ldots\circ L_n\simeq \Sp_2(q)\circ\ldots\circ
\Sp_2(q)$. Moreover all such
subgroups $M$ are conjugate in $S$, while $H,K\in \Hall_\pi(M)$ are conjugate in $S$ if and only if $H,K$ are conjugate
in $M$.
Since $p$ is odd, $S$ does not possesses graph automorphisms. So $\Out(S)$ is abelian and it is generated by a diagonal
automorphism (of
order $2$) and a field automorphism. Clearly we may choose representatives $\delta$ of the diagonal automorphism and
$\varphi$ of the field
automorphism
so that both $\delta$ and $\varphi$ normalize each factor $L_i$, and $\delta$ is a $2$-element. Moreover, preserving
such a choice, we may
also assume that $\delta$ induces a
diagonal automorphism on each factor $L_i$,
while $\varphi$ induces a field automorphism on each~$L_i$.

Now $H=(H_1\circ\ldots\circ H_n)\arbitraryext \ov{H}$, where $H_i\in\Hall_\pi(L_i)$ and $\ov{H}\in\Hall_\pi(\Sym_n)$.
Since
$k_\pi(S)=k_\pi(M)=9$, Lemma \ref{Numberkpi(G)Inextension} and the proof of Lemma \ref{HallSubgroupsOfSymplecticGroups}
imply that
$\ov{H}$, acting by conjugation on $\{L_1,\ldots, L_n\}$, has $2$ orbits and $k_\pi(L_i)=3$ for each $i$. By Lemma
\ref{SLU23dim2}, we
obtain that $L_i$ possesses a class $\mathcal{K}_i$ of conjugate $\pi$-Hall subgroups satisfying Lemma
\ref{SLU23dim2}(a), and classes $\mathcal{L}_i'$ and $\mathcal{L}_i''$  consisting of $\pi$-Hall subgroups
satisfying either (c) or (d) of Lemma \ref{SLU23dim2}. If we choose $K_i\in \mathcal{K}_i$ and consider
$K=K_1\circ\ldots\circ K_n$, then
$K^S=K^G=K^{\Aut(S)}$, hence $k_\pi(G)\ge 1$.

If $k_\pi^G(S)\ge 2$, then there exists $F\in\Hall_\pi(N_G(M))\subseteq \Hall_\pi(G)$ such that either $F\cap
L_i\in\mathcal{L}_i'$, or
$F\cap L_i\in\mathcal{L}_i''$
for some $i$. Since $\delta$ is a $2$-element and normalizes~$M$, the Sylow theorem implies that we may assume
$\delta\in F$. Moreover, the
induced action of both $\delta$ and $\varphi$ on $M/L$ is trivial. Since $FL/L$ is a $\pi$-Hall subgroup of $N_G(M)/L$,
we obtain that $F$
contains an element $\psi$ such that $\psi$ normalizes each $L_j$, induces a field automorphism on each $L_j$, and
$F=\langle\delta,\psi,H\rangle$, where $H=F\cap
M\in \Hall_\pi(M)\subseteq \Hall_\pi(S)$. Consider $K=\langle \delta,\psi,F\cap L\rangle$. Then either $K\cap
L_i\in\mathcal{L}_i'$, or
$K\cap
L_i\in\mathcal{L}_i''$, hence both classes $\mathcal{L}_i'$ and $\mathcal{L}_i''$ are invariant under $K$.  By
construction, both $\delta$
and $\psi$ induce automorphisms of the same type on each $L_j$, so
we obtain that $\mathcal{L}_j'$ and $\mathcal{L}_j''$ are invariant under $K$ for $j=1,\ldots,n$. In particular, each
class of conjugate $\pi$-Hall
subgroups of $L$ is invariant under $\delta$ and $\psi$. Lemma \ref{simpleepi}(e) implies that each $\pi$-Hall subgroup
of
$L$ is embedded into a $\pi$-Hall subgroup of $\langle\delta,\psi,L\rangle$, hence of $\langle
\delta,\varphi, L\rangle$. Since both $\delta$ and $\psi$ induce a trivial automorphism on $M/L\simeq \Sym_n$, Lemma
\ref{simpleepi}(d), (e)
implies that a $\pi$-Hall subgroup $K$ of $L$ is embedded into a $\pi$-Hall subgroup of $\langle \delta,\psi,M\rangle$
if and only if
$K$ is embedded into a $\pi$-Hall subgroup of $M$. Thus $\Hall_\pi(M)=\Hall_\pi^{N_G(M)}(M)$. Since all subgroups $M$
are conjugate in
$S$ and $H_1,H_2\in\Hall_\pi(M)$ are conjugate in $S$ if and only if they are conjugate in $M$, we obtain that
$9=k_\pi(S)=k_\pi(M)=k_\pi^{N_G(M)}(M)\le k_\pi^G(S)\le 9$.
\qed

Now we are able to prove Theorem~\ref{IndClassNumb}.

Without lost of generality we may assume that $G\in E_\pi$, and by Lemma \ref{simpleepi}(a) we have $S\in E_\pi$.

If $2\not\in\pi$, then  $\pi$-Hall subgroups of $S$ are conjugate by \cite[Theorem~B]{GrossConjugacy}, i.e.,
$k_\pi(S)=1$, whence~${k_\pi^G(S)=1}$.

If $2\in\pi$ and  $3\not\in\pi$, then \cite[Lemma~5.1 and Theorem~5.2]{RevVdoContemp} implies that $S\in C_\pi$, unless $S= {}^2G_2(q)$.
Therefore $k_\pi^G(S)=k_\pi(S)=1$ if $S\not={}^2G_2(q)$. If  $S= {}^2G_2(q)$, $q=3^{2n+1}$ and either
$\pi\cap\pi(S)\ne\{2,7\}$, or $n\equiv 3\pmod 7$, then again $k_\pi^G(S)=k_\pi(S)=1$. Suppose that $S={}^2G_2(q)$, $q=3^{2n+1}$,
$\pi\cap\pi(S)=\{2,7\}$, $n\not\equiv 3\pmod 7$ and $S\in E_\pi$. By \cite[Lemma~5.1]{RevVdoContemp}, we have that $7\in\pi(q+1)$ and every
$\pi$-Hall subgroup of  $S$ is either included in the normalizer of a maximal torus of order $q+1$, or is a  Frobenius group of order $56$.
In the first case a $\pi$-Hall subgroup has a Sylow tower of complexion $2\prec 7$, and in the second case a $\pi$-Hall subgroup has a
Sylow tower of complexion $7\prec 2$. Thus $k_\pi(S)=2$ and, since
$k_\pi^G(S)\leqslant k_\pi(S)$, it follows that $k_\pi^G(S)\in\{1,2\}$.\footnote{Notice that Lemma \ref{simpleepi}(e) implies
$k_\pi^G(S)=2$, since
classes of nonisomorphic subgroups cannot be interchanged by an automorphism.}

Since the
inequality $k_\pi(S)\geqslant k_\pi^G(S)$ is evident, we may assume that $k_\pi(S)> 4$. Lemma  \ref{ClassNumb} implies that
$k_\pi(S)=9$ and
$S\simeq \PSp_{2n}(q)$. Now Theorem \ref{IndClassNumb} follows from Lemma~\ref{InducedSubgroupsInSymplectic}.

\section{Proof of Corollary~1.3}

Lemma~\ref{simpleepi}(a) implies that~$A$ satisfies~$E_\pi$. Each $\pi$-subgroup~$K$ of~$A$ is included in a $\pi$-Hall subgroup
$H$ of~$G$. Hence $K$ is included in a $\pi$-Hall subgroup $H\cap A$ of~$A$. So it remains to show that $A\in C_\pi$. We proceed by
induction on~$|A|$.

Suppose there exists a normal subgroup~$N$ of~$G$ such that $1<N<A$.  Then Lemma~\ref{simpleepi}(c) implies that $G/N$  satisfies~$D_\pi$
and $A/N$ is normal in~$G/A$. By induction both $A/N$ and $N$ satisfy $C_\pi$, whence $A\in C_\pi$ by Lemma~\ref{simpleepi}(f). Thus we may
assume that $A$ is a minimal normal subgroup of~$G$, so
$$A=S_1\times\dots\times S_n,$$ where $S_1,\dots, S_n$ are conjugate simple subgroups of~$G$. In view of the Hall theorem we may also assume
that $S_i$ is nonabelian for~${i=1,\ldots,n}$.

By Theorem~\ref{IndClassNumb} it follows that  $k=k_\pi(S_i)$ is a $\pi$-number. Clearly each $\pi$-Hall subgroup of~$A$ is equal to
$$\langle K_1,\dots, K_n\rangle=K_1\times\dots\times K_n,$$ where
$K_i$ is a $\pi$-Hall subgroup of~$S_i$ for $i=1,\dots, n$, and vice versa. So $k_\pi(A)=k^n$ is a $\pi$-number.

The group~$G$ acts by conjugation on the set~$\Omega$ of classes of conjugate $\pi$-Hall subgroups of $A$. Moreover this action is
transitive. Indeed, since~$G$ satisfies~$D_\pi$, each $\pi$-Hall subgroup of~$A$ has the form~$H\cap A$, where $H$ is a $\pi$-Hall subgroup
of $G$. Since $G$ acts by conjugation transitively on $\Hall_\pi(G)$, it follows that $G$ acts by conjugation transitively on
$\Hall_\pi(A)$, hence on
$\Omega$. We fix $H\in\Hall_\pi(G)$, let $\Delta\in\Omega$ be the class of conjugate $\pi$-Hall subgroups of $A$, containing~${H\cap A}$.
Since $H$ normalizes~$H\cap A$, it follows that $H$ is included in the stabilizer~$G_\Delta$ of this class. So
$k_\pi(A)=|\Omega|=|G:G_\Delta|$ divides $|G:H|.$ Hence $k_\pi(A)$ is a~$\pi'$-number.
Thus $k_\pi(A)$ is a $\pi$- and a $\pi'$- number, so $k_\pi(A)=1$ and $A\in C_\pi$.

\end{document}